\documentclass[11pt,a4]{article}
\usepackage[T2A]{fontenc}
\usepackage[cp866]{inputenc}
\usepackage[english]{babel}
\usepackage[mathscr]{eucal}
\usepackage{amssymb}
\usepackage{amsmath}
\usepackage{fancyhea}
\renewcommand{\leq}{\leqslant}
\renewcommand{\geq}{\geqslant}

\renewcommand{\C}{{\mathbb C}}

\newcommand{\R}{{\mathbb R}}

\renewcommand{\k}{\rule{0.7em}{0.7em}}
\begin{document}
\sloppy
\fussy
\headrulewidth = 2pt
\pagestyle{fancy}
\lhead[\scriptsize V.N.Gorbuzov,  A.F. Pranevich]{\scriptsize V.N.Gorbuzov,  A.F. Pranevich}
\rhead[\it \scriptsize First integrals of linear differential systems]{\it \scriptsize First integrals of linear differential systems}
\headrulewidth=0.25pt

{\normalsize 
\mbox{}
\\[-5.5ex]
\centerline{
{\large\bf 
FIRST  INTEGRALS  OF  LINEAR  DIFFERENTIAL SYSTEMS}
}
\\[2.5ex]
\centerline{
\bf 
V.N.Gorbuzov$\!{}^{*},$   A.F. Pranevich$\!{}^{**}$
}
\\[2.5ex]
\centerline{
\it 
$\!{}^{*}\!\!$Department of Mathematics and Computer Science, 
Yanka Kupala Grodno State University,
}
\\[1ex]
\centerline{
\it 
Ozeshko 22, Grodno, 230023, Belarus
}
\\[1ex]
\centerline{
E-mail: gorbuzov@grsu.by
}
\\[2.5ex]
\centerline{
\it 
$\!{}^{**}\!\!$Department of Economics and Management, Yanka Kupala Grodno State University,
}
\\[1ex]
\centerline{
\it 
 Ozeshko 22, Grodno, 230023, Belarus
}
\\[1ex]
\centerline{
E-mail: pronevich@tut.by
}
\\[7ex]
\centerline{
{\large\bf Abstract}}
\\[0.75ex]
\indent
We investigate the problem of the existence of first integrals
for multidimensional and ordinary linear differential systems with constant co\-ef\-fi\-ci\-ents.
The spectral method of the first integrals basis construction
for these systems of linear differential equations is developed.
\\[2ex]
\\[0.75ex]
\\[5.5ex]
{\large\bf Contents} 
\\[1ex]
{\bf  
1. ${\mathbb R}\!$-differentiable integrals of  ${\mathbb R}\!$-linear systems in total differentials}                 \hfill\ 2
\\[0.5ex]
\mbox{}\hspace{1em}
1.1. ${\mathbb R}\!$-linear homogeneous systems in total differentials\footnote[1]{
This Subsection has been published in {\it Dokl. Akad. Nauk Belarusi} 48 (2004), No. 1, 49-52.
} 
                                                                                                                                                                         \dotfill\ 2
\\[0.5ex]
\mbox{}\hspace{2.75em}
1.1.1. ${\mathbb R}\!$-linear partial integral                                                                                                \dotfill\ 2
\\[0.5ex]
\mbox{}\hspace{2.75em}
1.1.2. Autonomous ${\mathbb R}\!$-differentiable  first integrals                                                             \dotfill\ 3
\\[0.5ex]
\mbox{}\hspace{2.75em}
1.1.3. Nonautonomous ${\mathbb R}\!$-differentiable  first integrals                                                       \dotfill\ 7
\\[0.5ex]
\mbox{}\hspace{1em}
1.2. $\R\!$-linear nonhomogeneous systems in total differentials                                                             \dotfill\ 8
\\[1ex]
{\bf  
2. First integrals of  linear real systems in total differentials}                                                                     \hfill\ 9
\\[0.5ex]
\mbox{}\hspace{1em}
2.1. Linear real homogeneous systems in total differentials\footnote[2]{
The definitive version of this Subsection has been published in 
{\it Differential equations and control processes}, 2001, No.3, 17-45 (http://www.neva.ru/journal);
Deponent VINITI of {\it Differentsial Uravneniya} ({\it Differential Equations}), 02.10.2002, No. 1667-B2002;
{\it Mathematical research}, Vol. 10 (2003), 143-152.
}                                                                                                                                                    \dotfill\ 9
\\[0.5ex]
\mbox{}\hspace{2.75em}
2.1.1. Partial integrals                                                                                                                                   \dotfill\ 9
\\[0.5ex]
\mbox{}\hspace{2.75em}
2.1.2. Autonomous first integrals                                                                                                                \dotfill\ 10
\\[0.5ex]
\mbox{}\hspace{2.75em}
2.1.3. Nonautonomous first integrals                                                                                                           \dotfill\ 20
\\[0.5ex]
\mbox{}\hspace{1em}
2.2. Linear real nonhomogeneous systems in total differentials                                                                  \dotfill\ 22
\\[1ex]
{\bf  
3.  First integrals of  linear systems of ordinary differential equations}\footnote[3]{
Section 3 has been published in {\it Vestnik of the Yanka Kupala Grodno State Univ.}, 
2002, Ser. 2, No. 2(11), 23-29; 2003, Ser. 2, No. 2(22), 50-60;
2008, Ser. 2, No. 2(68), 5-10.}                                                                                                            \hfill\ 26
\\[0.5ex]
\mbox{}\hspace{1em}
3.1. Linear homogeneous systems of ordinary differential equations                                                       \dotfill\ 26
\\[0.5ex]
\mbox{}\hspace{2.75em}
3.1.1. Autonomous first integrals                                                                                                                 \dotfill\ 26
\\[0.5ex]
\mbox{}\hspace{2.75em}
3.1.2. Nonautonomous first integrals                                                                                                          \dotfill\ 30
\\[0.5ex]
\mbox{}\hspace{1em}
3.2. Linear nonhomogeneous systems of ordinary differential equations                                                 \dotfill\ 32
\\[0.75ex]
{\bf  
References}                                                                                                                                          \dotfill\ 34

\thispagestyle{empty}

\newpage

\mbox{}
\\[1ex]
\centerline{
\large\bf  
1. ${\mathbb R}\!$-differentiable integrals of  $\R\!$-linear systems in total differentials
}
\\[1.5ex]
\indent
{\bf  
1.1. ${\mathbb R}\!$-linear homogeneous systems in total differentials
}
\\[1ex]
\indent
We consider the system of equations in total differentials
\\[1.75ex]
\mbox{}\hfill                           
$
dw = X_1(w)dz+X_2(w)d\,\overline{z}\,,
$
\hfill (1.1)
\\[1.75ex]
where 
\vspace{0.5ex}
$w\!=\mbox{colon}(w_1,\ldots,w_n)\!\in\! {\mathbb C}^n, \,
z\!=\mbox{colon}(z_1,\ldots,z_m)\!\in\! {\mathbb C}^m; \
dw\!=\mbox{colon}(dw_1,\ldots,dw_n),\!\!$ 
$dz=\mbox{colon}(dz_1,\ldots,dz_m),$ and
\vspace{0.5ex}
$d\,\overline{z}=\mbox{colon}(d\,\overline{z}_1,\ldots,d\,\overline{z}_m)$ 
are vector columns; $\overline{z}_j$ is complex conjugate of $z_j;$
\vspace{0.5ex}
the entries of the matrices
$X_1(w)=\|X_{\tau j}(w)\|$ and $X_2(w)=\|X_{\tau, m+j}(w)\|,$
$\tau=1,\ldots, n,\ j=1,\ldots, m$
are the $\R\!$-linear functions [1, p. 21]
\\[1.75ex]
\mbox{}\hfill
$
\displaystyle
X_{{}_{\scriptstyle \tau k}}\colon w\to
\sum\limits_{\xi=1}^{n}
\bigl( a_{{}_{\scriptstyle \tau k\xi}} w_{{}_{\scriptstyle \xi}}+
a_{{}_{\scriptstyle \tau k,n+\xi}}
\overline{w}_{{}_{\scriptstyle \xi}}\bigr)
\ \ \text{for all}\  w\in {\mathbb C}^n,
\ \ \ k=1,\ldots, 2m,
\ \ \theta=1,\ldots, n,
\hfill
$
\\[1.75ex]
with constant coefficients $a_{{}_{\scriptstyle \tau k\varrho}}\in {\mathbb C},\
\varrho=1,\ldots, 2n,\ k=1,\ldots,2m, \ \tau=1,\ldots,n.$
We assume that the linear differential operators 
\\[1.5ex]
\mbox{}\hfill                                       
$
\displaystyle
{\frak x}_j(w) = \sum\limits_{\xi=1}^{n}
\bigl(X_{\xi j}(w)\partial_{{}_{\scriptstyle w_\xi}} +
\overline{X}_{\xi,m+j}(w)\partial_{{}_{\scriptstyle \overline{w}_\xi}}
\bigr)
\ \ \ \text{for all}\  w\in {\mathbb C}^n,
\ \ j=1,\ldots, m,
$
\hfill (1.2)
\\[1ex]
and
\\[1ex]
\mbox{}\hfill                                       
$
\displaystyle
{\frak x}_{m+j}(w) = \sum\limits_{\xi=1}^{n}
\bigl(X_{\xi,m+j}(w)\partial_{{}_{\scriptstyle w_\xi}} +
\overline{X}_{\xi j}(w)\partial_{{}_{\scriptstyle \overline{w}_\xi}}
\bigr)
\ \ \ \text{for all}\  w\in {\mathbb C}^n,
\ \ j=1,\ldots, m,
$
\hfill (1.3)
\\[1.5ex]
induced by this system are related by the Frobenius conditions [2; 3]. 
These conditions are represented via Poisson brackets as the system of identities
\\[1.5ex]
\mbox{}\hfill                                       
$
\bigl[ {\frak x}_{k}(w), {\frak x}_{l}(w)\bigr]  = {\frak O}
\ \ \ \text{for all}\  w\in {\mathbb C}^n,
\ \ k=1,\ldots, 2m,\ l=1,\ldots, 2m,
$
\hfill (1.4)
\\[1.5ex]
i.e., system (1.1) is completely solvable [4; 5, pp. 15 -- 25]. 

A general integral of the completely solvable system in total differentials (1.1)
is $n$ functionally independent ${\mathbb R}\!$-differentiable first integrals of (1.1).
The completely solvable dif\-fe\-ren\-t\-ial system (1.1) has also 
$n-m$ autonomous ${\mathbb R}\!$-differentiable first integrals (see [6]).

In this paper we study Darboux's problem of finding first in\-teg\-rals 
in case that partial integrals are known [7]. Using method of partial integrals [5; 8; 9], we obtain the 
spectral method for building first integrals of linear differential systems [5, pp. 239 -- 272; 10 -- 17].
 
The problems related to those studied in the present paper were intensively investigated by many 
people. For example, see for systems of ordinary differential equations [7; 12; 16 -- 54],
for  partial differential systems [5; 8; 10; 15;  55 -- 59],
for systems of equations in total differentials 
\linebreak
$[2 - 6; 8; 9; 11; 13; 14; 16; 32; 60 - 63],$
and this list is very far from being complete.
\\[1.5ex]
\indent                       
{\bf 1.1.1. ${\mathbb R}\!$-linear partial integral}.
The $\R\!$-linear function
\\[1.5ex]
\mbox{}\hfill
$
\displaystyle
p\colon w\to 
\sum\limits_{\xi=1}^{n}
\bigl(
b_{{}_{\scriptstyle \xi}} w_{{}_{\scriptstyle \xi}} +
b_{{}_{\scriptstyle n+\xi}}
\overline{w}_{{}_{\scriptstyle \xi}}\bigr)
\ \ \text{for all}\  \,w\in {\mathbb C}^n
\quad 
(b_{\varrho}\in {\mathbb C},\ \varrho=1,\ldots, 2n)
\hfill
$
\\[1.5ex]
is a {\it partial integral} of the system in total differentials (1.1) iff
\\[1.75ex]
\mbox{}\hfill
$
{\frak x}_{{}_{\scriptstyle k}} p(w)= p(w)\lambda^k
\ \ \text{for all}\  w\in {\mathbb C}^n,
\quad
\lambda^k\in {\mathbb C}, \ k=1,\ldots, 2m.
\hfill
$
\\[1.75ex]
This system of identities is equivalent to the linear homogeneous system 
\\[1.75ex]
\mbox{}\hfill                                      
$
\bigl(A_{k} - \lambda^{k} E\bigr)\,b= 0,
\ \ \ k=1,\ldots, 2m,
$
\hfill (1.5)
\\[1.75ex]
where $A_{j}= \|a_{1j}\ldots a_{nj}\,
\overline{a}_{1,m+j}\ldots \overline{a}_{n,m+j}\|$ and
\vspace{1ex}
$A_{m+j}= \|a_{1,m+j}\ldots a_{n,m+j}
\overline{a}_{1j}\ldots\overline{a}_{nj}\|\!$
are the $\!\!2n\!\times\! 2n\!$-matrices with   
\vspace{0.75ex}
$\!a_{\tau k}\!\!=\!\mbox{colon}(\!a_{\tau k1},\ldots,a_{\tau k,2n}\!),
\overline{a}_{\tau k}\!\!=\!
\mbox{colon}(\overline{a}_{\tau k,n+1},\ldots,\overline{a}_{\tau k,2n},
\overline{a}_{\tau k1},\ldots,\overline{a}_{\tau kn}\!)\!,\!\!$
$\tau =1,\ldots,n,\ k=1,\ldots, 2m, \ j=1,\ldots,m,\  
E$ is the $2n\times 2n$ identity matrix, and
\linebreak
$b=\mbox{colon}(b_1,\ldots,b_{2n})$ is a vector column. 
\vspace{0.5ex}

The Frobenius conditions (1.4) for system (1.1) are equivalent [60, p. 73]
\\[1.5ex]
\mbox{}\hfill
$
A_{k}A_{l} = A_{l}A_{k},
\ \ \ 
k=1,\ldots,2m,
\ \ l=1,\ldots, 2m.
\hfill
$
\\[1.5ex]
\indent
Then there exists a relation [64, pp. 193 -- 194; 65] between eigenvectors and 
eigenvalues of the matrices $A_{k},\ k=1,\ldots, 2m.$
\vspace{0.5ex}

{\bf  Lemma 1.1}.
\vspace{0.5ex}
{\it
Let $\nu\in {\mathbb C}^{2n}$ be a common eigenvector of the matrices 
$A_{k},\ k=1,\ldots, 2m.$ Then the ${\mathbb R}\!$-linear function
\\[1.25ex]
\mbox{}\hfill
$
p\colon w\to \nu\gamma
$
\ 
for all \ $w\in {\mathbb C}^n,
\hfill
$
\\[1.5ex]
where 
$\gamma=\mbox{\rm colon}(w_1,\ldots,w_n,\overline{w}_{1},\ldots,\overline{w}_n),$
\vspace{0.5ex}
is a partial integral of the system} (1.1).

{\sl Proof}. If $\nu$ is a common eigenvector of the matrices 
$A_{k},\ k=1,\ldots, 2m,$ then $\nu$ is a solution to system (1.5),
where $\lambda^{k}$ is an eigenvalue of the matrix $A_{k}$ 
corresponding to the eigenvector $\nu.$ We obtain
$
{\frak x}_{k} (\nu\gamma) = \lambda^{k} \nu\gamma$
for all $w\in {\mathbb C}^n,\ k=1,\ldots,2m.$
\vspace{0.5ex}
Therefore  the ${\mathbb R}\!$-li\-ne\-ar function $p$ is a partial integral 
of the system in total differential (1.1).\ \k
\vspace{3.75ex}

{\bf 1.1.2. Autonomous ${\mathbb R}\!$-differentiable first integrals}
\\[1ex]
\indent    
{\bf Theorem 1.1}.
{\it 
Let $\nu^{\theta},\, \theta=1,\ldots,2m+1$ be common eigenvectors of  the matrices  
$A_{k},$ $k=1,\ldots, 2m.$
Then the system {\rm (1.1)} has the ${\mathbb R}\!$-differentiable autonomous first integral
\\[1.5ex]
\mbox{}\hfill                                         
$
\displaystyle
F\colon w\to
\prod\limits_{\theta=1}^{2m+1}
\bigl(\nu^{\theta}\gamma\bigr)^{h_{\theta}}$ \ 
for all $w\in\Omega,
\ \ \ \Omega\subset {\rm D}(F),
$
\hfill {\rm (1.6)}
\\[1ex]
where $h_{1},\ldots,h_{2m+1}$ is a nontrivial solution to the system  
$
\sum\limits_{\theta=1}^{2m+1}\,\lambda_{\theta}^{k}\,h_{\theta} =0,
\ k=1,\ldots,2m,
$
and $\lambda_{\theta}^{k}$ are the eigenvalues of the matrices $A_{k},\ k=1,\ldots, 2m,$ 
corresponding to the common eigenvectors $\nu^{\theta},\ \theta=1,\ldots, 2m+1.$
}
\vspace{0.5ex}

{\sl Proof}.
Suppose $\nu^{\theta}$ are common eigenvectors of  the matrices  $A_{k}$ 
corresponding to the eigenvalues $\lambda_{\theta}^{k},\ k=1,\ldots, 2m,\ \theta=1,\ldots,2m+1,$ respectively.
By lemma 1.1, it follows that the 
${\mathbb R}\!$-linear functions 
$w\to\nu^{\theta}\gamma$ for all $w\in {\mathbb C}^n,\ \theta=1,\ldots, 2m+1$
are partial integrals of the system of equations in total differentials (1.1). Hence,
\\[1.5ex]
\mbox{}\hfill                                         
$
{\frak x}_{{}_{\scriptstyle k}}\,\nu^{\theta}\gamma =
\lambda_{\theta}^{k}\,\nu^{\theta}\gamma
$
\ for all $w\in {\mathbb C}^n,
\ \
k=1,\ldots, 2m,
\quad \theta=1,\ldots, 2m+1.
$
\hfill (1.7)
\\[1.5ex]
\indent
We form the function
\\[1.5ex]
\mbox{}\hfill
$
\displaystyle
F\colon w\to
\prod\limits_{\theta=1}^{2m+1}
\bigl(\nu^{\theta}\gamma\bigr)^{h_{\theta}}
$ 
\ \ for all \ $ w\in\Omega,
\quad 
\Omega\subset {\mathbb C}^n,
\hfill
$  
\\[1.5ex]
where $\Omega$ is a domain (open arcwise connected set) in ${\mathbb C}^n$ 
and $h_{\theta},\,\theta=1,\ldots, 2m+1$ are complex numbers 
with $\sum\limits_{\theta=1}^{2m+1}|h_{\theta}|\ne 0.$
The Lie derivative of $F$ by virtue of (1.1) is equal to 
\\[1.75ex]
\mbox{}\hfill
$
\displaystyle
{\frak x}_{{}_{\scriptstyle k}} F(w)=
\prod\limits_{\theta=1}^{2m+1}
\bigl(\nu^{\theta}\gamma\bigr)^{h_{\theta}-1}\,
\sum\limits_{\theta=1}^{2m+1} h_{\theta}\,
\prod\limits_{l=1,l\ne \theta}^{2m+1}(\nu^{l}\gamma) \
{\frak x}_{{}_{\scriptstyle k}}\,\nu^{\theta}\gamma
$
\ \ for all $w\in\Omega, 
\quad k=1,\ldots, 2m.
\hfill
$
\\[1.5ex]
\indent
Using (1.7), we get
\\[1.5ex]
\mbox{}\hfill
$
\displaystyle
{\frak x}_{{}_{\scriptstyle k}}F(w)=
\sum\limits_{\theta=1}^{2m+1}\lambda_{\theta}^{k}h_{\theta}\,F(w)
$
\ \ for all $w\in\Omega, 
\quad k=1,\ldots, 2m.
\hfill
$
\\[1.5ex]
\indent
If $\sum\limits_{\theta=1}^{2m+1}\lambda_{\theta}^{k}h_{\theta}=0,\
k=1,\ldots 2m,$ then the function (1.6) is an autonomous ${\mathbb R}\!$-differentiable  first integral 
of the system in total differentials (1.1).\ \k
\vspace{1ex}

{\bf Corollary 1.1}.
\vspace{0.35ex}
{\it
Let $\nu^{\theta}$ be common eigenvectors of  the matrices  $A_{k}$ 
corresponding to the eigenvalues $\lambda_{\theta}^{k},\ k=1,\ldots, 2m,\ \theta=1,\ldots,2m+1,$ respectively.
\vspace{0.35ex}
Then the system of equations in total differentials {\rm (1.1)} has the ${\mathbb R}\!$-differentiable autonomous first integral
\\[1.75ex]
\mbox{}\hfill
$
F_{12\ldots 2m(2m+1)}\colon w\to
\prod\limits_{\theta=1}^{2m}
\bigl(\nu^{\theta}\gamma\bigr)^{{}-\delta_{\theta}}
\bigl(\nu^{2m+1}\gamma\bigr)^{\delta}$
\ \ for all $w\in\Omega,
\quad \Omega\subset {\rm D}(F_{12\ldots 2m(2m+1)}),
\hfill
$
\\[1.75ex]
where the determinants $\delta_{\theta},\, \theta=1,\ldots,2m$ 
\vspace{0.35ex}
are obtained by replacing the $\theta\!$-th column of the determinant $\delta=\big|\lambda_{\theta}^{k}\big|$  by  
${\rm colon}\!\bigl(\lambda_{2m+1}^{1},
\ldots\!,\lambda_{2m+1}^{2m}\bigr),$ respectively.
}
\vspace{1ex}

For example, the ${\mathbb R}\!$-linear autonomous system of equations in total differentials
\\[1.5ex]
\mbox{}\hfill                              
$
dw_1=
(2w_1-i(w_2+\overline{w}_1)+(1-i)\overline{w}_2)dz
+(w_1+(2-i)(w_2+\overline{w}_1)+(1-i)\overline{w}_2)d\,\overline{z},
\hfill
$
\\[0.1ex]
\mbox{}\hfill (1.8)
\\
\mbox{}\hfill
$
dw_2=
{}-((2+i)(w_1+\overline{w}_2)+iw_2)dz
-(i(w_1+\overline{w}_2)+(i-1)w_2)d\,\overline{z}
\hfill
$
\\[1.5ex]
has the commuting matrices 
\\[1.5ex]
\mbox{}\hfill
$
A_1 =
\left\|\!
\begin{array}{rrcc}
   2 & -2-i\, & 2+i & 0
\\
  -i &   -i\, & 1+i & i
\\
  -i &    0\, &   1 & i
\\
 1-i & -2-i\, & 2+i & 1+i
\end{array}\!
\right\|
$
\ and \
$
A_2 =
\left\|\!
\begin{array}{crcc}
    1 &  -i\, &   i & 0
\\
  2-i & 1-i\, & 1+i & -2+i
\\
  2-i &   0\, &   2 & -2+i
\\
  1-i &  -i\, &   i & i
\end{array}\!
\right\|.
\hfill
$
\\[2ex]
Therefore the system of equations in total differentials (1.8) is completely solvable. 

The matrices $A_1$ and $A_2$ have the eigenvalues 
\vspace{0.35ex}
$\lambda_{1}^1=1+i, 
\ \lambda_{2}^{1}={}-i,\
\lambda_{3}^{1}=1,
\ \lambda_{4}^1=2,
$
and
\vspace{0.35ex}
$
\lambda_{1}^2=i,
\
\lambda_{2}^{2}=1-i,
\
\lambda_{3}^{2}=2,
\
\lambda_{4}^2=1
$
corresponding to the eigenvectors 
$\nu^1=(0,1,1,1),$
$\nu^2\!=\!(1,1,0,1),\ \nu^3\!=\!(0,1,1,0),\ \nu^4\!=\!(1,0,0,1),$ respectively. 
\vspace{0.35ex}

The solution to the linear homogeneous system 
\\[1.5ex]
\mbox{}\hfill
$
\left\{\!\!
\begin{array}{l}
(1+i)h_1-ih_2+h_3=0,
\\[1ex]
ih_1+(1-i)h_2+2h_3=0
\end{array}
\right.
\iff
\left\{\!\!
\begin{array}{l}
h_1={}-(1+i)h_3,
\\[1ex]
h_2={}-(2+i)h_3
\end{array}
\right.
\hfill
$
\\[2ex]
is $h_1=1+i,\ h_2=2+i, \ h_3={}-1.$

The ${\mathbb R}\!$-differentiable function (by Theorem 1.1)
\\[2ex]
\mbox{}\hfill                                    
$
F\colon w\to
\dfrac{(w_2+\overline{w}_1+\overline{w}_2)^{1+i}
(w_1+w_2+\overline{w}_2)^{2+i}}{w_2+\overline{w}_1}$
\ \ for all 
$w\in \Omega,
$
\hfill (1.9)
\\[2.5ex]
where a domain $\Omega\subset \{w\colon w_2+\overline{w}_1\ne 0\},$ is an autonomous first integral of the system (1.8).

The ${\mathbb R}\!$-differentiable first integral (1.9) is an autonomous general integral of the completely solvable 
system of equations in total differentials (1.8).
\vspace{1ex}

From the entire set of ordinary differential systems induced by the completely solvable 
system of  equations in total differentials (1.1), we extract system
\\[1.5ex]
\mbox{}\hfill                                                           
$
dw_\tau= X_{\tau\zeta}(w)\,dz_{\zeta}+
X_{\tau,m+\zeta}(w)\,d\,\overline{z}_{\zeta},
\ \ \ \tau=1,\ldots, n,
\quad
\zeta\in \{1,\ldots,m\}
\hfill (1.1.\zeta)
$
\\[1.5ex]
such that the matrix $A_{\zeta}$ has the smallest number of elementary divisors [64, p. 147]. 
\vspace{0.75ex}

{\bf Definition 1.1}.
{\it
Let $\nu^{0l}$ be an eigenvector of the matrix $A_{\zeta}$ corresponding to the eigen\-va\-lue $\lambda_{l}^{\zeta}$ with 
elementary divisor of multiplicity $s_{l}.$ 
A non-zero vector $\nu^{\eta l}\in {\mathbb C}^{2n}$ is called a
\textit{\textbf{generalized eigenvector of order}} {\boldmath $\eta$} for $\lambda_{l}^{\zeta}$ if and only if
\\[1.5ex]
\mbox{}\hfill                          
$
(A_{\zeta}-\lambda_{l}^{\zeta} E)\,\nu^{\eta}=\eta\cdot \nu^{\eta-1},
\quad 
\eta=1,\ldots, s_{l}-1,
\hfill (1.10)
$
\\[1.5ex]
where $E$ is the $2n\times 2n$ identity matrix.
}
\vspace{1ex}

Using Lemma 1.1 and (1.10), we obtain
\\[1.5ex]
\mbox{}\hfill                                
$
{\frak x}_{\zeta}\,\nu^{0l}\gamma =
\lambda_{l}^{\zeta}\,\nu^{0l}\gamma,
\ \ \ \ 
{\frak x}_{\zeta}\,\nu^{\eta l}\gamma =
\lambda_{l}^{\zeta}\,\nu^{\eta l}\gamma +
\eta\,\nu^{\eta-1{,}\,l}\gamma$
\ for all $w\in {\mathbb C}^n,
\ \ \eta=1,\ldots, s_{l}-1.
$
\hfill (1.11) 
\\[2ex]
\indent
The following lemmas are needed for the sequel.
\vspace{1.25ex}

{\bf Lemma 1.2}.                                           
{\it
Let $\nu^{\,0l}$ be a common eigenvector of the matrices $A_k$ corresponding to the eigenvalues
$\lambda^k_{l}, \ k=1,\ldots,2m,$ respectively. Let $\nu^{\,\eta l},\ \eta =1,\ldots, s_l-1$ be generalized eigenvectors
of the matrix $A_\zeta$ 
\vspace{0.35ex}
corresponding to the eigenvalue $\lambda^{\zeta}_{l}$ with elementary divisor of multiplicity 
$s_l\ (s_{l}\geq 2).$ If the system $(1.1.\zeta)$ hasn't the first integrals 
\\[2ex]
\mbox{}\hfill                               
$
\displaystyle
F_{k\eta l}^{\,\zeta}\colon w\to
{\frak x}_k\, \Psi_{\eta l}^{\zeta}(w)$ 
\ for all $w\in \Omega, 
\ \ k=1,\ldots, 2m, \ \, k\ne\zeta,
\ \ \  \eta =1,\ldots, s_l-1,
$
\hfill {\rm (1.12)}
\\[2.25ex]
then 
\\[1.5ex]
\mbox{}\hfill                             
$
{\frak x}_{\zeta}\, \Psi_{\eta l}^{\zeta}(w)=
\left[\!\!
\begin{array}{lll}
1\! & \text{for all}\ \ w\in \Omega, & \eta =1,
\\[1ex]
0\! & \text{for all}\ \ w\in \Omega, & \eta =2,\ldots, s_{l}-1,
\end{array}
\right.
$
\hfill {\rm (1.13)}
\\[2ex]
\mbox{}\hfill                               
$
{\frak x}_{k}\, \Psi_{\eta l}^{\zeta}(w)=\mu_{\eta l}^{k\zeta}={\rm const}$
\ for all $w\in \Omega,
\ \ k=1,\ldots, 2m,\ \ k\ne \zeta,
\ \ \  \eta =1,\ldots, s_l-1,
\hfill
$ 
\\[2.5ex]
where $\Psi_{\eta l}^{\zeta}\colon \Omega\to {\mathbb C},\  \eta =1,\ldots, s_l-1$ is a solution to the system 
\\[1.75ex]
\mbox{}\hfill                             
$
\nu^{\,\eta l}\gamma=
{\displaystyle \sum\limits_{\delta=1}^{\eta} }
\binom{\eta -1}{\delta-1}\Psi_{\delta l}^{\zeta}(w)\cdot \nu^{\,\eta-\delta,l}\gamma,
\quad 
\eta=1,\ldots, s_l-1,
\quad
\Omega\subset \{w\colon \nu^{0l}\gamma\ne 0\}.
$
\hfill {\rm (1.14)}
\\[1.75ex]
}
\indent
{\sl Proof}. 
\vspace{0.35ex}
The system (1.14) has the determinant $(\nu^{0l}\gamma)^{s_{l}-1}.$ Therefore
there exists the solution $\Psi_{\eta l}^{\zeta}, \eta=1,\ldots, s_l-1$ on a domain $\Omega\subset \{w\colon \nu^{0l}\gamma\ne 0\}$ 
of the system (1.14).  
\vspace{0.55ex}

The proof of the lemma is by induction on $\eta.$

For $\eta=1$ and $\eta=2,$ the assertion (1.13) follows from (1.11).
\vspace{0.35ex}

Assume that (1.13) for $\eta=1,\ldots,\varepsilon-1$ is true. Using (1.11) and (1.14), we get
\\[1.75ex]
\mbox{}\hfill
$
\displaystyle
{\frak x}_{\zeta}\,\nu^{\varepsilon l}\gamma =
\lambda_{l}^{\zeta}\sum\limits_{\delta=1}^{\varepsilon}
{\textstyle\binom{\varepsilon -1}{\delta-1}}\,\Psi_{\delta l}^{\zeta}(w)\,
\nu^{\varepsilon-\delta{,}\,l}\gamma +  
(\varepsilon -1)\sum\limits_{\delta=1}^{\varepsilon -1}
{\textstyle\binom{\varepsilon -2}{\delta-1}}\,\Psi_{\delta l}^{\zeta}(w)\,
\nu^{\varepsilon-\delta-1{,}\,l}\gamma  \ +
\hfill
$
\\[2ex]
\mbox{}\hfill
$
\displaystyle
+\ \nu^{\varepsilon -1{,}\,l}\gamma \,+\,
\nu^{0l}\gamma\ {\frak x}_{\zeta}\Psi_{\varepsilon l}^{\zeta}(w)$
\ \ for all $w\in \Omega.
\hfill
$
\\[2ex]
\indent
Combining (1.14) for 
\vspace{0.35ex}
$\eta=\varepsilon-1$ and $\eta=\varepsilon,$ (1.11) for $\eta=\varepsilon,$ and 
$\nu^{0l}\gamma\not\equiv 0$ in $\C^n,$ we obtain 
${\frak x}_{\zeta}\,\Psi_{\varepsilon l}^{\zeta}(w)=0$ for all $w\in\Omega.$
\vspace{0.35ex}
So by the principle of mathematical induction, the statement (1.13) is true for every $\eta=1,\ldots, s_{l}-1$ and 
$\zeta\in \{1,\ldots, m\}.$ 
\vspace{0.35ex}

Taking into account (1.4) and (1.12), we have the statement (1.13) is true for $k\ne \zeta.$ \k
\vspace{1.5ex}

{\bf Lemma 1.3}
{\it
Under the conditions of Lemma {\rm 1.2}, we have
\\[2ex]
\mbox{}\hfill                                                                
$
{\frak x}_{k}\,\nu^{\,\eta l}\gamma=
{\displaystyle \sum\limits_{\delta=0}^{\eta}}
\binom{\eta}{\delta}\,\mu_{\delta l}^{ k\zeta}\cdot\nu^{\,\eta-\delta, l}\gamma$
\ \ for all $w\in \Omega,
\ \ \ k=1,\ldots, 2m, 
\ \ \eta=1,\ldots, s_{l}-1,
$
\hfill {\rm (1.15)}
\\[2ex]
where $\mu_{0l}^{k\zeta}=\lambda^k_{l}, \ \mu_{\eta l}^{k\zeta}={\frak x}_{k}\,\Psi_{\eta l}^{\zeta}(w),
\ \eta=1,\ldots, s_{l}-1, \ k=1,\ldots, 2m.$ 
}
\vspace{1.5ex}

{\sl Proof}. The proof of Lemma 1.3 is by induction on $\eta.$
 \vspace{0.35ex}

Let $\eta=1.$ Using (1.14), we get 
\\[1.5ex]
\mbox{}\hfill                                  
$
\nu^{\,1 l}\gamma=
\Psi_{1 l}^{\zeta}(w)\cdot\nu^{0 l}\gamma$
\  for all $w\in \Omega, 
\ \ \ \Omega\subset \{w\colon \nu^{0l}\gamma\ne 0\}.
$
\hfill {\rm(1.16)}
\\[1ex]
\indent
Then
\\[1.25ex]
\mbox{}\hfill                                 
$
{\frak x}_{k}\nu^{\,1l}\gamma=
{\frak x}_{k}\Psi_{1 l}^{\zeta}(w)\cdot\nu^{0l}\gamma +
\Psi_{1 l}^{\zeta}(w)\cdot {\frak x}_{k}\, \nu^{0l}\gamma$
\  for all $w\in \Omega, 
\ \ \ k=1,\ldots, 2m. 
\hfill
$
\\[-1.5ex]

\newpage

Taking into account Lemma 1.1, Lemma 1.2 and (1.16),we obtain 
\\[1.75ex]
\mbox{}\hfill                                                            
$
{\frak x}_{k}\nu^{\,1 l}\gamma=
\mu_{1 l}^{k\zeta}\cdot\nu^{0l}\gamma +
\mu_{0 l}^{k\zeta}\,\Psi_{1l }^{\zeta}(w)\cdot \nu^{0l}\gamma=
\mu_{0l}^{k\zeta}\cdot\nu^{1l}\gamma+
\mu_{1l}^{k\zeta}\cdot \nu^{0l}\gamma,
\quad
k=1,\ldots, 2m, 
\hfill 
$
\\[2ex]
where $\mu_{0 l}^{k\zeta}=\lambda^k_l,\  
 {\frak x}_{k}\,\Psi_{1 l}^{\zeta}(w)=\mu_{1l}^{k\zeta}, \ k=1,\ldots, 2m.$ So (1.15)  for $\eta=1$ is true.
\vspace{0.75ex}

Suppose that the assertion of the lemma is valid for $\eta=1,\ldots, \varepsilon-1.$ From (1.14), we get
\\[2ex]
\mbox{}\hfill                               
$
{\frak x}_{k}\nu^{\varepsilon l}\gamma= 
{\displaystyle \sum\limits_{\delta=1}^{\varepsilon} }
\binom{\varepsilon-1}{\delta-1}\,{\frak x}_{k}\Psi_{\delta l}^{\zeta}(w)\,\nu^{\,\varepsilon-\delta, l}\gamma +
{\displaystyle \sum\limits_{\delta=1}^{\varepsilon} }
\binom{\varepsilon-1}{\delta-1}\,\Psi_{\delta l}^{\zeta}(w)\, {\frak x}_{k}\nu^{\,\varepsilon-\delta,l}\gamma$
for all $w\in \Omega,
\, k\!=\!1,\ldots, 2m. 
\hfill
$
\\[1.75ex]
\indent
By the induction hypothesis, we have
\\[2ex]
\mbox{}\hfill                               
$
{\frak x}_{k}\nu^{\varepsilon l}\gamma =
{\displaystyle \sum\limits_{\delta=1}^{\varepsilon} }
\binom{\varepsilon-1}{\delta-1}\,\mu_{\delta l}^{k\zeta}\cdot \nu^{\,\varepsilon-\delta,l}\gamma  +
{\displaystyle \sum\limits_{\delta=1}^{\varepsilon} }
\binom{\varepsilon-1}{\delta-1}\,\Psi_{\delta l}^{\zeta}(w)\cdot 
{\displaystyle \sum\limits_{\varkappa=0}^{\varepsilon-\delta} }
\binom{\varepsilon-\delta}{\varkappa}\,\mu_{\varkappa l}^{k\zeta}\cdot \nu^{\,\varepsilon-\delta-\varkappa,l}\gamma =
\hfill
$
\\[2ex]
\mbox{}\hfill                               
$
={\displaystyle \sum\limits_{\delta=1}^{\varepsilon} }
\binom{\varepsilon-1}{\delta-1}\,\mu_{\delta l}^{k\zeta}\cdot \nu^{\,\varepsilon-\delta,l}\gamma  +
{\displaystyle \sum\limits_{\alpha=0}^{\varepsilon-1} }
\binom{\varepsilon-1}{\alpha}\,\mu_{\alpha l}^{k\zeta}\cdot
{\displaystyle \sum\limits_{\beta=1}^{\varepsilon-\alpha} }
\binom{\varepsilon-\alpha-1}{\beta-1}\,\Psi_{\beta l}^{\zeta}(w)\cdot 
\nu^{(\varepsilon-\alpha)-\beta,l}\gamma=
\hfill
$
\\[2ex]
\mbox{}\hfill                               
$
={\displaystyle \sum\limits_{\delta=1}^{\varepsilon} }
\binom{\varepsilon-1}{\delta-1}\,\mu_{\delta l}^{k\zeta}\cdot \nu^{\,\varepsilon-\delta, l}\gamma +
{\displaystyle \sum\limits_{\alpha=0}^{\varepsilon-1} }
\binom{\varepsilon-1}{\alpha}\,\mu_{\alpha l}^{k\zeta}\cdot \nu^{\varepsilon-\alpha,l}\gamma=
\hfill
$
\\[2ex]
\mbox{}\hfill                               
$
={\displaystyle \sum\limits_{\delta=1}^{\varepsilon-1} }
\binom{\varepsilon-1}{\delta-1}\,\mu_{\delta l}^{k\zeta}\cdot \nu^{\,\varepsilon-\delta,l}\gamma +
\mu_{\varepsilon l}^{k\zeta}\cdot \nu^{0l}\gamma +\mu_{0l}^{k\zeta}\cdot \nu^{\varepsilon l}\gamma +
{\displaystyle \sum\limits_{\delta=1}^{\varepsilon-1} }
\binom{\varepsilon-1}{\delta}\,\mu_{\delta l}^{k\zeta}\cdot \nu^{\varepsilon-\delta,l}\gamma=
\hfill
$
\\[2.5ex]
\mbox{}\hfill                               
$
=\mu_{0l}^{k\zeta}\cdot \nu^{\varepsilon l}\gamma +
{\displaystyle \sum\limits_{\delta=1}^{\varepsilon-1} }
\binom{\varepsilon}{\delta}\,\mu_{\delta l}^{k\zeta}\cdot \nu^{\,\varepsilon-\delta,l}\gamma +
\mu_{\varepsilon l}^{k\zeta}\cdot \nu^{0l}\gamma =
{\displaystyle \sum\limits_{\delta=0}^{\varepsilon} }
\binom{\varepsilon}{\delta}\,\mu_{\delta l}^{k\zeta}\cdot \nu^{\varepsilon-\delta,l}\gamma,
\ \ \ k=1,\ldots, 2m.
\hfill
$
\\[2ex]
\indent
Thus by the principle of induction, the statement (1.15) is true 
\vspace{2.25ex}
for $\eta=1,\ldots, s_{l}-1.\!\!$ \k

{\bf Theorem 1.2}.
{\it
\vspace{0.5ex}
Let the assumptions of Lemma {\rm 1.2} with $l=1,\ldots r\ \Bigl(\, \sum\limits_{l=1}^{r}s_{l}\geq m+1\Bigr)$ hold.
Then the completely solvable system {\rm (1.1)} has the autonomous first integral
\\[1.75ex]
\mbox{}\hfill                                                                  
$
\displaystyle
F\colon w\to
\prod\limits_{\xi=1}^{\alpha}\bigl(\nu^{0\xi} \gamma\bigr)^{h_{0 \xi}}
\exp\sum\limits_{q=1}^{\varepsilon_{\xi}}\,h_{q\xi} \Psi_{q\xi}^{\zeta}(w)$
\ \ for all $w\in \Omega,
\quad \Omega\subset {\rm D}(F),
$
\hfill {\rm (1.17)}
\\[2ex]
where $\sum\limits_{\xi=1}^{\alpha}\varepsilon_{\xi}=2m-\alpha+1, \
\varepsilon_{\xi}\leq s_{\xi}\!-\!1,\, \xi=1,\ldots,\alpha,\, \alpha\leq r,\!$ and 
$h_{q\xi},\, q\!=\!0,\ldots,\varepsilon_{\xi},\, \xi\!=\!1,\ldots, \alpha$ is a nontrivial solution to
the linear homogeneous algebraic system of equations 
\\[1.75ex]
\mbox{}\hfill
$
\displaystyle
\sum\limits_{\xi=1}^{\alpha}
\bigl(\lambda_{\xi}^{k}\,h_{0\xi}
+ \sum\limits_{q=1}^{\varepsilon_{\xi}}
\mu_{q\xi}^{k\zeta}\,h_{q\xi}\big) = 0,
\ \ k=1,\ldots, 2m.
\hfill
$
\\[2.25ex]
}
\indent
{\sl Proof}.
The Lie derivative of (1.17) by virtue of (1.1) is equal to 
\\[1.5ex]
\mbox{}\hfill
$
\displaystyle
{\frak x}_{k}\,F(w) =
\sum\limits_{\xi=1}^{\alpha}
\bigl(\lambda_{\xi}^{k}h_{0\xi} +
\sum\limits_{q=1}^{\varepsilon_{\xi}}\mu_{q\xi}^{k\zeta}h_{q\xi}\bigr)F(w)$
\ \ for all $w\in\Omega,
\quad 
k=1,\ldots, 2m.
\hfill
$
\\[2ex]
\indent
If
\vspace{0.5ex}
$
\sum\limits_{\xi=1}^{\alpha}
\bigl(\lambda_{\xi}^{k}h_{0\xi} +
\sum\limits_{q=1}^{\varepsilon_{\xi}}
\mu_{q\xi}^{k\zeta}h_{q\xi}\bigr) = 0,\ k=1,\ldots, 2m,$ then 
the ${\mathbb R}\!$-differentiable function (1.17) is an autonomous first integral of the 
completely solvable system (1.1). \ \k

\newpage

As an example, the completely solvable ${\mathbb R}\!$-linear system of equations in total differentials
\\[1.5ex]
\mbox{}\hfill                                       
$
\begin{array}{l}
dw_1=
((1+i)w_1+iw_2-\overline{w}_1-\overline{w}_2)dz
+(w_1+iw_2-\overline{w}_1-\overline{w}_2)d\,\overline{z},
\\[1ex]
dw_2=
(w_2+\overline{w}_1+\overline{w}_2)dz
+((1-i)w_2+\overline{w}_1+\overline{w}_2)d\,\overline{z},
\\[1ex]
dw_3=
({}-w_1+w_2+w_3-i\,\overline{w}_2)dz
+({}-w_1+w_2+(1-i)w_3-i\,\overline{w}_2)d\,\overline{z}
\end{array}
$
\hfill  (1.18)
\\[1.5ex]
has two eigenvalues: $\!\lambda_1^1\!=\!1+i$ with elementary divisor $(\lambda^1\!-\!1\!-\!i)^3\!\!$ and 
$\lambda_2^1\!=\!1\!$ with elementary divisor $(\lambda^1\!-\!1)^3.$
The $\lambda_1^1\!=\!1+i$ corresponding to the  eigenvector $\nu^{01}\!=\!(1, 1, 0, 0, 0, 0)\!$
and to the generalized eigenvectors $\!\nu^{11}\!=\!(0, 0, 0, 0, 1, 0), \, \nu^{21}\!=\!(0, 1, 0, 0, 0, 1).\!\!$
The eigenvalue $\!\lambda_2^1\!=\!1$ corresponding to 
the  $\nu^{02}\!=(0, 0, 0, 1, 1, 0),\,  \nu^{12}\!=(0, 1, 0, 0, 0, 0), \, \nu^{22}\!=(0, 0, 1, 0, 1, 0).$

The scalar functions (see (1.14))
\\[1.5ex]
\mbox{}\hfill
$
\Psi_{11}^1\colon w\to
\dfrac{\overline{w}_2}{w_1+w_2}\,,
\quad
\Psi_{21}^1\colon w\to
\dfrac{(w_1+w_2)(w_2+\overline{w}_3)-\overline{w}_2^{\,2}}{(w_1+w_2)^2}$
\ \ for all $w\in\Omega,
\hfill
$
\\[1.75ex]
where $\Omega\subset \{w\colon w_1+w_2\ne 0\}\subset {\mathbb C}^3.$
The ${\mathbb R}\!$-differentiable first integrals (by Theorem 1.2)
\\[1.5ex]
\mbox{}\hfill                             
$
F_{1}\colon w\to \Psi_{21}^1(w)$ 
\ \ for all $w\in \Omega
$
\hfill (1.19)
\\[1ex]
and
\\[1ex]
\mbox{}\hfill                                    
$
F_{2}\colon w\to
\dfrac{\overline{w}_1+\overline{w}_2}{w_1+w_2} \
\exp\Bigl(i\,\dfrac{\overline{w}_2}{w_1+w_2}\Bigr)$
\ \ for all $w\in\Omega
$
\hfill (1.20)
\\[2ex]
are an autonomous general integral of the system of equations in total differentials (1.18).
\\[3ex]
\indent
{\bf 1.1.3. Nonautonomous ${\mathbb R}\!$-differentiable  first integrals}
\\[1ex]
\indent
{\bf Theorem 1.3.}
{\it 
Suppose $\nu$ is a common eigenvector of  the matrices  $A_{k}$ cor\-res\-pon\-ding to the eigenvalues 
$\lambda^{k},\ k=1,\ldots, 2m,$ respectively.
Then the ${\mathbb R}\!$-differentiable function
\\[1.5ex]
\mbox{}\hfill                     
$
\displaystyle
F\colon (z,w)\to
(\nu\gamma)\exp\biggl({}-
\sum_{j=1}^{m}(\lambda^{j}z_{j}+
\lambda^{m+j}\,\overline{z}_{j})\biggr)$
\ \ for all $(z,w)\in {\mathbb C}^{m+n}
\hfill 
$
\\[1.5ex]
is a first integral of the system of equations in total differentials {\rm(1.1)}.
}
\vspace{0.5ex}

{\sl Proof}.  Using Lemma 1.1, we obtain 
$
{\frak X}_{k}F(z,w) =0$ for all $(z,w)\in {\mathbb C}^{m+n},
\ k=1,\ldots, 2m,$
where the linear nonautonomous differential operators 
\\[1.5ex]
\mbox{}\hfill
$
{\frak X}_j(z,w)=\partial_{z_j}+{\frak x}_j(w),
\ 
{\frak X}_{m+j}(z,w)=
\partial_{{}_{\scriptstyle\overline{z}_j}}+ {\frak x}_{m+j}(w)$
for all $(z,w)\!\in\! {\mathbb C}^{m+n},
\,  j\!=\!1,\ldots, m.\k
\hfill
$
\\[1.5ex]
\indent
Consider the system (1.8). 
Using the eigenvector $\nu^1=(0,1,1,1)$ corresponding to the eigenvalues 
$\lambda_1^1=1+i$ and $\lambda_1^2=i,$
we can build the first integral (by Theorem 1.3)
\\[1.5ex]
\mbox{}\hfill                            
$
F\colon (z,w)\to
(w_2+\overline{w}_1+\overline{w}_2)\exp({}-(1+i)z-i\,\overline{z}\,)$
\ \ for all $(z,w)\in {\mathbb C}^3.
$
\hfill (1.21)
\\[1.5ex]
\indent
The first integrals (1.9) and (1.21) are a general integral of the system (1.8).
\vspace{1.5ex}

{\bf Theorem 1.4.}
{\it 
Suppose the system {\rm (1.1)} satisfies the conditions of Lemma {\rm 1.2.} Then the completely solvable 
system {\rm (1.1)} has the ${\mathbb R}\!$-differentiable first integrals
\\[1.5ex]
\mbox{}\hfill                        
$
\displaystyle
F_{\eta}\colon (z,w)\to
\Psi_{\eta l}^{\zeta}(w) -
\sum_{j=1}^{m}\bigl(
\mu_{\eta l}^{j\zeta}z_{j}+\mu_{\eta l}^{m+j,\zeta}\overline{z}_{j}\bigr)$
\ for all $(z,w)\in G,
\ \ \eta=1,\ldots, s_{l}-1,
$
\hfill {\rm (1.22)}
\\[1.5ex]
where a domain $G\subset {\mathbb C}^{m+n},$ the functions 
\vspace{0.25ex}
$\Psi_{\eta l}^{\zeta}$ are the solution to system {\rm (1.14)}, the numbers 
$\mu_{\eta l}^{k\zeta}={\frak x}_{k}\Psi_{\eta l}^{\zeta}(w),\
\eta=1,\ldots, s_{l}-1,\ k=1,\ldots, 2m.$
}
\vspace{1.25ex}

{\sl Proof.}
The Lie derivative of (1.22) by virtue of (1.1) is  
\\[1.5ex]
\mbox{}\hfill
$
{\frak X}_k F_{\eta}(z,w) =
{}-\mu_{\eta l}^{k\zeta} + {\frak x}_k \Psi_{\eta l}^{\zeta}(w)$
\ for all $(z,w)\in G,
\ \ k=1,\ldots, 2m,
\ \ \eta=1,\ldots, s_{l}-1.
\hfill
$
\\[1.5ex]
\indent
Taking into account Lemma 1.2, we get the functions (1.22) are first integrals of (1.1).\k

For example, the system (1.18) has the numbers $\mu_{11}^{11}=1, \  \mu_{11}^{21}=1,$ and the first integral 
\\[2ex]
\mbox{}\hfill                      
$
F\colon (z,w)\to
\dfrac{\overline{w}_2}{w_1+w_2} -z -\overline{z}$
\ \ for all $(z,w)\in\C\times \Omega
$
\ \ \ (by Theorem 1.4).
\hfill\mbox{}
\\[2ex]
\indent
The ${\mathbb R}\!$-differentiable first integrals (1.19),  (1.20), and $F$ are a
general integral on a domain $\C\times \Omega$ of the system (1.18), 
where a domain $\Omega\subset \{w\colon w_1+w_2\ne 0\}.$
\\[3ex]
\indent
{\bf 1.2. $\R\!$-linear nonhomogeneous systems in total differentials}
\\[1ex]
\indent
Let us consider a nonhomogeneous system of equations in total differentials
\\[1.75ex]
\mbox{}\hfill                                          
$
\displaystyle
dw =\sum\limits_{j=1}^m\bigl((B_j\,\gamma+f_j(z))\,dz_j+
(B_{m+j}\,\gamma+f_{m+j}(z))\,d\,\overline{z}_j\bigr)
$
\hfill (1.23)
\\[1.75ex]
corresponding to the $\R\!$-linear homogeneous system (1.1), where
the matrices $B_1,\ldots, B_{2m}$ are transpose of $A_1,\ldots, A_{2m},$ respectively, 
and the vector functions
\\[1.5ex]
\mbox{}\hfill
$
f_{k}\colon z\to {\rm colon}(f_{k1}(z),\ldots, f_{kn}(z))$
\ for all $z\in V,
\quad
k=1,\ldots, 2m
\hfill
$
\\[1.5ex]
are continuously ${\mathbb R}\!$-differentiable on a domain $V\subset {\mathbb C}^{m}.$
\vspace{0.35ex}

The Frobenius conditions for the total solvability of system 
of equations in total differentials (1.23) are the relations (1.4) and 
\\[1.5ex]
\mbox{}\hfill                                                         
$
\partial_{z_j}f^{\zeta}(z)+B_{\zeta}f^j(z)=
\partial_{z_\zeta}f^{j}(z)+B_{j}f^{\zeta}(z)$
\ for all $z\in V,
\quad 
j,\zeta=1,\ldots, m,
\hfill 
$
\\[2.25ex]
\mbox{}\hfill                                                         
$
\partial_{\overline{z}_j}f^{m+\zeta}(z)+B_{m+\zeta}f^{m+j}(z)=
\partial_{\overline{z}_\zeta}f^{m+j}(z)+B_{m+j}f^{m+\zeta}(z)$
\ for all $z\in V,
\ \ 
j,\zeta=1,\ldots, m,
\hfill 
$
\\[2.25ex]
\mbox{}\hfill                                                         
$
\partial_{z_j}f^{m+\zeta}(z)+B_{m+\zeta}f^{j}(z)=
\partial_{\overline{z}_\zeta}f^{j}(z)+B_{j}f^{m+\zeta}(z)$
\ for all $z\in V,
\quad 
j,\zeta=1,\ldots, m,
\hfill 
$
\\[2ex]
where the vector functions 
\\[1.75ex]
\mbox{}\hfill
$
f^{j}\colon z\to {\rm colon}(f_{j1}(z),\ldots, f_{jn}(z),\overline{f}_{m+j,1}(z),\ldots, \overline{f}_{m+j,n}(z))$
\ for all $z\in V,
\quad
j=1,\ldots, m,
\hfill
$
\\[2.25ex]
\mbox{}\hfill
$
f^{m+j}\colon z\to {\rm colon}(f_{m+j,1}(z),\ldots, f_{m+j,n}(z),\overline{f}_{j1}(z),\ldots, \overline{f}_{jn}(z))$
for all $z\in V,
\ j=1,\ldots, m.
\hfill
$
\\[3ex]
\indent
{\bf Theorem 1.5}.\!
{\it 
Let the assumptions of Lemma {\rm 1.2} hold. 
Then the completely solvable sys\-tem  of equations in total differentials {\rm (1.23)} has the first integrals
\\[1.75ex]
\mbox{}\hfill                            
$
F_{\eta}\colon (z,w)\to \nu^{\,\eta l}\gamma\cdot \varphi(z) -
{\displaystyle \sum\limits_{\tau=1}^{\eta} }
K_{\tau-1}^{\eta}(z)\cdot F_{\tau-1}(z,w)-C_{\eta}(z)
\hfill                                  
$
\\[1.5ex]
\mbox{}\hfill                                  
for all $(z,w)\in \widetilde{V}\times \Omega,
\ \ \ \
\eta=0,\ldots, s_{l}-1,
\quad
 \Omega\subset \{w\colon \nu^{0l}\gamma\ne 0\},
\ \ \ \widetilde{V}\subset V,
\hfill 
$
\\[2.25ex]
where the functions 
\\[2ex]
\mbox{}\hfill
$
\displaystyle
\varphi\colon z\to \exp\biggl({}-\sum\limits_{k=1}^{2m} \mu_{0l}^{k\zeta}\,u_k\biggr)$
\ \ for all $z\in \widetilde{V},
\quad
u=(z_1,\ldots,z_m,\overline{z}_1,\ldots,\overline{z}_m),
\hfill
$
\\[2.25ex]
\mbox{}\hfill
$
K_{\tau-1}^{\eta}\colon z\to 
{\displaystyle \int} 
{\displaystyle \sum\limits_{k=1}^{2m}}
\biggl(\binom{\eta}{\tau-1}\,\mu_{\eta-\tau+1,l}^{\,k\zeta}+
{\displaystyle \sum\limits_{\delta=1}^{\eta-\tau}}
\binom{\eta}{\delta}\,\mu_{\delta l}^{k\zeta}\cdot K_{\tau-1}^{\eta-\delta}(z)\biggr)\,du_k,
\ 
\tau=1,\ldots, \eta,
\ \ \eta=1,\ldots, s_{l}-1,
\hfill
$
\\[2.25ex]
\mbox{}\hfill
$
C_{\eta}\colon z\to 
{\displaystyle \int\sum\limits_{k=1}^{2m}} 
\biggl(\nu^{\,\eta l}f^k(z)\cdot \varphi(z)+
{\displaystyle \sum\limits_{\tau=1}^{\eta}}
\binom{\eta}{\tau}\,\mu_{\tau l}^{k\zeta}\cdot C_{\eta-\tau}(z)\!\!\biggr)\,du_k$
for all $z\in \widetilde{V},
\ \eta=0,\ldots, s_{l}-1,
\hfill
$
\\[2ex]
the numbers $\mu_{0l}^{k\zeta}=\lambda^k_{l}, \ \mu_{\eta l}^{k\zeta}={\frak x}_{k}\Psi_{\eta l}^{\zeta}(w),\
\eta=1,\ldots, s_{l}-1, \ k=1,\ldots, 2m.$ 
}

\newpage

\mbox{}
\\
\centerline{
\large\bf  
2. First integrals of  linear real systems in total differentials}
\\[1.5ex]
\indent
{\bf 2.1. Linear real homogeneous systems in total differentials}
\\[1ex]
\indent
Consider an autonomous system of equations in total differentials
\\[1.75ex]
\mbox{}\hfill                                  
$
dx = A(x)\,dt,
$
\hfill (2.1)
\\[1.75ex]
where $x=(x_{1},\ldots,x_{n})\!\in {\mathbb R}^n,\, t=(t_{1},\ldots,t_{m})\!\in {\mathbb R}^{m},$ 
the entries of the matrix \vspace{0.35ex}
$\!A(x)\!=\!\| a_{{}_{\scriptstyle ij}}(x)\|\!$
(with $n$ rows and $m$ columns) are linear homogeneous functions
\\[1.75ex]
\mbox{}\hfill
$
\displaystyle
a_{{}_{\scriptstyle ij}}\colon x\to
\sum\limits_{\xi=1}^{n}
a_{{}_{\scriptstyle ij\xi}} x_{{}_{\scriptstyle \xi}}$
\  for all $x\in {\mathbb R}^{n}
\ \ \ (a_{{}_{\scriptstyle ij\xi}}\in {\mathbb R}, \
\xi=1,\ldots, n,\ j=1,\ldots, m,\ i=1,\ldots, n),
\hfill
$
\\[1.75ex]
and $dx=\mbox{colon}(dx_{1},\ldots,dx_{n})$ and $dt=\mbox{colon}(dt_{1},\ldots,dt_{m})$ are vector columns.
\vspace{0.35ex}

Assume that this system is completely solvable. 
\!The Frobenius condition [5, p.\! 19]
for the total solvability of system (2.1) in 
terms of the Poisson brackets is given by the relations
\\[1.75ex]
\mbox{}\hfill                                                                          
$
[{\frak p}_{j}(x), {\frak p}_{\zeta}(x)]  = 0$
\ \ for all $x\in {\mathbb R}^{n},
\quad 
j=1,\ldots,m,
\ \ \zeta=1,\ldots, m,
$
\hfill (2.2)
\\[1.75ex]
where the linear autonomous differential operators
\\[1.75ex]
\mbox{}\hfill
$
\displaystyle
{\frak p}_{j}(x)=
\sum\limits_{i=1}^{n}
a_{{}_{\scriptstyle ij}}(x)
\partial_{{}_{\scriptstyle x_{{}_{\scriptsize i}}}}$
\ \ for all $x\in {\mathbb R}^{n},
\quad j=1,\ldots, m.
\hfill
$
\\[1.75ex]
\indent
The Frobenius conditions (2.2) for system (2.1) are equivalent [60, p. 73]:
\\[2ex]
\mbox{}\hfill
$
A_{j}A_{\zeta} = A_{\zeta}A_{j},
\quad 
j=1,\ldots,m,
\ \ \zeta  =1,\ldots, m,
\hfill
$
\\[2ex]
where $A_{j}= \bigl\| a_{{}_{\scriptstyle \xi j i}} \bigr\|, \ j=1,\ldots,m$ are real $n\times n$ matrices.
\vspace{1.5ex}

{\bf 2.1.1. Partial integrals}. The complex-valued linear homogeneous function
\\[1.75ex]
\mbox{}\hfill
$
\displaystyle
p \colon x\to
\sum\limits_{\xi=1}^{n} b_{\xi} x_{\xi}$
\ \ for all  $x\in {\mathbb R}^{n}
\quad (b_{\xi}\in {\mathbb C}, \ \xi=1,\ldots, n)
\hfill
$
\\[1.75ex]
is a {\it partial integral} of the system in total differentials (2.1) iff
\\[1.75ex]
\mbox{}\hfill                                
$
{\frak p}_{j}\,p(x) =  \lambda^{j}\,p(x)$
\ \ for all $x\in {\mathbb R}^{n},
\ \ \ \lambda^{j}\in {\mathbb C},
\ \ j=1,\ldots, m.
\hfill 
$
\\[1.75ex]
\indent
\!This system of identities 
\vspace{0.35ex}
is equivalent to the linear homogeneous system of equa\-ti\-ons 
$\bigl(A_{j} -  \lambda^{j} E\bigr)b = 0,\,  j\!=\!1,\ldots, m,$ where
\vspace{0.5ex}
$E$ is the $n\times n$ identity matrix, $b\!=\!\mbox{colon}(b_{1},\ldots,b_{n}).$ 

The proof of the following statement is similar to that of Lemma 1.1.
\vspace{0.75ex}

{\bf  Lemma 2.1}.
\vspace{0.5ex}
{\it
Suppose $\nu\in {\mathbb C}^{n}$ is a common eigenvector of the matrices 
$A_{j},\ j=1,\ldots, m.$ Then the linear function
$p\colon x\to \nu x$ for all $x\in {\mathbb R}^n
$
\vspace{1ex}
is a partial integral of the system} (2.1).

The following properties are needed for the sequel.
\vspace{0.75ex}

{\bf  Property 2.1} ([11]).
{\it
Suppose  $\nu={\stackrel{*}{\nu}}+\widetilde{\nu}\,i\ \,
({\stackrel{*}{\nu}}={\rm Re}\,\nu,\ \widetilde{\nu}={\rm Im}\,\nu)$
is a common eigenvector of the matrices $A_{j}$
cor\-res\-pon\-ding to the eigenvalues
\vspace{0.35ex}
$\lambda^j={\stackrel{*}{\lambda}}{}^j+\widetilde{\lambda}{}^j\,i\ 
({\stackrel{*}{\lambda}}{}^j={\rm Re}\,\lambda^j,\ \widetilde{\lambda}{}^j={\rm Im}\,\lambda{}^j),$ $j=1,\ldots, m,$
respectively. Then the real-valued function
\\[1.5ex]
\mbox{}\hfill                                      
$
P\colon x\to ({\stackrel{*}{\nu}}x)^2\, +\,  (\widetilde{\nu}x)^2$
\ for all $x\in {\mathbb R}^n
\hfill 
$
\\[1.5ex]
is a partial integral of the sys\-tem of  equations in total differentials {\rm (2.1)} and
\\[1.5ex]
\mbox{}\hfill                                   
$
{\frak p}_j P(x) =
2\,{\stackrel{*}{\lambda}}{}^j\, P(x)$
\ for all $x\in {\mathbb R}^n,
\quad
j= 1,\ldots,m.
\hfill 
$
\\[-1.5ex]
}

\newpage

{\bf  Property 2.2} ([11]).
{\it
Let  $\nu={\stackrel{*}{\nu}}+\widetilde{\nu}\,i$
be a common eigenvector of the matrices $A_{j}$
cor\-res\-pon\-ding to the eigenvalues
\vspace{0.35ex}
$\lambda^j={\stackrel{*}{\lambda}}{}^j+\widetilde{\lambda}{}^j\,i,\ j=1,\ldots, m,$ respectively. Then 
\\[1.5ex]
\mbox{}\hfill                                      
$
{\frak p}_j\,{\rm arctg}\,\dfrac{\widetilde{\nu}x}{{\stackrel{*}{\nu}}x}\,=\,\widetilde{\lambda}{}^j$
\ \ for all $x\in {\mathscr X},
\quad   j=1,\ldots, m,
\quad  {\mathscr X}\subset \{x\colon {\stackrel{*}{\nu}}x\ne 0\}.
\hfill 
$
}
\\[2ex]
\indent
{\bf 2.1.2. Autonomous first integrals}.
The proof of the following assertions is similar to those of Theorem 1.1 and Corollary 1.1.
\vspace{0.5ex}

{\bf Theorem 2.1}.
{\it 
Suppose  $\nu^{k}$ are real common eigenvectors of  the matrices  
$A_{j}$ cor\-res\-pon\-ding to the eigenvalues 
\vspace{0.35ex}
$\lambda_{k}^{j},\ j=1,\ldots, m, \ k=1,\ldots, m+1,$ respectively.
Then the system of equations in total differentials {\rm (2.1)} has the autonomous first integral
\\[1.5ex]
\mbox{}\hfill                                      
$
\displaystyle
F\colon x\to
\prod\limits_{k=1}^{m+1}\bigl|\nu^{k} x\bigr|^{h_{k}}$ 
\ \ for all $x\in {\mathscr X},
\quad {\mathscr X}\subset {\rm D}(F),
\hfill 
$
\\[1.5ex]
where $h_{1},\ldots,h_{m+1}$ is a real nontrivial solution to the system  
$
\sum\limits_{k=1}^{m+1}\!\lambda_{k}^{j}h_{k} = 0,
\, j=1,\ldots, m.
$
}
\vspace{0.75ex}

{\bf Corollary 2.1}.
{\it
Let $\nu^{k}$ be real common eigenvectors of  the matrices  
$A_{j}$ cor\-res\-pon\-ding to the eigenvalues 
\vspace{0.35ex}
$\lambda_{k}^{j},\ j=1,\ldots, m, \ k=1,\ldots, m+1,$ respectively.
Then an autonomous first integral of the system of equations in total differentials {\rm (2.1)} is
the function 
\\[1.75ex]
\mbox{}\hfill
$
\displaystyle
F_{12\ldots m(m+1)}\colon x\to
\prod\limits_{k=1}^{m}
\bigl|\nu^{k} x\bigr|^{{}^{\scriptstyle {}-\triangle_{k}}}
\bigl|\nu^{m+1} x\bigr|^{{}^{\scriptstyle\triangle}}$
\ \ for all $x\in {\mathscr X},
\quad
{\mathscr X}\subset {\rm D}(F_{12\ldots m(m+1)}),
\hfill
$
\\[1.75ex]
where 
\vspace{0.35ex}
the determinants $\triangle_{k},\, k=1,\ldots, m$ 
are obtained by replacing the $k\!$-th column of the determinant $\triangle = \big|\lambda_{k}^{j}\big|$  by  
${\rm colon}\left(\lambda_{m+1}^{1},\ldots,\lambda_{m+1}^{m}\right),$ respectively.
}
\vspace{1ex}

As an example, the linear autonomous system of equations in total differentials
\\[1.5ex]
\mbox{}\hfill                                     
$
\begin{array}{ll}
dx_1 = {}-x_1\,dt_2, &
dx_2 = 2(x_3+x_4)\,dt_1 +  x_2\,dt_2,
\\[1.25ex]
dx_3 = x_2\,dt_1 + x_4\,dt_2,\quad &
dx_4 = x_2\,dt_1 + x_3\,dt_2
\end{array}
\hfill (2.3)
$
\\[1.75ex]
has the commuting matrices 
\\[1.75ex]
\mbox{}\hfill
$
A_1 =
\left\|
\begin{array}{cccc}
 0\, & 0\, & 0\, & 0
\\
 0\, & 0\, & 1\, & 1
\\
 0\, & 2\, & 0\, & 0
\\
 0\, & 2\, & 0\, & 0
\end{array}
\right\|
$
\quad and \quad
$
A_2 =
\left\|\!
\begin{array}{rccc}
 {}-1\, & 0\, &  0\, & 0
\\
    0\, & 1\, &  0\, & 0
\\
    0\, & 0\, &  0\, & 1
\\
    0\, & 0\, &  1\, & 0
\end{array}
\right\|.
\hfill
$
\\[1.75ex]
\indent
Therefore the system of equations in total differentials (2.3) is completely solvable. 
\vspace{0.35ex}

The matrices $A_1$ and $A_2$ have the eigenvalues 
\vspace{0.35ex}
$\lambda_1^1\!=\!-2,\, \lambda_2^1=\lambda_3^1=0,\, \lambda_4^1=2,\!$
and
$\lambda_1^2\!=\!1,$ $\lambda_2^2=\lambda_3^2=-1, \, \lambda_4^2=1$
corresponding to the eigenvectors 
\vspace{0.35ex}
$
\nu^{1}\! = (0, -1, 1, 1),\,
\nu^{2}\! = (1, 0, 0, 0),$\
$\nu^{3} = (0, 0, 1,{}-1),\
\nu^{4} = (0, 1, 1, 1),
$
respectively. 
\vspace{0.35ex}

The determinants 
\\[1.5ex]
\mbox{}\hfill
$
\triangle =
\left|\!\!
\begin{array}{rr}
-2 &  0
\\
1  & -1
\end{array}
\!\!\right| = 2,
\quad
\triangle_{11} =
\left|\!\!
\begin{array}{rr}
0 &  0
\\
-1  & -1
\end{array}
\!\!\right| = 0,
\quad
\triangle_{21} =
\left|\!\!
\begin{array}{rr}
-2 &  0
\\
1  & -1
\end{array}
\!\!\right| = 2,
\hfill
$
\\[2ex]
\mbox{}\hfill
$
\triangle_{12} =
\left|\!\!
\begin{array}{rr}
2 &  0
\\
1  & -1
\end{array}
\!\!\right| ={}-2,
\qquad
\triangle_{22} =
\left|\!\!
\begin{array}{rr}
-2 & 2
\\
1  & 1
\end{array}
\!\!\right| ={}-4.
\hfill
$
\\[1.5ex]
\indent
By Corollary 2.1 we have that the functionally independent functions 
\\[1.5ex]
\mbox{}\hfill                                    
$
F_{{}_{\scriptstyle 123}}\colon x\to \dfrac{(x_3-x_4)^2}{x_1^2}$
\ \ for all $x\in {\mathscr X}
\hfill (2.4)
$
\\[1ex]
and
\\[1ex]
\mbox{}\hfill                                    
$
F_{{}_{\scriptstyle 124}}\colon x\to
x_1^4\,(x_2^2-(x_3+x_4)^2\,)^2$
\ \ for all $x\in {\mathbb R}^4
\hfill (2.5)
$
\\[2ex]
are first integrals of the system (2.3), where a domain ${\mathscr X}\subset \{x\colon x_1\ne 0\}.$
\vspace{1.5ex}

{\bf Theorem 2.2.}
{\it
Let $\nu^{k} = {\stackrel{*}{\nu}}{}^{\,k}+\widetilde{\nu}{}^{\,k}\,i$ {\rm(}this set hasn't complex conjugate vectors{\rm)}
be common complex eigenvectors of  the matrices  
$A_{j}$ cor\-res\-pon\-ding to the eigenvalues 
\vspace{0.35ex}
$\lambda_{k}^{j} = {\stackrel{*}{\lambda}}{}_{k}^{j} +
\widetilde{\lambda}{}_{k}^{j}\,i$ $j=1,\ldots, m,\ k=1,\ldots, s,$ $s\leq (m+1)/2,$ respectively. 
\vspace{0.35ex}
Let $\nu^{\theta}$ be common real eigenvectors of  $A_{j}$ cor\-res\-pon\-ding to the
\vspace{0.35ex}
eigenvalues $\lambda_{\theta}^{j},\ j=1,\ldots, m,\ \theta=s+1,\ldots, m+1-s.$ 
Then the system of equations in total differentials {\rm (2.1)} has the autonomous first integral
\\[2ex]
\mbox{}\hfill                                       
$
\displaystyle
F\colon x\to
\prod\limits_{k=1}^{s}
(P_{k}(x))^{{\stackrel{*}{h}_{k}}}
\exp\Bigl({}-2\,\widetilde{h}_{k}\,\varphi_{k}(x)\Bigr)
\prod\limits_{\theta=s+1}^{m+1-s}
\bigl|\nu^{\theta} x\bigr|^{h_{\theta}}$
\ \ \ for all $x\in {\mathscr X},
\hfill (2.6)
$
\\[2ex]
where a domain ${\mathscr X}\subset {\rm D}(F),$ the functions 
\\[1.5ex]
\mbox{}\hfill
$
P_{k}\colon x\to
\bigl({\stackrel{*}{\nu}}{}^{\,k} x\bigr)^{2} +
\bigl(\widetilde{\nu}{}^{\,k} x\bigr)^{2}\!$
for all $x\in {\mathbb R}^n, 
\
\varphi_{k}\colon x\to
{\rm arctg}\,\dfrac{\widetilde{\nu}{}^{\,k} x}
{{\stackrel{*}{\nu}}{}^{\,k} x}$
for all $x\in {\mathscr X},
\,  k=1,\ldots, s,
\hfill
$
\\[1.5ex]
and $\stackrel{*}{h}_{k},\, \widetilde{h}_{k},\, k=1,\ldots, s,\
h_{\theta},\, \theta=s+1,\ldots, m+1-s$ is a real nontrivial solution to 
\\[1.5ex]
\mbox{}\hfill
$
\displaystyle 
2\, \sum\limits_{k=1}^{s}
\bigl(\,{\stackrel{*}{\lambda}}{}_{k}^{j}\,
{\stackrel{*}{h}}_{k} -
\widetilde{\lambda}{}_{k}^{j}\,
\widetilde{h}_{k}\bigr) \ +\ 
\sum\limits_{\theta=s+1}^{m+1-s}
\lambda_{\theta}^{j}h_{\theta}=0, 
\quad
j=1,\ldots, m.
\hfill
$
\\[2.25ex]
}
\indent
{\sl Proof}. Taking into account Property 2.1 and Lemma 2.1, we obtain
\\[1.75ex]
\mbox{}\hfill                                         
$
{\frak p}_{j}
\bigl( ({\stackrel{*}{\nu}}{}^{\,k}x)^2 +(\widetilde{\nu}{}^{\,k} x)^2\bigr) =
2\,{\stackrel{*}{\lambda}}{}_{k}^{j}\,
\bigl( ({\stackrel{*}{\nu}}{}^{\,k}x)^2 +(\widetilde{\nu}{}^{\,k} x)^2\bigr),
\ \ j=1,\ldots, m,
\quad  
k=1,\ldots, s,
\hfill
$
\\
\mbox{}\hfill (2.7)
\\
\mbox{}\hfill
$
{\frak p}_{j}\,\nu^{\theta} x\, =\,
\lambda_{\theta}^{j}\, \nu^{\theta} x$
\ \ for all $x\in {\mathbb R}^n,
\ \ \ j=1,\ldots, m,
\quad \theta=s+1,\ldots, m+1-s.
\hfill
$
\\[2.25ex]
\indent
The Lie derivative of (2.6) by virtue of (2.1) is equal to 
\\[2ex]
\mbox{}\hfill
$
\displaystyle
{\frak p}_{j} F(x) =
\biggl(\
\prod\limits_{k=1}^{s}
\bigl(P_{k}(x)\bigr)^{{\stackrel{*}{h}_{k}}-1}
\exp\Bigl({}-2\,\widetilde{h}_{k}\,\varphi_{k}(x)\Bigr)
\, \sum\limits_{k=1}^{s}{\stackrel{*}{h}}_{k}\,
\prod\limits_{l=1, l\ne k}^{s} P_{l}(x)\cdot 
{\frak p}_{j}P_{k}(x) \ +
\hfill
$
\\[2.75ex]
\mbox{}\hfill
$
\displaystyle
+ \ \prod\limits_{k=1}^{s}
\bigl(P_{k}(x)\bigr)^{{\stackrel{*}{h}}_{k}}
\exp\Bigl({}-2\,\widetilde{h}_{k}\,\varphi_{k}(x)\Bigr)\, 
\sum\limits_{k=1}^{s}
{\frak p}_{j}
\Bigl({}-2\,\widetilde{h}_{k}\,\varphi_{k}(x)\Bigr)\biggr)
\prod\limits_{\theta=s+1}^{m+1-s}
\bigl|\nu^{\theta} x\bigr|^{h_{\theta}} \ +
\hfill
$
\\[2.75ex]
\mbox{}\hfill
$
\displaystyle
+\ \prod\limits_{k=1}^{s}
\bigl(P_{k}(x)\bigr)^{{\stackrel{*}{h}}_{k}}
\exp\Bigl({}-2\,\widetilde{h}_{k}\,\varphi_{k}(x)\Bigr)\,
\prod\limits_{\theta=s+1}^{m+1-s}
\bigl|\nu^{\theta} x\bigr|^{h_{\theta}-1}\ \cdot
\hfill
$
\\[2.75ex]
\mbox{}\hfill
$
\displaystyle
\cdot\ 
\sum\limits_{\theta=s+1}^{m+1-s}
\mbox{sgn}\,\bigl(\nu^{\theta} x\bigr)\,h_{\theta}\,
\prod\limits_{l=s+1,l\ne \theta}^{m+1-s}
\bigl|\nu^{l} x\bigr|\cdot 
{\frak p}_{j}\bigl(\nu^{\theta} x\bigr)$
\ \ for all $x\in {\mathscr X},
\quad j=1,\ldots, m.
\hfill
$
\\[2ex]
\indent
Using Property 2.2 and (2.7), we get
\\[2ex]
\mbox{}\hfill
$
\displaystyle
{\frak p}_{j}F(x) =
\biggl(\,
\sum\limits_{k=1}^{s}
2\bigl({\stackrel{*}{\lambda}}{}_{k}^{j}\, {\stackrel{*}{h}}_{k} -
\widetilde{\lambda}{}_{k}^{j}\,
\widetilde{h}_{k}\bigr)   +
\sum\limits_{\theta=s+1}^{m+1-s}\lambda_{\theta}^{j}h_{\theta}
\biggr) F(x)$
\ for all $x\in {\mathscr X},
\quad 
j=1,\ldots, m.
\hfill
$
\\[2ex]
\indent
If 
\vspace{0.75ex}
$
2\sum\limits_{k=1}^{s}
\bigl({\stackrel{*}{\lambda}}{}_{k}^{j}\, {\stackrel{*}{h}}_{k} -
\widetilde{\lambda}{}_{k}^{j}\, \widetilde{h}_{k}\bigr)  +
\sum\limits_{\theta=s+1}^{m+1-s}
\lambda_{\theta}^{j}h_{\theta} = 0,
\ j=1,\ldots, m,$ 
then the function (2.6) is an autonomous first integral of 
the system of equations in total differentials {\rm (2.1)}.\ \k

\newpage

For example, the completely solvable linear system of total differential equations 
\\[1.75ex]
\mbox{}\hfill                                         
$
dx_1 = x_1\,dt_1 + x_2\,dt_2,
\quad \
dx_2 = x_2\,dt_1 - x_1\,dt_2,
\quad \
dx_3 = x_3\,dt_1 - x_3\,dt_2
\hfill (2.8)
$
\\[1.75ex]
has the eigenvalues $\lambda_1^1 =\lambda_2^1 =\lambda_3^1 = 1,\  
\vspace{0.35ex}
\lambda_1^2 = -\,i, \ \lambda_2^2 = i,\ \lambda_3^2 = -\,1
$
cor\-res\-pon\-ding to the eigenvectors
$
\nu^{1}\! = (1, i, 0),
\, \nu^{2}\! =(1,-i, 0),
\, \nu^{3}\! = (0, 0, 1),
$
\vspace{0.35ex}
respectively. According to Theorem~2.2, we can construct the autonomous first integral of the system (2.8):
\\[1.75ex]
\mbox{}\hfill                                         
$
F\colon x\to
\dfrac{x_{1}^{2} + x_{2}^{2}}{x_{3}^{2}}\
\exp\Bigl( 2\,{\rm arctg}\, \dfrac{x_2}{x_1}\, \Bigr)$
\ \ for all $x\in {\mathscr  X},
\quad
{\mathscr  X}\subset \{ x\colon x_1\ne 0,\ x_3\ne 0\}.
\hfill (2.9)
$
\\[2.25ex]
\indent
{\bf Theorem 2.3}.
{\it
Suppose $\nu^{\tau} ={\stackrel{*}{\nu}}{}^{\,\tau}+
\widetilde{\nu}{}^{\,\tau}\,i, \
\nu^{s+\tau} ={\stackrel{*}{\nu}}{}^{\,\tau} -
\widetilde{\nu}{}^{\,\tau}\,i,
\ \tau=1,\ldots, s,\ s\leq m/2,$ and
$\nu^{2s+1} ={\stackrel{*}{\nu}}{}^{\,2s+1} +
\widetilde{\nu}{}^{\,2s+1}\,i$ are 
common complex eigenvectors of  the matrices  
$A_{j}$ cor\-res\-pon\-ding to the eigenvalues 
\vspace{0.35ex}
$
\lambda_{\tau}^{j}\! =\!{\stackrel{*}{\lambda}}{}_{\tau}^{j}\, +
\widetilde{\lambda}{}_{\tau}^{j}\,i,\
\lambda_{s+\tau}^{j}\! =\!
{\stackrel{*}{\lambda}}{}_{\tau}^{j}\,-
\widetilde{\lambda}{}_{\tau}^{j}\,i, \,
\tau\!=\!1,\ldots, s,\
\lambda_{2s+1}^{j}=
{\stackrel{*}{\lambda}}{}_{2s+1}^{j} +
\widetilde{\lambda}{}_{2s+1}^{j}\,i,$
$j=1,\ldots, m,$ respectively.
Let $\nu^{\theta}$ be common real eigenvectors of  $A_{j}$ cor\-res\-pon\-ding to the
\vspace{0.35ex}
eigenvalues $\lambda_{\theta}^{j},\ j=1,\ldots, m,\ \theta=2s+2,\ldots, m+1,$ respectively. 
Then the system of equations in total differentials {\rm (2.1)} has the autonomous first integrals
\\[2ex]
\mbox{}\hfill                                  
$
\displaystyle
F_{1}\colon  x\to
\prod\limits_{k=1}^{s}
\bigl(P_{k}(x)\bigr)^{{\stackrel{*}{h}}_{k}+{\stackrel{*}{h}}_{s+k}}
\exp\Bigl({}-2\,\bigl(\widetilde{h}_{k}-
\widetilde{h}_{s+k}\bigr)\varphi_{k}(x)\Bigr) \cdot
\hfill
$
\\
\mbox{}\hfill {\rm (2.10)}
\\
\mbox{}\hfill
$
\displaystyle
\cdot
\bigl(P_{2s+1}(x)\bigr)^{{\stackrel{*}{h}}_{2s+1}}
\exp\Bigl({}-2\,\widetilde{h}_{2s+1}\,\varphi_{2s+1}(x)\Bigr)
\prod\limits_{\theta=2s+2}^{m+1}
\bigl(\nu^{\theta} x\bigr)^{2\,{\stackrel{*}{h}}_{\theta}}$
\ \ for all $x\in {\mathscr X}
\hfill
$
\\[1.5ex]
and
\\[1.5ex]
\mbox{}\hfill                                        
$
\displaystyle
F_{2}\colon  x\to
\prod\limits_{k=1}^{s}
\bigl(P_{k}(x)\bigr)^{
\widetilde{h}_{k}+\widetilde{h}_{s+k}}
\exp\Bigl(2\,\bigl({\stackrel{*}{h}}_{k} -
{\stackrel{*}{h}}_{s+k}\bigr)\varphi_{k}(x) \Bigr) \cdot
\hfill
$
\\
\mbox{}\hfill {\rm (2.11)}
\\
\mbox{}\hfill
$
\displaystyle
\cdot
\bigl(P_{2s+1}(x)\bigr)^{\widetilde{h}_{2s+1}}
\exp\Bigl(2\,{\stackrel{*}{h}}_{2s+1}\,\varphi_{2s+1}(x)\Bigr)
\prod\limits_{\theta=2s+2}^{m+1}
\bigl(\nu^{\theta} x\bigr)^{2\,\widetilde{h}_{\theta}}$
\ \ for all $x\in {\mathscr X},
\hfill
$
\\[2ex]
where a domain ${\mathscr X}\subset {\rm D}(F_1)\cap{\rm D}(F_2),$ the functions 
\\[2ex]
\mbox{}\hfill
$
P_{k}\colon x\to
({\stackrel{*}{\nu}}{}^{\,k} x)^{2}  +
(\widetilde{\nu}{}^{\,k} x)^{2}$
\ \ for all $x\in {\mathbb R}^n, 
\quad
k=1,\ldots, s, 2s+1,
\hfill
$
\\[2.5ex]
\mbox{}\hfill
$
\varphi_{k}\colon x\to
{\rm arctg}\,\dfrac{\widetilde{\nu}{}^{\,k} x}{{\stackrel{*}{\nu}}{}^{\,k} x}$
\ \ \ for all $x\in {\mathscr X},
\quad
k=1,\ldots, s, 2s+1,
\hfill
$
\\[2ex]
and 
\vspace{0.75ex}
$h_{k} = {\stackrel{*}{h}}_{k} + \widetilde{h}_{k}\,i,\, k=1,\ldots, m+1$ is a nontrivial solution to 
$\sum\limits_{k=1}^{m+1}{\lambda}_{k}^{j}\,h_{k}=0,\,  j=1,\ldots, m.$
}

{\sl Proof.}
We form two complex-valued functions
\\[2ex]
\mbox{}\hfill
$
\displaystyle
{\stackrel{*}{F}}\colon x\to
\prod\limits_{k=1}^{2s}\bigl(\nu^{k} x\bigl)^{h_{k}}
\bigl(\nu^{2s+1} x\bigr)^{h_{2s+1}}
\prod\limits_{\theta=2s+2}^{m+1}
\bigl(\nu^{\theta} x\bigr)^{h_{\theta}}$
\ \ for all $x\in {\mathscr X}
\hfill
$
\\[1.5ex]
and
\\[1.5ex]
\mbox{}\hfill
$
\displaystyle
\stackrel{**}{F}\colon x\to
\prod\limits_{k=1}^{2s}\bigl( \nu^{k} x\bigr)^{l_{k}}
\bigl(\overline{\nu^{2s+1}} x\bigr)^{l_{2s+1}}
\prod\limits_{\theta=2s+2}^{m+1}
\bigl(\nu^{\theta} x\bigr)^{l_{\theta}}$
\ \ for all $x\in {\mathscr X},
\hfill
$
\\[2ex]
where a domain ${\mathscr X}\subset {\mathbb R}^n$ and $h_{k},\ l_{k},\, k=1,\ldots, m+1$ are complex numbers.
\vspace{0.35ex}
 
Using Lemma 2.1, we get
\\[2ex]
\mbox{}\hfill
$
\displaystyle
{\frak p}_{j}\,{\stackrel{*}{F}}(x) =
\sum\limits_{k=1}^{m+1}
\lambda_{k}^{j}\,h_{k}\,{\stackrel{*}{F}}(x)$
\ \ for all $x\in {\mathscr X},
\quad  j=1,\ldots, m,
\hfill
$
\\[2ex]
\mbox{}\hfill
$
\displaystyle
{\frak p}_{j}\,{\stackrel{**}{F}}(x) =
\biggl(\ \sum\limits_{k=1}^{2s}
\lambda_{k}^{j}\,l_{k} +
\overline{\lambda_{2s+1}^{j}}\, l_{2s+1} +
\sum\limits_{\theta=2s+2}^{m+1}
\lambda_{\theta}^{j}\,l_{\theta}
\biggr){\stackrel{**}{F}}(x)$
\ for all $x\in {\mathscr X},
\quad  j=1,\ldots, m.
\hfill
$
\\[2ex]
\indent
Let $h_{k}\! =\!{\stackrel{*}{h}}_{k}+
\widetilde{h}_{k}\,i,\, k\!=\!1,\ldots, m+1$ be a nontrivial solution to 
$\sum\limits_{k=1}^{m+1}\lambda_{k}^{j}h_{k} =  0, \, j\!=\!1,\ldots, m.$
Then $l_{k}={\stackrel{*}{h}}_{s+k} -\widetilde{h}_{s+k}\,i, \
\vspace{0.5ex}
l_{s+k} ={\stackrel{*}{h}}_{k} - \widetilde{h}_{k}\,i,\,
k\!=\!1,\ldots, s, \
l_{2s+1} ={\stackrel{*}{h}}_{2s+1} - \widetilde{h}_{2s+1}\,i,\ 
l_{\theta} ={\stackrel{*}{h}}_{\theta} - \widetilde{h}_{\theta}\,i,$
$\theta=2s+2,\ldots, m+1$ is a solution to 
$
\sum\limits_{k=1}^{2s}\lambda_{k}^{j}\,l_{k} +
\overline{\lambda_{2s+1}^{j}}\,l_{2s+1} +
\sum\limits_{\theta=2s+2}^{m+1}
\lambda_{\theta}^{j}\,l_{\theta} = 0,
\, j\!=\!1,\ldots, m$
and the functions ${\stackrel{*}{F}}\colon {\mathscr X}\to {\mathbb C},\
{\stackrel{**}{F}}\colon {\mathscr X}\to {\mathbb C}$ are first integrals of the system (2.1).
\vspace{0.5ex}

Since $F_1 = {\stackrel{*}{F}}\,{\stackrel{**}{F}}$ and
\vspace{0.35ex}
$F_{2} = \bigl({\stackrel{**}{F}}/{\stackrel{*}{F}}\bigr)^{i},$ we see that 
the functions (2.10) and (211) are autonomous first integrals of the system of equations in total differentials (2.1). \ \k
\vspace{1ex}

In particular, the completely solvable linear system of total differential equations 
\\[1.75ex]
\mbox{}\hfill                                         
$
dx_1 = x_1\,dt_1 + x_3\,dt_2,
\qquad
dx_2 = {}-x_2\,dt_1 + x_4\,dt_2,
\hfill
$
\\[-0.25ex]
\mbox{}\hfill (2.12)
\\[-0.25ex]
\mbox{}\hfill
$
dx_3 = x_3\,dt_1 - x_1\,dt_2,
\qquad
dx_4 = {}-x_4\,dt_1 - x_2\,dt_2
\hfill
$
\\[2ex]
has the eigenvalues 
\vspace{0.35ex}
$
\lambda_1^1=\lambda_2^1={}-1,\ 
\lambda_3^1=\lambda_4^1=1,\ 
\lambda_1^2=\lambda_3^2={}-i,\ \lambda_2^2=\lambda_4^2=i$
cor\-res\-pon\-ding to the linearly independent eigenvectors
\vspace{0.35ex}
$\nu^{1} = (0, {}-i, 0, 1),
\ \nu^{2} = (0, i,0, 1),
\ \nu^{3} = ({}-i,0, 1, 0),$ and
$\nu^{4} = (i, 0,1, 0),$ respectively.
The functions (by Theorem 2.3)
\\[2ex]
\mbox{}\hfill                                         
$
F_1\colon x\to \dfrac{x_1 x_2 + x_3 x_4}{x_1 x_4 - x_2 x_3}$
\ for all $x\in {\mathscr X},
$
\ 
$
F_2\colon x\to (x_1^2+ x_3^2)(x_2^2+ x_4^2)$
for all $x\in {\mathbb R}^4
\hfill (2.13)
$
\\[2ex]
are first integrals of the system (2.12), 
where a domain ${\mathscr X}\subset \{x\colon x_1 x_4 - x_2 x_3\ne 0\}\subset {\mathbb R}^4.$ 
\vspace{1ex}

From the entire set of ordinary differential systems induced by the completely solvable 
system of  equations in total differentials (2.1), we extract system
\\[1.75ex]                                  
\mbox{}\hfill
$
dx=A^{\zeta}(x)\,dt_{\zeta},
\quad
A^{\zeta}(x)=\mbox{colon}(a_{1\zeta}(x),\ldots,a_{n\zeta}(x))$
for all $x\in {\mathbb R}^n,
\hfill (2.1.\zeta)
$
\\[1.75ex]
such that the matrix $A_{\zeta}$ has the smallest number of elementary divisors. 
\vspace{0.75ex}

{\bf Definition 2.1}.
{\it
Let $\nu^{0l}$ be an eigenvector of the matrix $A_{\zeta}$ corresponding to the eigen\-va\-lue $\lambda^{\zeta}_{l}$ with 
elementary divisor of multiplicity $s_l.$
A non-zero vector $\nu^{\theta l}\in {\mathbb C}^n$ is called a
\textit{\textbf{generalized eigenvector of order}} {\boldmath $\theta$} for $\lambda^{\zeta}_{l}$ if and only if
\\[1.5ex]                                                      
\mbox{}\hfill
$
(A_{\zeta}-\lambda^{\zeta}_{l} E)\,\nu^{\theta l}=\theta  \cdot \nu^{\theta-1,l},
\quad 
\theta =1,\ldots, s_l-1,
\hfill (2.14)
$
\\[1.5ex]
where $E$ is the $n\times n$ identity matrix.
}
\vspace{0.5ex}

Using Lemma 2.1 and (2.14), we obtain
\\[1.5ex]
\mbox{}\hfill                                                      
$
{\frak p}_{\zeta}\,\nu^{0l}x =
\lambda^{\zeta}_{l}\,\nu^{0l}x,
\ \ \ \ 
{\frak p}_{\zeta}\,\nu^{\,\theta l}x =
\lambda^{\zeta}_{l}\,\nu^{\,\theta l}x+\theta\,\nu^{\,\theta -1,l}x$
\ for all $x\in {\mathbb R}^n,
\ \ \ \theta=1,\ldots, s_l-1.
$
\hfill (2.15)
\\[1.5ex]
\indent
The following lemmas are needed for the sequel.
\vspace{0.75ex}

{\bf Lemma 2.2}.                                           
{\it
Let $\nu^{\,0l}$ be a common eigenvector of the matrices $A_j$ corresponding to the eigenvalues
$\lambda^j_{l}, \ j=1,\ldots,m,$ respectively. Let $\nu^{\,\theta l},\ \theta =1,\ldots, s_l-1$ be generalized eigenvectors
of the matrix $A_\zeta$ 
\vspace{0.35ex}
corresponding to the eigenvalue $\lambda^{\zeta}_{l}$ with elementary divisor of multiplicity 
$s_l\ (s_{l}\geq 2).$ If the system $(2.1.\zeta)$ hasn't the first integrals 
\\[2ex]
\mbox{}\hfill                               
$
\displaystyle
F_{j\theta l}^{\,\zeta}\colon x\to
{\frak p}_j\, v_{\theta l}^{\zeta}(x)$ 
\ for all $x\in {\mathscr X}, 
\ \ j=1,\ldots, m, \ \, j\ne\zeta,
\ \ \  \theta =1,\ldots, s_l-1,
\hfill (2.16)
$
\\[2ex]
then 
\\[1.25ex]
\mbox{}\hfill                             
$
{\frak p}_{\zeta}\,v_{\theta l}^{\zeta}(x)=
\left[\!\!
\begin{array}{lll}
1\! & \text{for all}\ \ x\in {\mathscr X}, & \theta =1,
\\[1ex]
0\! & \text{for all}\ \ x\in {\mathscr X} & \theta =2,\ldots, s_{l}-1,
\end{array}
\right.
\hfill
$
\\[2ex]
\mbox{}\hfill                               
$
{\frak p}_{j}\,v_{\theta l}^{\zeta}(x)=\mu_{\theta l}^{j\zeta}={\rm const}$
\ for all $x\in {\mathscr X},
\ \ j=1,\ldots, m,\ \ j\ne \zeta,
\ \ \  \theta =1,\ldots, s_l-1,
\hfill
$ 
\\[2ex]
where $v_{\theta l}^{\zeta}\colon {\mathscr X}\to {\mathbb R},\  \theta =1,\ldots, s_l-1$ is a solution to the system 
\\[1.5ex]
\mbox{}\hfill                             
$
\nu^{\,\theta l}x=
{\displaystyle \sum\limits_{\delta=1}^{\theta} }
\binom{\theta -1}{\delta-1}v_{\delta l}^{\zeta}(x)\cdot \nu^{\,\theta-\delta,l}x,
\quad 
\theta=1,\ldots, s_l-1,
\quad
{\mathscr X}\subset \{x\colon \nu^{0l}x\ne 0\}.
\hfill (2.17)
$
\\[2ex]
}
\indent
{\bf Lemma 2.3}
{\it
Under the conditions of Lemma {\rm 2.2}, we have
\\[1.5ex]
\mbox{}\hfill                                                                 
$
{\frak p}_{j}\,\nu^{\,\theta l}x=
{\displaystyle \sum\limits_{\delta=0}^{\theta}}
\binom{\theta}{\delta}\,\mu_{\delta l}^{ j\zeta}\cdot\nu^{\,\theta-\delta, l}x$
\ \ for all $x\in {\mathscr X},
\ \ \ j=1,\ldots, m, 
\ \ \theta=1,\ldots, s_{l}-1,
\hfill
$
\\[1.75ex]
where $\mu_{0l}^{j\zeta}=\lambda^j_{l}, \ \mu_{\theta l}^{j\zeta}={\frak p}_{j}\,v_{\theta l}^{\zeta}(x),
\ \theta=1,\ldots, s_{l}-1, \ j=1,\ldots, m.$ 
}
\vspace{1.5ex}

{\bf Theorem 2.4}.
{\it
\vspace{0.5ex}
Let the assumptions of Lemma {\rm 2.2} with $l=1,\ldots r\ \Bigl(\, \sum\limits_{l=1}^{r}s_{l}\geq m+1\Bigr)$ hold.
Then the completely solvable system {\rm (2.1)} has the autonomous first integral
\\[1.5ex]
\mbox{}\hfill
$
\displaystyle
F\colon x\to
\prod\limits_{\xi=1}^{k}\bigl(\nu^{0\xi} x\bigr)^{h_{0 \xi}}
\exp\sum\limits_{q=1}^{\varepsilon_{\xi}}\,h_{q\xi} v_{q\xi}^{\zeta}(x)$
\ \ for all $x\in {\mathscr X},
\quad {\mathscr X}\subset {\rm D}(F),
\hfill
$
\\[1.5ex]
where $\sum\limits_{\xi=1}^{k}\varepsilon_{\xi}=m-k+1, \,
\varepsilon_{\xi}\leq s_{\xi}-1, \, \xi=1,\ldots, k,\, k\leq r,$ and 
$h_{q\xi},\, q=0,\ldots,\varepsilon_{\xi}, \, \xi=1,\ldots, k$ is a nontrivial solution to
the linear homogeneous algebraic system of equations 
\\[1.5ex]
\mbox{}\hfill
$
\displaystyle
\sum\limits_{\xi=1}^{k}
\bigl(\lambda_{\xi}^{j}\,h_{0\xi}
+ \sum\limits_{q=1}^{\varepsilon_{\xi}}
\mu_{q\xi}^{j\zeta}\,h_{q\xi}\big) = 0,
\ \ j=1,\ldots, m.
\hfill
$
\\[2ex]
}
\indent
The completely solvable linear homogeneous system of total differential equations 
\\[1.5ex]
\mbox{}\hfill                                       
$
\begin{array}{l}
dx_1 = x_2\,dt_1 + (2x_1-x_3)\,dt_2,
\\[1ex]
dx_2 = (2x_2-x_3-x_4)\,dt_1 + ({}-x_1+2x_2+x_4)\,dt_2,
\\[1ex]
dx_3 = (x_1-x_4)\,dt_1 + ({}-x_1+3x_3+x_4)\,dt_2,
\\[1ex]
dx_4 = ({}-x_1+2x_3+2x_4)\,dt_1 + (x_2-3x_3+x_4)\,dt_2
\end{array}
\hfill (2.18)
$
\\[1.5ex]
has the eigenvalue $\lambda_1^1=1$ with elementary divisor $(\lambda^1-1)^4$ 
corresponding to the  eigenvector $\nu^{0}=({}-1,1,{}-1,0)$ and to the generalized eigenvectors 
$\nu^{1}=(1,0,{}-1,{}-1),\ \nu^{2}=(1,{}-1,3,0),$ $\nu^{3}=({}-3,0,9,9).$  The functions (see (2.17))
\\[1.5ex]
\mbox{}\hfill                                                                      
$
v_{11}^1\colon x\to
\dfrac{x_1-x_3-x_4}{{}-x_1+x_2-x_3}$
\ \ for all $x\in {\mathscr X},
\hfill
$
\\[2ex]
\mbox{}\hfill
$
v_{21}^1\colon x\to
\dfrac{({}-x_1+x_2-x_3)(x_1-x_2+3x_3)-(x_1-x_3-x_4)^2}
{({}-x_1+x_2-x_3)^{2}}$
\  \ for all $x\in {\mathscr X},
\hfill (2.19)
$
\\[2ex]
\mbox{}\hfill
$
v_{31}^1\colon x\to
\dfrac{1}{({}-x_1+x_2-x_3)^{3}}\,
\bigl(({}-3x_1+9x_3+9x_4)({}-x_1+x_2-x_3)^2\, -
\hfill
$
\\[1.5ex]
\mbox{}\hfill
$
-\,3({}-x_1+x_2-x_3)(x_1-x_3-x_4)(x_1-x_2+3x_3) +
2(x_1-x_3-x_4)^3\,\bigr)$
\ for all $x\in {\mathscr X},
\hfill
$
\\[1.75ex]
where a domain ${\mathscr X}\subset \{ x\colon x_1-x_2+x_3\ne 0\}.$ 

\newpage

Autonomous first integrals of the system (2.18) are the functions (by Theorem 2.4)
\\[1.5ex]
\mbox{}\hfill                                       
$
F_{1}\colon x\to v_{21}^{1}(x),
\ \
F_{2}\colon x\to ({}-x_{1}+x_{2}-x_{3})^{2}
\exp\bigl( {}-2v_{11}^{1}(x)-v_{31}^{1}(x)\bigr)$
\ for all $x\in {\mathscr X}.
\hfill (2.20)
$
\\[2.5ex]
\indent
In the complex case, we shall have two logical possibilities:

1.\! Any function from the set 
$\!V\!=\!\{\nu^{0\xi}x, v_{q\xi}^{\zeta}(x)\colon\!  q\!=\!1,\ldots,\varepsilon_{\xi},\, \xi\!=\!1,\ldots, k,
\sum\limits_{\xi=1}^{k}\!\varepsilon_{\xi}\!=\!m\!-\!k\!+\!1\}\!\!$ 
has 
the complex conjugate function in the set $V.$
\vspace{0.25ex}

2.\! At least one function from the set $\!V\!$ has not the complex conjugate function 
\vspace{0.5ex}
in the $\!V.\!\!$
 
{\sl Case}\ 1. The completely solvable system (2.1) has the autonomous first integral
\\[1.75ex]
\mbox{}\hfill
$
\displaystyle
F\colon x\to
\prod\limits_{\xi=1}^{k_1}
\bigl(\bigl({\stackrel{*}{\nu}}{}^{0\xi} x\bigr)^{2} +
\bigl(\,\widetilde{\nu}\,{}^{0\xi} x\bigr)^{2}\,
\bigr)^{{\stackrel{*}{h}}_{0\xi}}\, 
\exp\Bigl({}-2\,\widetilde{h}_{0\xi}\ 
{\rm arctg}\,\dfrac{\widetilde{\nu}\,{}^{0\xi} x}
{{\stackrel{*}{\nu}}{}^{0\xi} x} \ +
\hfill
$
\\[2.25ex]
\mbox{}\hfill
$
\displaystyle
+ \ 2\sum\limits_{q=1}^{\varepsilon_{\xi}}
\bigl(\,{\stackrel{*}{h}}_{q\xi}\,{\stackrel{*}{v}}{}_{q\xi}^{\,\zeta}(x) -
\widetilde{h}_{q\xi}\,\widetilde{v}{}_{q\xi}^{\,\,\zeta}(x)
\bigr)\Bigr)
\prod\limits_{\theta=1}^{k_2}
\bigl|\nu^{0\theta} x\bigr|^{h_{0\theta}}
\exp\sum\limits_{q=1}^{\varepsilon_{\theta}}
h_{q\theta}\,v_{q\theta}^{\,\zeta}(x)$
\ for all $x\in {\mathscr X},
\ \ {\mathscr X}\subset {\rm D}(F),
\hfill
$
\\[1.75ex]
where 
\vspace{0.35ex}
$\stackrel{*}{h}_{q\xi},\ \widetilde{h}_{q\xi}, q=0,\ldots, \varepsilon_\xi,\ \xi=1,\ldots, k_1,\ 
h_{q\theta},\ q=0,\ldots, \varepsilon_\theta,\ \theta=1,\ldots, k_2$ is a real nontrivial solution to
the linear homogeneous algebraic system of equations 
\\[1.75ex]
\mbox{}\hfill
$
\displaystyle
2\sum\limits_{\xi=1}^{k_1}
\Bigl( \bigl(
{\stackrel{*}{\lambda}}{}_{\xi}^{j}\,{\stackrel{*}{h}}_{0\xi} -
\widetilde{\lambda}{}_{\xi}^{j}\,\widetilde{h}_{0\xi}
\bigr) +
\sum\limits_{q=1}^{\varepsilon_{\xi}}
\bigl(
{\stackrel{*}{\mu}}{}_{q\xi}^{\,j\zeta}\,{\stackrel{*}{h}}_{q\xi} -
\widetilde{\mu}{}_{q\xi}^{\,j\zeta}\,\widetilde{h}_{q\xi}
\bigr)
\Bigr) +
\sum\limits_{\theta=1}^{k_2}
\Bigl(
\lambda_{\theta}^{j}\,h_{0\theta} +
\sum\limits_{q=1}^{\varepsilon_{\theta}}
\mu_{q\theta}^{j\zeta}\,h_{q\theta} \Bigr)\! =\! 0,
\ \ j=1,\ldots, m.
\hfill
$
\\[1.75ex]
\indent
Here $\nu^{0\xi}={\stackrel{*}{\nu}}{}^{0\xi} +\widetilde{\nu}\,{}^{0\xi}\,i$ are complex 
common eigenvectors of the matrices $A_j$ corresponding to the eigenvalues
\vspace{0.35ex}
$\lambda_{\xi}^{j} ={\stackrel{*}{\lambda}}{}_{\xi}^{j} +\widetilde{\lambda}{}_{\xi}^{j}\,i,\ 
j=1,\ldots, m,\ \xi=1,\ldots, k_1,$ respectively; $\nu^{0\theta}$ are real 
common eigenvectors of $A_j$ corresponding to the eigenvalues
\vspace{0.35ex}
$\lambda_{\theta}^{j}, \ j=1,\ldots, m,\ \theta=1,\ldots, k_2,$ respectively;
the functions $v_{q\xi}^{\,\zeta}\!=\!{\stackrel{*}{v}}{}_{q\xi}^{\,\zeta}\! +\!\widetilde{v}{}_{q\xi}^{\,\,\zeta}\,i, 
v_{q\theta}^{\,\zeta}\!$
is the solution to the system \!(2.17); the numbers 
\\[2ex]
\mbox{}\hfill
$
{\stackrel{*}{\mu}}{}_{q\xi}^{\,j\zeta} =
{\frak p}_{j}\,{\rm Re}\,v_{q\xi}^{\zeta}={\frak p}_{j}\,{\stackrel{*}{v}}{}_{q\xi}^{\,\zeta},
\ \ \ \
\widetilde{\mu}{}_{q\xi}^{\,j\zeta} =
{\frak p}_{j}\,{\rm Im}\,v_{q\xi}^{\zeta}={\frak p}_{j}\,\widetilde{v}{}_{q\xi}^{\,\,\zeta},
\ \ \ q=1,\ldots, \varepsilon_\xi,
\ \ \xi=1,\ldots, k_1,
\hfill
$
\\[2.5ex]
\mbox{}\hfill
$
{\mu}_{q\theta}^{j\zeta}={\frak p}_{j}\, v_{q\theta}^{\zeta},
\quad
q=1,\ldots,\varepsilon_\theta,
\ \ \theta=1,\ldots, k_2,
\ \ \ j=1,\ldots, m;
\hfill
$
\\[2ex]
the numbers $\varepsilon_{\xi},\ \varepsilon_{\theta}$ such that 
$
2\sum\limits_{\xi=1}^{k_1}\varepsilon_{\xi}
+ \sum\limits_{\theta=1}^{k_2}
\varepsilon_{\theta} = m-2k_{1}-k_{2}+1$
with
$
2k_1\!+\!k_2\leq r,\, \varepsilon_{\xi}\leq s_{\xi}\!-\!1,$ $\xi=1,\ldots, k_1,\,
\varepsilon_{\theta}\leq s_{\theta}\!-\!1,\, \theta=1,\ldots, k_2,$
\vspace{0.35ex}
where $k_1\!$ is a number of  complex common eigen\-vec\-tors {\rm(}this set hasn't complex conjugate vectors{\rm)} 
of  the matrices $A_j$  and  
$k_2$ is a number of real common eigenvectors of the matrices $A_j, \ j=1,\ldots, m.$
\vspace{1.5ex}

For example, the completely solvable system of total differential equations 
\\[1.75ex]                               
\mbox{}\qquad
$
dx_1=(3,-4,4,1,0,2)x\,dt_1+
(0,-4,2,1,-1,1)x\,dt_2
-(3,-2,4,3,0,2)x\,dt_3,
$
\\[1.5ex]
\mbox{}\qquad
$
dx_2={}-(1,-3,3,0,2,3)x\,dt_1 +
(1,3,0,0,1,-1)x\,dt_2 + (2,-3,3,3,-1,2)x\,dt_3,
\hfill (2.21)
$
\\[1.5ex]
\mbox{}\qquad
$
dx_3={}-(3,-5,5,1,2,4)x\,dt_1-
(0,-6,2,1,-2,1)xdt_2+(3,-3,5,4,0,2)x\,dt_3,
$
\\[1.5ex]
\mbox{}\qquad
$
dx_4= (3,-6,4,4,-1,5)x\,dt_1+ (2,-6,2,3,-4,2)x\,dt_2 -(3,-2,6,4,1,1)x\,dt_3,
\hfill
$
\\[1.5ex]
\mbox{}\qquad
$
dx_5=(5,-5,8,3,3,6)x\,dt_1 +
(1,-6,3,2,-2,2)x\,dt_2 - (3,-3,6,4,1,2)x\,dt_3,
$
\\[1.5ex]
\mbox{}\qquad
$
dx_6=-(2,-5,4,3,-1,2)x\,dt_1-(2,-4,3,3,-2,2)x\,dt_2
+(2,-1,4,2,1,0) x\,dt_3
$
\\[1.75ex]
has the eigenvalue 
\vspace{0.35ex}
$\lambda_{1}^{1}=1+2i$ with elementary divisor $\left(\lambda^1-1-2i\right)^3$  
corresponding to the eigenvector $\nu^{0}=(1, 0, 1+i, 1, i, 1)$ 
\vspace{0.35ex}
and to the generalized eigenvectors
$\nu^{1}=(1, 1+i, 0, 0, i, i),$ $\nu^{2}=(2+2i, 0, 2+2i, 0, 2i, 2i).$ The scalar functions of the vector argument
\\[1.5ex]
\mbox{}\hfill
$
{\stackrel{*}{v}}{}_{11}^{\,1}\colon x\to
\bigl( (x_1+x_2)(x_1+x_3+x_4+x_6)\,+\,(x_3+x_5)(x_2+x_5+x_6)\bigr) / P(x),
\hfill
$
\\[2ex]
\mbox{}\hfill
$
\widetilde{v}{}_{11}^{\,1}\colon x\to
\bigl( (x_1+x_3+x_4+x_6)(x_2+x_5+x_6)-(x_1+x_2)(x_3+x_5)\bigr) / P(x),
\hfill
$
\\[2ex]
\mbox{}\hfill
$
{\stackrel{*}{v}}{}_{21}^{\,1}\colon x\to
\bigl(
\bigl( (x_1+x_3+x_4+x_6)(x_2+x_5+x_6)-(x_1+x_2)(x_3+x_5)\bigr)^2 \ +
\hfill
$
\\[2ex]
\mbox{}\hfill
$
+\  2\,P(x)\bigl(
(x_1+x_3)(x_1+x_3+x_4+x_6)+(x_3+x_5)(x_1+x_3+x_5+x_6)
\bigr)\  -
\hfill
$
\\[1.75ex]
\mbox{}\hfill
$
- \, \bigl(
(x_1+x_2)(x_1+x_3+x_4+x_6)\,+\,(x_3+x_5)(x_2+x_5+x_6)
\bigr)^2 \bigr) / P(x),
\hfill (2.22)
$
\\[2ex]
\mbox{}\hfill
$
\widetilde{v}{}_{21}^{\,1}\colon x\to
2\bigl( P(x)
\bigl( (x_1+x_3+x_4+x_6)(x_1+x_3+x_5+x_6) - (x_1+x_3)(x_3+x_5) \bigr)\ +
\hfill
$
\\[1.75ex]
\mbox{}\hfill
$
 +\ \bigl((x_3+x_5)(x_1+x_2) -  (x_1+x_3+x_4+x_6)(x_2+x_5+x_6) \bigr)\, \cdot
\hfill
$
\\[1.75ex]
\mbox{}\hfill
$
\cdot
\bigl((x_1+x_2)(x_1+x_3+x_4+x_6)+(x_3+x_5)(x_2+x_5+x_6)
\bigr)\bigr)/ P(x)$
\ \ for all $x\in  {\mathscr X},
\hfill
$
\\[2ex]
where a domain ${\mathscr X}\subset \{ x\colon x_1+x_3+x_4+x_6\ne 0\},$ the polynomial 
\\[1.75ex]
\mbox{}\hfill                                        
$
P\colon x\to (x_1+x_3+x_4+x_6)^2+(x_3+x_5)^2$ 
\ \ for all $x\in {\mathbb R}^6.
\hfill (2.23)
$
\\[1.75ex]
\indent
Autonomous first integrals of the system (2.21) are the functions 
\\[1.75ex]
\mbox{}\hfill                                           
$
F_1\colon x\to P(x)\exp\bigl({}-4\varphi(x)
+6\,{\stackrel{*}{v}}{}_{11}^{\,1}(x)+2\,\widetilde{v}{}_{11}^{\,1}(x)\bigr)$
\ \ for all $x\in  {\mathscr X},
\hfill (2.24)
$
\\[2.25ex]
\mbox{}\hfill                                           
$
F_2\colon  x\to
P^2(x)\exp\bigl({}-2\varphi(x)+
{\stackrel{*}{v}}{}_{21}^{\,1}(x) - \widetilde{v}{}_{21}^{\,1}(x)\bigr)$
\ \ for all $x\in  {\mathscr X},
\hfill (2.25)
$
\\[2.25ex]
\mbox{}\hfill                                           
$
F_3\colon x\to
2\,\widetilde{v}{}_{11}^{\,1}(x) - 2\,{\stackrel{*}{v}}{}_{21}^{\,1}(x)
- \widetilde{v}{}_{21}^{\,1}(x)$
\ \ for all $x\in  {\mathscr X},
\hfill (2.26)
$
\\[1.5ex]
where
\\[1.25ex]
\mbox{}\hfill 
$
\varphi\colon x\to \ 
 {\rm arctg}\,\dfrac{x_3+x_5}{x_1+x_3+x_4+x_6}$
\ \ for all $x\in {\mathscr X}.
\hfill (2.27)
$
\\[3ex]
\indent
{\sl Case}\ 2A. 
A common complex eigenvector of the matrices $\!A_{j}, j\!=\!1,\ldots, m\!$ hasn't the com\-p\-lex conjugate vector.
\!First integrals of the completely solvable system (2.1) are the functions
\\[2ex]
\mbox{}\hfill
$
\displaystyle
F_{1}\colon  x\to
\prod\limits_{\xi=1}^{k_1}
\bigl(P_{\xi}(x)\bigr)^{{\stackrel{*}{h}}_{0\xi}+
{\stackrel{*}{h}}_{0,(k_1+\xi)}}
\exp\Bigl({}-2\bigl(\,\widetilde{h}_{0\xi}-
\widetilde{h}_{0,(k_1+\xi)}\bigr)\varphi_{\xi}(x)\ +
\hfill
$
\\[1.75ex]
\mbox{}\hfill
$
\displaystyle
+\ 2\sum\limits_{q=1}^{\varepsilon_{\xi}}
\Bigl(
\bigl(\,
{\stackrel{*}{h}}_{q\xi}+{\stackrel{*}{h}}_{q,(k_1+\xi)}
\bigr)\,{\stackrel{*}{v}}{}_{q\xi}^{\,\zeta}(x) +
\bigl(\,
\widetilde{h}_{q,(k_1+\xi)}-\widetilde{h}_{q\xi}
\bigr)\,\widetilde{v}{}_{q\xi}^{\,\zeta}(x)
\Bigr)\Bigr) \bigl(P_{2k_1+1}(x)\bigr)^{{\stackrel{*}{h}}_{0,(2k_1+1)}}\ 
\cdot
\hfill
$
\\[1.75ex]
\mbox{}\hfill
$
\displaystyle
\cdot\,
\exp\Bigl({}-2\,
\widetilde{h}_{0,(2k_1+1)}\,\varphi_{2k_1+1}(x)\Bigr)\,\cdot
\prod\limits_{\theta=1}^{k_2}
\bigl(\nu^{0\theta} x\bigr)^{2\,{\stackrel{*}{h}}_{0\theta}}\,
\exp\Bigl(\, 2\sum\limits_{q=1}^{\varepsilon_{\theta}}
{\stackrel{*}{h}}_{q\theta}\,v_{q\theta}^{\,\zeta}(x)\Bigr)$
\ \ for all $x\in {\mathscr X},
\hfill
$
\\[2.75ex]
\mbox{}\hfill
$
\displaystyle
F_{2}\colon x\to
\prod\limits_{\xi=1}^{k_1}
\bigl(P_{\xi}(x)\bigr)^{\widetilde{h}_{0\xi}+\widetilde{h}_{0,(k_1+\xi)}}
\exp\Big( 2\bigl(\, {\stackrel{*}{h}}_{0\xi}-
{\stackrel{*}{h}}_{0,(k_1+\xi)}\bigr)\,\varphi_{\xi}(x) \ +
\hfill
$
\\[1.75ex]
\mbox{}\hfill
$
\displaystyle
+ \ 2\sum\limits_{q=1}^{\varepsilon_{\xi}}
\Bigl(
\bigl(\,\widetilde{h}_{q\xi}+\widetilde{h}_{q,(k_1+\xi)}\bigr)\,{\stackrel{*}{v}}{}_{q\xi}^{\,\zeta}(x) +
\bigl(\,
{\stackrel{*}{h}}_{q\xi}-{\stackrel{*}{h}}_{q,(k_1+\xi)}
\bigr)\,\widetilde{v}{}_{q\xi}^{\,\zeta}(x)
\Bigr)
\Bigr) \bigl( P_{2k_1+1}(x)\bigr)^{\widetilde{h}_{0,(2k_1+1)}}\ \cdot
\hfill
$
\\[1.5ex]
\mbox{}\hfill
$
\displaystyle
\cdot
\exp\Bigl(\,2\,{\stackrel{*}{h}}_{0,(2k_1+1)}\,\varphi_{2k_1+1}(x)
\Bigr)\cdot
\prod\limits_{\theta=1}^{k_2}
\bigl( \nu^{0\theta} x \bigr)^
{2\,\widetilde{h}_{0\theta}}
\exp\Bigl(\, 2\sum\limits_{q=1}^{\varepsilon_{\theta}}\widetilde{h}_{q\theta}\,v_{q\theta}^{\,\zeta}(x)\Bigr)$
\ \ for all $x\in {\mathscr X},
\hfill
$
\\[2ex]
where a domain ${\mathscr X}\subset {\rm D}(F_1)\cap {\rm D}(F_2),$ the functions 
\\[1.75ex]
\mbox{}\hfill
$
P_{\xi}\colon x \to
({\stackrel{*}{\nu}}{}^{\,0\xi} x)^{2} +
(\,\widetilde{\nu}{}^{\,0\xi}x)^{2}$
\ \ for all $x\in {\mathbb R}^n,
\quad
\xi=1,\ldots, k_1, 2k_1+1,
\hfill
$
\\[2ex]
\mbox{}\hfill
$
\varphi_{\xi}\colon x\to\ 
{\rm arctg}\,\dfrac{\widetilde{\nu}{}^{\,0\xi} x}
{{\stackrel{*}{\nu}}{}^{\,0\xi} x}$
\ \ \, for all $x\in {\mathscr X},
\quad
\xi=1,\ldots, k_1, 2k_1+1,
\hfill
$
\\[1.75ex]
and $h_{q\xi}={\stackrel{*}{h}}_{q\xi}+\widetilde{h}_{q\xi}\,i,\ 
h_{q\theta}={\stackrel{*}{h}}_{q\theta}+\widetilde{h}_{q\theta}\,i$ is a nontrivial solution to the system 
\\[2ex]
\mbox{}\hfill
$
\displaystyle
\sum\limits_{\xi=1}^{2k_1}
\Bigl(
\lambda_{\xi}^{j}h_{0\xi} +
\sum\limits_{q=1}^{\varepsilon_{\xi}}
\mu_{q\xi}^{\,j\zeta}h_{q\xi}
\Bigr) +
\lambda_{2k_1+1}^{j}h_{0,(2k_1+1)} +
\sum\limits_{\theta=1}^{k_2}
\Bigl(
\lambda_{\theta}^{j}h_{0\theta} +\!\!
\sum\limits_{q=1}^{\varepsilon_{\theta}}
\mu_{q\theta}^{j\zeta}h_{q\theta}
\Bigr) = 0,
\ \
j=1,\ldots, m.
\hfill
$
\\[2ex]
\indent
Here $\nu^{0\xi}={\stackrel{*}{\nu}}{}^{\,0\xi}\, +\widetilde{\nu}{}^{\,0\xi}\,i, \
\nu^{0,(k_1+\xi)}=\overline{\nu^{0\xi}},\ 
\nu^{0,(2k_1+1)}={\stackrel{*}{\nu}}{}^{\,0,(2k_1+1)}\,+\widetilde{\nu}{}^{\,\, 0,(2k_1+1)}\,i$
are complex eigenvectors of the matrices $\!\!A_j\!\!$ corresponding to the eigenvalues
$\!\lambda_{\xi}^{j}\!\!=\!{\stackrel{*}{\lambda}}{}_{\xi}^{j}\! +\!\widetilde{\lambda}{}_{\xi}^{j}\,i, 
\lambda_{k_1+\xi}^{j}\!\!=\!\overline{\lambda_{\xi}^{j}},$
$\xi=1,\ldots, k_1,\ \lambda_{2k_1+1}^{j}={\stackrel{*}{\lambda}}{}_{2k_1+1}^{j}\,+
\widetilde{\lambda}{}_{2k_1+1}^{j}\,i,\ j=1,\ldots, m,$ respectively; 
\vspace{0.35ex}
$\nu^{0\theta}$ are real common eigenvectors of $A_j$ corresponding to the eigenvalues
\vspace{0.35ex}
$\lambda_{\theta}^{j},\ j=1,\ldots,m,\ \theta=1,\ldots, k_2,$ respectively;
the functions $v_{q\xi}^{\,\zeta}={\stackrel{*}{v}}{}_{q\xi}^{\,\zeta} +\widetilde{v}{}_{q\xi}^{\,\,\zeta}\,i,\ 
v_{q\theta}^{\,\zeta}$
is the solution to the system (2.17); the numbers 
\\[2ex]
\mbox{}\hfill
$
{\mu}_{q\xi}^{j\zeta}={\frak p}_{j}\,v_{q\xi}^{\,\zeta}(x),
\ \ \ \
{\stackrel{*}{\mu}}{}_{q\xi}^{\,j\zeta} =\, {\rm Re}\,{\mu}_{q\xi}^{j\zeta},
\ \ \
\widetilde{\mu}{}_{q\xi}^{\,j\zeta} =\,{\rm Im}\,{\mu}_{q\xi}^{j\zeta},
\ \ \
q=1,\ldots,\varepsilon_\xi,
\ \ 
\xi=1,\ldots, 2k_1,
\hfill
$
\\[2.5ex]
\mbox{}\hfill
$
\mu_{q\theta}^{j\zeta} ={\frak p}_{j}\,v_{q\theta}^{\,\zeta}(x),
\ \ \
q=1,\ldots,\varepsilon_\theta,
\ \ 
\theta=1,\ldots, k_2,
\ \ \ j=1,\ldots,m;
\hfill
$
\\[1.75ex]
the numbers $\varepsilon_{\xi},\ \varepsilon_{\theta}$ such that
$
2\sum\limits_{\xi=1}^{k_1} \varepsilon_{\xi} +
\sum\limits_{\theta=1}^{k_2} \varepsilon_{\theta} =m-2k_1-k_2$
with
$2k_1+1+k_2 \leq r,\ \varepsilon_{\xi}\leq s_{\xi}-1,$
$\xi=1,\ldots, k_1,\ \varepsilon_{\theta}\leq s_{\theta}-1,\, \theta=1,\ldots, k_2,$
where $k_1\!$ is a number of  complex common eigen\-vec\-tors {\rm(}this set hasn't complex conjugate vectors{\rm)} 
of  the matrices $A_j$  and  
$k_2$ is a number of real common eigenvectors of the matrices $A_j, \ j=1,\ldots, m.$
\vspace{1.5ex}

As an example, the completely solvable system of equations in total differentials
\\[1.75ex]                                          
$
dx_1 \!=\! (1,-2,2,0,1,1)xdt_1 + (0,2,0,0,1,1)xdt_2 +
(3,0,0,0,-1,-1)xdt_3 + (1,-2,4,0,2,2)xdt_4,
$
\\[1.5ex]
$
dx_2\! =\! (0,2,-2,0,-2,-2)xdt_1 - (1,3,0,0,1,1)xdt_2 +
(-1,2,0,0,1,1)xdt_3 - (2,1,4,0,4,4)xdt_4,
$
\\[1.5ex]
$
dx_3 = (0,3,-2,0,-2,-2)x\,dt_1 - (1,3,1,0,2,2)x\,dt_2\ +
$
\\[1.25ex]
\mbox{}\hfill
$
+ \ (-2,-1,2,0,1,1)x\,dt_3 - (3,-2,7,0,5,5)x\,dt_4,
\qquad  (2.28)
$
\\[1.75ex]
$\!\!
dx_4\! =\! x\bigl(\! (0,\!-4,0,2,\!-2,2)dt_1 \!+ (2,2,0,1,0,4)dt_2 +
(1,2,\!-2,1,\!-1,\!-1)dt_3 + (3,\!-4,10,2,7,7)dt_4\!\bigr),\!\!
$
\\[1.5ex]
$
dx_5 \!=\! (2,-3,4,2,2,4)xdt_1 + (3,3,2,2,1,4)xdt_2 +
(2,1,-1,0,0,-1)xdt_3 + (3,-2,9,0,7,5)xdt_4,
$
\\[1.5ex]
$\!
dx_6 \!=\! x\bigl(\!-(1,\!-3,2,2,\!-1,1)dt_1\! - (2,1,2,2,1,4)dt_2\! +
(\!-1,\!-1,1,0,1,2)dt_3\! - (\!-1,\!-4,5,0,4,2)dt_4\!\bigr)\!\!
$
\\[1.75ex]
has the eigenvalue 
\vspace{0.35ex}
$\lambda_{1}^{1}=1+i$ with elementary divisor $(\lambda^1-1-i)^2$  
corresponding to the eigen\-vec\-tor $\nu^{01}=(1, 1+i, 0, 0, i, i)$ and 
\vspace{0.35ex}
to the generalized eigenvector
$\nu^{11}=(1+i, 0, 1+i, 0, i, i)$
and the simple eigenvalue $\lambda_{2}^1=2i$ with common eigenvector $\nu^{02}=(1, 0, 1+i, 1,i, 1).$
\vspace{0.5ex}

The real-valued scalar functions
\\[2ex]
\mbox{}\ \ \ 
$
{\stackrel{*}{v}}{}_{11}^{\,1} \colon x\to
\dfrac{(x_1+x_2)(x_1+x_3)+(x_2+x_5+x_6)(x_1+x_3+x_5+x_6)}{P_1(x)}$
\ \ for all $x\in {\mathscr X},
\hfill
$
\\
\mbox{}\hfill (2.29)
\\
\mbox{}\ \ \ 
$
\widetilde{v}{}_{11}^{\,1} \colon x\to
\dfrac{(x_1+x_2)(x_1+x_3+x_5+x_6)-(x_1+x_3)(x_2+x_5+x_6)}{P_1(x)}$
\ \ for all $x\in {\mathscr X},
\hfill
$
\\[2ex]
where $P_1\colon x\to (x_1+x_2)^2+(x_2+x_5+x_6)^2$ for all $x\in {\mathbb R}^6.$ 
\vspace{0.5ex}

The autonomous general integral of the system (2.28) is the functions
\\[1.75ex]
\mbox{}\hfill                                         
$
F_1\colon x\to P_1(x)\bigl(P_2(x)\bigr)^2
\exp\bigl({}-10\varphi_1(x) + 8\,{\stackrel{*}{v}}{}_{11}^{\,1}(x) +
6\,\widetilde{v}{}_{11}^{\,1}(x)\bigr)$
\ \ for all $x\in {\mathscr X}
\hfill (2.30)
$
\\[1.5ex]
and
\\[1.5ex]
\mbox{}\hfill                                          
$
F_2\colon x\to
\bigl(P_1(x)\bigr)^3 \exp\bigl({}-10\varphi_1(x) -
4\varphi_2(x) + 12\,{\stackrel{*}{v}}{}_{11}^{\,1}(x) +
14\,\widetilde{v}{}_{11}^{\,1}(x)\bigr)$
\ \ for all $x\in {\mathscr X},
\hfill (2.31)
$
\\[2ex]
where a domain ${\mathscr X}\subset \{ x\colon x_1+x_2\ne 0,\ x_1+x_3+x_4+x_6\ne 0\},$ the scalar functions 
\\[2ex]
\mbox{}\hfill
$
P_2\colon x\to (x_1+x_3+x_4+x_6)^2+(x_3+x_5)^2$
\ \ for all $x\in {\mathbb R}^6,
\hfill
$
\\[2.25ex]
\mbox{}\hfill
$
\varphi_1\colon x\to \ 
{\rm arctg}\,\dfrac{x_2+x_5+x_6}{x_1+x_2}\,, 
\quad
\varphi_2\colon x\to \ 
{\rm arctg}\,\dfrac{x_3+x_5}{x_1+x_3+x_4+x_6}$
\ \ for all $x\in {\mathscr  X}.
\hfill
$
\\[3ex]
\indent
{\sl Case} 2B. 
\vspace{0.35ex}
A function $\!v_{l\gamma}^{\,\zeta},\, \gamma\!\in\! \{1,\ldots,k_1\},\, l\!\in\! \{1,\ldots,\varepsilon_{\gamma}\}\!$ 
hasn't the complex conjugate fun\-c\-tion.
Autonomous first integrals of the completely solvable system (2.1) are the functions
\\[2ex]
\mbox{}\hfill
$
\displaystyle
F_{1}\colon x\to
\prod\limits_{\xi=1}^{k_1}
\bigl(P_{\xi}(x)\bigr)^{{\stackrel{*}{h}}_{0\xi}+
{\stackrel{*}{h}}_{0,(k_1+\xi)}}
\exp\Bigl({}-2\bigl(\,\widetilde{h}_{0\xi}-
\widetilde{h}_{0,(k_1+\xi)}\bigr)\,\varphi_{\xi}(x) \ +
\hfill
$
\\[2ex]
\centerline{
$
\displaystyle
+\ 2\sum\limits_{q=1}^{\varepsilon_{\xi}}
(1-\delta_{ql}\delta_{\xi\gamma})\Bigl(
\bigl(\,
{\stackrel{*}{h}}_{q\xi}+{\stackrel{*}{h}}_{q,(k_1+\xi)}
\bigr)\,{\stackrel{*}{v}}{}_{q\xi}^{\,\zeta}(x) +
\bigl(\,
\widetilde{h}_{q,(k_1+\xi)}-\widetilde{h}_{q\xi}
\bigr)\,\widetilde{v}{}_{q\xi}^{\,\zeta}(x)
\Bigr) \ +
$
}
\\[2ex]
\mbox{}\hfill
$
\displaystyle
+ \ 2\,\Bigl(\,{\stackrel{*}{h}}_{l\gamma}\,
{\stackrel{*}{v}}{}_{l\gamma}^{\,\zeta}(x)  -
\widetilde{h}_{l\gamma}\,
\widetilde{v}{}_{l\gamma}^{\,\zeta}(x)\Bigr)
\Bigr)
\prod\limits_{\theta=1}^{k_2}
\bigl(\nu^{0\theta} x \bigr)^{2\,{\stackrel{*}{h}}_{0\theta}}
\exp\biggl( 2\sum\limits_{q=1}^{\varepsilon_{\theta}}
{\stackrel{*}{h}}_{q\theta}\, v_{q\theta}^{\,\zeta}(x)\biggr)$
\ \ for all $x\in {\mathscr X},
\hfill
$
\\[1.75ex]
and
\\[2ex]
\mbox{}\hfill
$
\displaystyle
F_{2}\colon x\to
\prod\limits_{\xi=1}^{k_1}
\bigl(P_{\xi}(x)\bigr)^{\widetilde{h}_{0\xi}+
\widetilde{h}_{0,(k_1+\xi)}}
\exp\Bigl( 2\bigl(\, {\stackrel{*}{h}}_{0\xi}-
{\stackrel{*}{h}}_{0,(k_1+\xi)}\bigr)\,\varphi_{\xi}(x)\ +
\hfill
$
\\[2ex]
\centerline{
$
\displaystyle
+\ 2\sum\limits_{q=1}^{\varepsilon_{\xi}}
(1-\delta_{ql}\delta_{\xi\gamma})\Bigl(
\bigl(\,
\widetilde{h}_{q\xi} +
\widetilde{h}_{q,(k_1+\xi)}\bigr)\,
{\stackrel{*}{v}}{}_{q\xi}^{\,\zeta}(x) +
\bigl(\,
{\stackrel{*}{h}}_{q\xi}-{\stackrel{*}{h}}_{q,(k_1+\xi)}
\bigr)\,
\widetilde{v}{}_{q\xi}^{\,\zeta}(x)
\Bigr) \ +
$
}
\\[2ex]
\mbox{}\hfill
$
\displaystyle
+\ 2\bigl(\, {\stackrel{*}{h}}_{l\gamma}\,
\widetilde{v}{}_{l\gamma}^{\,\zeta}(x) +
\widetilde{h}{}_{l\gamma}\,
{\stackrel{*}{v}}{}_{l\gamma}^{\,\zeta}(x)\bigr)
\Bigr)
\prod\limits_{\theta=1}^{k_2}
\bigl(\nu^{0\theta}x\bigr)^{2\,\widetilde{h}_{0\theta}}
\exp\biggl( 2\sum\limits_{q=1}^{\varepsilon_{\theta}}
\widetilde{h}_{q\theta}\,v_{q\theta}^{\,\zeta}(x)\biggr)$
\ \ for all $x\in {\mathscr X},
\hfill
$
\\[2.5ex]
where a domain ${\mathscr X}\subset {\rm D}(F_1)\cap {\rm D}(F_2),\ \delta$ is the Kronecker delta, the functions
\\[2ex]
\mbox{}\hfill
$
P_{\xi}\colon\! x \to
({\stackrel{*}{\nu}}{}^{\,0\xi} x)^{2}  +
(\,\widetilde{\nu}{}^{\,0\xi} x)^{2}\!$
for all $\!x\in {\mathbb R}^n,
\
\varphi_{\xi}\colon x\to
{\rm arctg}\,\dfrac{\widetilde{\nu}{}^{\,0\xi} x}
{{\stackrel{*}{\nu}}{}^{\,0\xi} x}\!$
for all $\!x\!\in\! {\mathscr X},\,
\xi\!=\!1,\ldots,k_1,
\hfill
$
\\[2ex]
and $h_{q\xi}={\stackrel{*}{h}}_{q\xi}+\widetilde{h}_{q\xi}\,i,\ 
h_{q\theta}={\stackrel{*}{h}}_{q\theta}+\widetilde{h}_{q\theta}\,i$ is a nontrivial solution to the linear system 
\\[2.25ex]
\centerline{
$
\displaystyle
\sum\limits_{\xi=1}^{2k_1}
\Bigl(
\lambda_{\xi}^{j}\,h_{0\xi}+
\sum\limits_{q=1}^{\varepsilon_{\xi}}
\mu_{q\xi}^{j\zeta}\,h_{q\xi}\Bigr) -
\mu_{l,(k_1+\gamma)}^{j\zeta}h_{l,(k_1+\gamma)} +
\sum\limits_{\theta=1}^{k_2}
\Bigl(
\lambda_{\theta}^{j}\,h_{0\theta} +
\sum\limits_{q=1}^{\varepsilon_{\theta}}
\mu_{q\theta}^{j\zeta}\,h_{q\theta}\Bigr) \!= 0,
\, j=1,\ldots, m.
$
}
\\[2.25ex]
\indent
Here $\nu^{0\xi}={\stackrel{*}{\nu}}{}^{\,0\xi}\, +\widetilde{\nu}{}^{\,0\xi}\,i, \
\nu^{0,(k_1+\xi)}=\overline{\nu^{0\xi}}$
are complex common eigenvectors of the matrices $A_j$ corresponding to the eigenvalues
\vspace{0.35ex}
$\lambda_{\xi}^{j}={\stackrel{*}{\lambda}}{}_{\xi}^{j} +\widetilde{\lambda}{}_{\xi}^{j}\,i,\ 
\lambda_{k_1+\xi}^{j}=\overline{\lambda_{\xi}^{j}},\ j=1,\ldots, m,\ \xi=1,\ldots, k_1,$ respectively; 
$\nu^{0\theta}$ are real common eigenvectors of the matrices $A_j$ corresponding to the eigenvalues
\vspace{0.35ex}
$\lambda_{\theta}^{j},\ j=1,\ldots,m,\ \theta=1,\ldots, k_2,$ respectively;
the functions $v_{q\xi}^{\,\zeta}={\stackrel{*}{v}}{}_{q\xi}^{\,\zeta} +\widetilde{v}{}_{q\xi}^{\,\,\zeta}\,i,\ 
v_{q\theta}^{\,\zeta}$
is the solution to the system (2.17); the numbers 
\\[2ex]
\mbox{}\hfill
$
{\mu}_{q\xi}^{j\zeta}={\frak p}_{j}\,v_{q\xi}^{\,\zeta}(x),
\ \ \ \
{\stackrel{*}{\mu}}{}_{q\xi}^{\,j\zeta} =\, {\rm Re}\,{\mu}_{q\xi}^{j\zeta},
\ \ \
\widetilde{\mu}{}_{q\xi}^{\,j\zeta} =\,{\rm Im}\,{\mu}_{q\xi}^{j\zeta},
\ \ \
q=1,\ldots,\varepsilon_\xi,
\ \ 
\xi=1,\ldots, 2k_1,
\hfill
$
\\[2ex]
\mbox{}\hfill
$
\mu_{q\theta}^{j\zeta} ={\frak p}_{j}\,v_{q\theta}^{\,\zeta}(x),
\ \ \
q=1,\ldots,\varepsilon_\theta,
\ \ 
\theta=1,\ldots, k_2,
\ \ \ j=1,\ldots,m;
\hfill
$
\\[1.75ex]
the numbers $\varepsilon_{\xi},\ \varepsilon_{\theta}$ such that
$
2\sum\limits_{\xi=1}^{k_1} \varepsilon_{\xi} +
\sum\limits_{\theta=1}^{k_2} \varepsilon_{\theta} =m-2k_{1}-k_{2}+2
$
with
$2k_1+k_2 \leq r, \ \varepsilon_{\xi}\leq s_{\xi}-1,$
$\xi=1,\ldots, k_1,\ \varepsilon_{\theta}\leq s_{\theta}-1, \ \theta=1,\ldots, k_2,$
where $k_1\!$ is a number of  complex common eigen\-vec\-tors {\rm(}this set hasn't complex conjugate vectors{\rm)} 
of  the matrices $A_j$  and  
$k_2$ is a number of real common eigenvectors of the matrices $A_j, \ j=1,\ldots, m.$
\vspace{0.75ex}

In particular, the completely solvable system of total differential equations
\\[1.5ex]
\mbox{}\qquad                                           
$
dx_1 = (3,{}-4,4,1,0,2) x\,dt_1 +  (0,{}-4,2,1,{}-1,1) x\,dt_2,
\hfill
$
\\[1.5ex]
\mbox{}\qquad
$
dx_2 = ({}-1,3,{}-3,0,{}-2,{}-3) x\,dt_1 + (1,3,0,0,1,{}-1) x\,dt_2,
\hfill
$
\\[1.5ex]
\mbox{}\qquad
$
dx_3 = ({}-3,5,{}-5,{}-1,{}-2,{}-4) x\,dt_1 +
(0,6,{}-2,{}-1,2,{}-1) x\,dt_2,
\hfill  (2.32)
$
\\[1.5ex]
\mbox{}\qquad
$
dx_4 = (3,{}-6,4,4,{}-1,5,) x\,dt_1 + (2,{}-6,2,3,{}-4,2) x\,dt_2,
\hfill
$
\\[1.5ex]
\mbox{}\qquad
$
dx_5 = (5,{}-5,8,3,3,6) x\,dt_1 + (1,{}-6,3,2,{}-2,2) x\,dt_2,
\hfill
$
\\[1.5ex]
\mbox{}\qquad
$
dx_6 =
({}-2,5,{}-4,{}-3,1,{}-2) x\,dt_1 + ({}-2,4,{}-3,{}-3,2,{}-2) x\,dt_2
\hfill
$
\\[1.75ex]
has the eigenvalue 
\vspace{0.35ex}
$\lambda_{1}^{1}\!=\!1+2i$ with elementary divisor $\!(\lambda^1\!-1-2i)^3\!$  
corresponding to the eigen\-vec\-tor $\nu^{01}\!=\!(1, 0, 1+i, 1, i, 1)$ and 
\vspace{0.35ex}
to the generalized eigenvectors
$\!\nu^{11}\!=\!(1, 1+i, 0, 0, i, i),$ $\nu^{21}\!=\!(2+2i, 0, 2+2i, 0, 2i, 2i).\!$
Autonomous first integrals of the system (2.32) are the functions
\\[2ex]
\mbox{}\hfill                                          
$
F_1\colon x\to
P(x)\exp\bigl({} -\varphi(x)- \widetilde{v}{}_{11}^{\,1}(x)\bigr)$
\ for all $x\in {\mathscr X},
\hfill (2.33)
$
\\[2ex]
\mbox{}\hfill                                          
$
F_2\colon x\to
P(x)\exp\bigl({}-2\varphi(x) + 2\,{\stackrel{*}{v}}{}_{11}^{\,1}(x)\bigr)$
\ for all $x\in {\mathscr X},
\hfill (2.34)
$
\\[2ex]
\mbox{}\hfill                                          
$
F_3\colon x\to
P^2(x)\exp\big({}-2\varphi(x)- \widetilde{v}{}_{21}^{\,1}(x)\big)$
\ for all $x\in {\mathscr X},
\hfill (2.35)
$
\\[1.25ex]
and
\\[1.25ex]
\mbox{}\hfill                                          
$
F_4\colon x\to {\stackrel{*}{v}}{}_{21}^{\,1}(x)$
\ for all $x\in {\mathscr X},
\quad {\mathscr X}\subset \{ x\colon x_1+x_3+x_4+x_6\ne 0\},
\hfill (2.36)
$
\\[1.75ex]
where the functions $\widetilde{v}{}_{11}^{\,1},\ {\stackrel{*}{v}}{}_{11}^{\,1},\
\widetilde{v}{}_{21}^{\,1},\ {\stackrel{*}{v}}{}_{21}^{\,1},\ P,$ and $\varphi$ are given by (2.22), (2.23), and (2.27).
\vspace{1.25ex}

{\bf Remark 2.1}
Consider the completely solvable system of equations in total differentials 
\\[1.75ex]
\mbox{}\hfill                                     
$
dx_1 = (2x_1+x_3)\,dt_1 +(2x_1+3x_2+3x_3)\,dt_2,
\quad
dx_2 = (x_1+x_2+x_3)\,dt_1 +2x_2\,dt_2,
\hfill                                    
$
\\[0.15ex]
\mbox{}\hfill (2.37)                                    
\\
\mbox{}\hfill                                     
$
dx_3 = {}-x_1\,dt_1 -(3x_2+x_3)\,dt_2,
\quad
dx_4 = x_4\,dt_1 + (x_2+x_3-x_4)\,dt_2.
\hfill 
$
\\[1.75ex]
\indent
The matrices $A_1$ and $A_2$ of (2.37) have the elementary divisors 
\vspace{0.35ex}
$(\lambda^1-1)^2,\  \lambda^1-1,\ \lambda^1-1,$ and
$\lambda^2-2,\ \lambda^2-2, \ (\lambda^2+1)^2,$ respectively.
\vspace{0.35ex}
Using the common eigenvector $\nu^{01}=(1,0,1,0),$ generalized eigenvector $\nu^{11}=(0,1,0,0)$ of $A_1$ and 
\vspace{0.35ex}
the common eigenvector $\nu^{03}=(0,1,1,0),$ generalized eigenvector $\nu^{13}=(0,0,0,1)$ of the matrix $A_2,$
we can build the functions
\\[2ex]
\mbox{}\hfill                                       
$
v_{11}^{1}\colon x\to \dfrac{x_2}{x_1+x_3}\,,
\quad
v_{13}^{2}\colon x\to \dfrac{x_4}{x_2+x_3}$
\ for all $x\in {\mathscr X},
\ \  {\mathscr X}\subset \{ x\colon x_1+x_3\ne 0, x_2+x_3\ne 0\}.
\hfill
$
\\[2ex]
\indent
The system (2.37) has the autonomous first integrals on a domain ${\mathscr X}\colon$
\\[2ex]
\mbox{}\hfill                                       
$
F_1\colon (t,x)\to (x_1+x_3)(x_2+x_3)^2\exp ({}-3v_{11}^{1}(x)),
\quad
F_2\colon (t,x)\to \dfrac{x_2+x_3}{x_1+x_3}\,\exp (3v_{13}^{2}(x)).
\hfill                                       
$
\\[-5ex]

\newpage

{\bf 2.1.3. Nonautonomous first integrals}.
\\[1ex]
\indent
{\bf Theorem 2.5}.
{\it 
Suppose $\nu$ is a real common eigenvector of  the matrices  $A_{j}$ cor\-res\-pon\-ding to the eigenvalues 
$\lambda^{j},\ j=1,\ldots, m,$ respectively.
Then the scalar function
\\[1.25ex]
\mbox{}\hfill
$
\displaystyle
F\colon (t,x)\to
(\nu x) \exp\biggl({}-\sum\limits_{j=1}^{m}\lambda^{j}\,t_{j}\biggr)$
\ \ for all $(t,x)\in {\mathbb R}^{n+m}
\hfill
$
\\[1.25ex]
is a first integral of the linear system of equations in total differentials {\rm (2.1)}.
} 
\vspace{0.75ex}

Consider the system (2.3).  
\vspace{0.25ex}
Using the eigenvalues $\lambda_1^1={}-2,\ \lambda_1^2=1$ 
corresponding to the com\-mon eigenvector $\nu^1=(0,-1,1,1)$ and 
\vspace{0.25ex}
the eigenvalues $\lambda_2^1=0,\ \lambda_2^2={}-1$
corresponding to the eigenvector $\nu^2=(1,0,0,0),$
we can build the first integrals of the system (2.3):
\\[1.5ex]
\mbox{}\hfill
$
F_1\colon (t,x)\to ({}-x_2+x_3+x_4)\exp( 2t_1-t_2),
\ \ \  
F_2\colon (t,x)\to x_1\exp t_2$
\ \ for all $(t,x)\in {\mathbb R}^6.
\hfill
$
\\[1.5ex]
\indent
The functions $F_1,\ F_2,$ (2.4), and (2.5) are the general integral of the system (2.3).
\vspace{1ex}

{\bf  Corollary 2.2}.
{\it 
Let $\nu={\stackrel{*}{\nu}}+\widetilde{\nu}\,i\ 
({\stackrel{*}{\nu}}={\rm Re}\,\nu,\ \widetilde{\nu}={\rm Im}\,\nu)$
be a common complex eigenvector of the matrices $A_j$ cor\-res\-pon\-ding to the eigenvalues
\vspace{0.35ex}
$\lambda^j={\stackrel{*}{\lambda}}{}^{j}+\widetilde{\lambda}{}^{j}\,i\ 
({\stackrel{*}{\lambda}}{}^{j}={\rm Re}\,\lambda^j,\ \widetilde{\lambda}{}^j={\rm Im}\,\lambda^j),$
$j=1,\ldots, m,$ respectively.
Then the sys\-tem {\rm (2.1)} has the first integrals
\\[1.5ex]
\mbox{}\hfill
$
\displaystyle
F_1\colon (t,x)\to
\bigl( ({\stackrel{*}{\nu}}x)^{2} + (\widetilde{\nu}x)^{2}\bigr)
\exp\biggl({}-2\sum\limits_{j=1}^{m}
{\stackrel{*}{\lambda}}{}^{j}\, t_{j}\biggr)$
\ \ for all $(t,x)\in {\mathbb R}^{n+m}
\hfill
$
\\[1ex]
and
\\[1ex]
\mbox{}\hfill
$
\displaystyle
F_2\colon (t,x)\to\ 
{\rm arctg}\,\dfrac{\widetilde{\nu}x}{{\stackrel{*}{\nu}x}}\ - \
\sum\limits_{j=1}^{m} \widetilde{\lambda}{}^{j}\, t_{j}$
\ \ for all $(t,x)\in {\mathscr D},
\quad {\mathscr D}\subset {\mathbb R}^{n+m}.
\hfill
$
}
\\[2ex]
\indent
For example, the completely solvable system (2.8) has 
the eigenvector $\nu^{1}=(1, i, 0)$ corres\-pon\-ding to the eigenvalues $\lambda_1^1=1, \lambda_1^2={}-i$
and the first integrals (by Corollary 2.2)
\\[1.5ex]
\mbox{}\hfill
$
F_1\colon (t,x)\to  (x_1^2+x_2^2)\exp({}-2t_1)$
\ for all $(t,x)\in {\mathbb R}^5,
\hfill
$
\\[2.5ex]
\mbox{}\hfill
$
F_2\colon (t,x)\to \ 
{\rm arctg}\,\dfrac{x_2}{x_1} +t_2$
\ for all $(t,x)\in {\mathbb R}^2\times {\mathscr X},
\quad {\mathscr X}\subset \{ x\colon x_1\ne 0,\ x_3\ne 0\}.
\hfill
$
\\[1.5ex]
\indent
The functionally independent first integrals (2.9), $F_1,$ and $F_2$ are 
the general integral on a domain ${\mathbb R}^2\times {\mathscr X}$ for the system 
of equations in total differentials (2.8). 
\vspace{0.5ex}

The completely solvable system (2.12) has the eigenvector $\nu^{1}\!=\! (0, -i, 0, 1)$
corresponding to the eigenvalues $\lambda_1^1={}-1$ and $\lambda_1^2={}-i.$
The first integrals (2.13) and (by Corollary 2.2)
\\[1.75ex]
\mbox{}\hfill
$
F_1\colon (t,x)\to  (x_2^2+x_4^2)\exp(2t_1)$
\ for all $(t,x)\in {\mathbb R}^6,
\hfill
$
\\[2.5ex]
\mbox{}\hfill
$
F_2\colon (t,x)\to \ 
{\rm arctg}\,\dfrac{x_2}{x_4} - t_2$
\ for all $(t,x)\in {\mathbb R}^2\times {\mathscr X},
\ \ \
{\mathscr X}\subset \{x\colon x_1x_4-x_2x_3\ne 0,\, x_4\ne 0\},
\hfill
$
\\[1.5ex]
are the general integral on a domain ${\mathbb R}^2\times {\mathscr X}$ of the linear system (2.12). 
\vspace{1.25ex}

{\bf Theorem 2.6}.                                           
{\it
Let $\nu^{\,0}$ be a real common eigenvector of the matrices $A_j$ corresponding to the eigenvalues
$\lambda^j, \ j=1,\ldots,m,$ respectively. Let $\nu^{\,\theta},\ \theta =1,\ldots, s-1$ be real generalized eigenvectors
of the matrix $A_\zeta$ 
\vspace{0.35ex}
corresponding to the eigenvalue $\lambda^{\zeta}$ with elementary divisor of multiplicity  $s\geq 2.$
Then the completely solvable sys\-tem {\rm (2.1)} has the first integrals
\\[1.5ex]
\mbox{}\hfill
$
\displaystyle
F_q\colon (t, x)\to
v_{q}^{\zeta}(x)\ -\
\sum\limits_{j=1}^{m}\mu_{q}^{j\zeta}\, t_{j}$
\ for all $(t,x)\in {\mathbb R}^{m}\times {\mathscr X},
\ \ \ q=1,\ldots, s-1,
\hfill
$
\\[1.5ex]
where the set of functions 
\vspace{0.25ex}
$v_{q}^{\zeta}\colon {\mathscr X}\to {\mathbb R}$ is the solution to the system {\rm (2.17)}, 
the numbers $\mu_{q}^{j\zeta} = {\frak  p}_{j}\,v_{q}^{\zeta}(x),\
q=1,\ldots, s-1,\ j=1,\ldots, m,$ and a domain ${\mathscr X}\subset \{x\colon \nu^0x\ne 0\}.$
}

\newpage

The completely solvable system (2.18) has 
$\mu_{1}^{11}=1, \ \mu_{3}^{11}=0,\  \mu_{1}^{21}={}-1, \ \mu_{3}^{21}=6,$ and 
the first integrals (by Theorem 2.6) 
\\[1.5ex]
\mbox{}\hfill
$
F_1\colon (t,x)\to v_{11}^{\,1}(x) -t_1+t_2$
\ \ for all $(t,x)\in {\mathbb R}^2\times {\mathscr X},
\hfill
$
\\[2ex]
\mbox{}\hfill
$
F_2\colon (t,x)\to v_{31}^{\,1}(x) -6t_2$
\ for all $(t,x)\in {\mathbb R}^2\times {\mathscr X},
\quad {\mathscr X}\subset \{x\colon x_1-x_2+x_3\ne 0\},
\hfill
$  
\\[1.5ex]
where the functions 
\vspace{0.35ex}
$v_{11}^{1},\ v_{31}^{1}$ are given by (2.19). The functionally independent
first integrals (2.20), $F_1,$ and $F_2$
are the general integral on a domain ${\mathbb R}^2\times {\mathscr X}$ of the system (2.18).
\vspace{1.25ex}

{\bf Corollary 2.3}.                                           
{\it
Let $\nu^{\,0}$ be a common eigenvector of the matrices $A_j$ corresponding to the eigenvalues
$\lambda^j, \ j=1,\ldots,m.$ Suppose $\nu^{\,\theta},\ \theta =1,\ldots, s-1$ are generalized eigenvectors
of $A_\zeta$ 
\vspace{0.35ex}
corresponding to the complex eigenvalue $\lambda^{\zeta}\ ({\rm Im}\,\lambda^{\zeta}\ne 0)$ with elementary divisor of multiplicity  $s\geq 2.$
Then the completely solvable sys\-tem {\rm (2.1)} has the first integrals
\\[1.5ex]
\mbox{}\hfill
$
\displaystyle
F_{1 q}\colon (t,x)\to
{\stackrel{*}{v}}{}_{q}^{\zeta}(x) -
\sum\limits_{j=1}^{m}\,
{\stackrel{*}{\mu}}\!{}_{q}^{\,j\zeta}\, t_{j}$
\ \ for all $(t,x)\in {\mathbb R}^{m}\times {\mathscr X},
\quad q=1,\ldots, s-1,
\hfill
$
\\[1ex]
and
\\[1ex]
\mbox{}\hfill
$
\displaystyle
F_{2 q}\colon (t, x)\to
\widetilde{v}{}_{q}^{\,\zeta}(x) -
\sum\limits_{j=1}^{m}\,
\widetilde{\mu}{}_{q}^{\,j\zeta}\, t_{j}$
\ \ for all $(t,x)\in {\mathbb R}^{m}\times {\mathscr X},
\quad  q=1,\ldots, s-1,
\ \ \ \ {\mathscr X}\subset {\mathbb R}^{n},
\hfill
$
\\[1.5ex]
where the set of functions 
\vspace{0.25ex}
$v_{q}^{\zeta}\colon x\to {\stackrel{*}{v}}{}_{q}^{\,\zeta}(x) + \widetilde{v}{}_{q}^{\,\zeta}(x)\,i$ 
for all  $x\in {\mathscr X}$ is the solution to the system {\rm (2.17)}, the numbers 
${\stackrel{*}{\mu}}{}_{q}^{\,j\zeta}= {\frak p}_j\,{\stackrel{*}{v}}{}_{q}^{\,\zeta}(x), \
\widetilde{\mu}{}_{q}^{\,j\zeta}={\frak p}_j\,\widetilde{v}{}_{q}^{\,\zeta}(x), \
j=1,\ldots, m,\ q=1,\ldots, s-1.$
}
\vspace{1.5ex}

The completely solvable system (2.21) has the first integrals (by Corollary 2.3) 
\\[1ex]
\mbox{}\hfill
$
\displaystyle
F_{11}\colon (t,x)\to
{\stackrel{*}{v}}{}_{11}^{\,1}(x) -t_1-t_2,
\ \ \,
F_{21}\colon (t,x)\to
\widetilde{v}{}_{11}^{\,1}(x) +t_2-t_3$
\ for all $(t,x)\in {\mathbb R}^3\times {\mathscr X},
\hfill
$
\\[1ex]
and
\\[1ex]
\mbox{}\hfill
$
F_{12}\colon (t,x)\to
{\stackrel{*}{v}}{}_{21}^{\,1}(x) -2t_3$
\ for all $(t,x)\in {\mathbb R}^3\times {\mathscr X},
\quad
{\mathscr X}\subset \{x\colon x_1+x_3+x_4+x_6\ne 0\},
\hfill
$
\\[1.5ex]
where ${\stackrel{*}{v}}{}_{11}^{\,1},\ \widetilde{v}{}_{11}^{\,1},$ and ${\stackrel{*}{v}}{}_{21}^{\,1}$ 
\vspace{0.35ex}
are given by (2.22). The general integral for the linear system (2.21) is the
functionally independent first integrals (2.24),\, (2.25),\, (2.26), $F_{11},\ F_{21},$ and $F_{12}.$  
\vspace{1ex}

For example, the system (2.28) has the  
\vspace{0.35ex}
eigenvector $\nu^{01}=(1,1+i,0,0,i,i)$
corresponding to the eigenvalues $\lambda_1^1=1+i,\ \lambda_1^2={}-1,\ \lambda_1^3=2,\ \lambda_1^4=-1+2i$
and the first integrals 
\\[2ex]
\mbox{}\hfill
$
F_1\colon (t,x)\to
\bigl( (x_1+x_2)^2+(x_2+x_5+x_6)^2\bigr)
\exp({}-2t_1+2t_2-4t_3+2t_4)$
\ \ (by Corollary 2.2),
\hfill\mbox{}
\\[2.5ex]
\mbox{}\hfill
$
F_2\colon (t,x)\to\ 
{\rm arctg}\,\dfrac{x_2+x_5+x_6}{x_1+x_2} -t_1-2t_4$
\ for all $(t,x)\in\R^{4}\times {\mathscr X},
\hfill
$
\\[2ex]
where ${\mathscr X}\subset \{x\colon x_1+x_2\ne 0,\ x_1+x_3+x_4+x_6\ne 0\}\subset {\mathbb R}^6.$
Using Corollary 2.3, we get 
\\[1.75ex]
\mbox{}\hfill
$
F_3\colon  (t,x)\to
{\stackrel{*}{v}}{}_{11}^{\,1}(x)-t_1+t_3-t_4$
\ for all $(t,x)\in {\mathbb R}^{4}\times {\mathscr X}
\hfill
$
\\[0.5ex]
and
\\[0.5ex]
\mbox{}\hfill
$
F_4\colon (t,x)\to
\widetilde{v}{}_{11}^{\,1}(x)-t_2-t_4$
\ for all $(t,x)\in {\mathbb R}^{4}\times {\mathscr X},
\hfill
$
\\[1.75ex]
where ${\stackrel{*}{v}}{}_{11}^{\,1}$ and $\widetilde{v}{}_{11}^{\,1}$ are given by (2.29). 
\vspace{0.25ex}
The first integrals (2.30),\, (2.31), $F_1,\ldots, F_4$ are   
the general integral on a domain ${\mathbb R}^4\times {\mathscr X}$
for the completely solvable linear system (2.28).
\vspace{1ex}

Consider the completely solvable linear system of equations in total differentials (2.32).  
According to Corollary 2.3, we can construct the first integrals
\\[1.75ex]
\mbox{}\hfill
$
F_{11}\colon (t,x)\to
{\stackrel{*}{v}}{}_{11}^{\,1}(x)-t_1-t_2,
\ \ \ 
F_{21}\colon (t,x)\to
\widetilde{v}{}_{11}^{\,1}(x)+t_2$
\ for all $(t,x)\in {\mathbb R}^2\times  {\mathscr X},
\hfill
$
\\[1.75ex]
where the functions ${\stackrel{*}{v}}{}_{11}^{\,1}$ and $\widetilde{v}{}_{11}^{\,1}$ are given by (2.22), 
\vspace{0.5ex}
a domain ${\mathscr X}\subset \{x\colon x_1+x_3+x_4+x_6\ne 0\}.$

\newpage

The functionally independent first integrals (2.33),\, (2.34),\, (2.35),\, (2.36), $F_{11},$ and $F_{21}$ are 
the general integral on a domain ${\mathbb R}^2\times {\mathscr X}$
\vspace{1ex}
for the completely solvable system (2.32).

{\bf Remark 2.2}. The system of equations in total differentials 
\\[1.5ex]
\mbox{}\hfill                                     
$
dx_1 = {}-(x_1+x_2+3x_3)\,dt_1 -3(x_1+2x_3)\,dt_2,
\hfill (2.38)                                    
$
\\[1.5ex]
\mbox{}\hfill                                     
$
dx_2 = x_2\,dt_1 -2x_2\,dt_2,
\quad
dx_3 = (x_2+2x_3)\,dt_1 + 3x_3\,dt_2
\hfill 
$
\\[1.5ex]
is not completely solvable (this system has the defect [5, p. 54] of order 1):
\\[1.5ex]
\mbox{}\hfill
$
[{\frak p}_1(x),{\frak p}_2(x)]=[{}-(x_1+x_2+3x_3)\,\partial_{x_1}+x_2\,\partial_{x_2}+(x_2+2x_3)\,\partial_{x_3},
\hfill                                     
$
\\[1.5ex]
\mbox{}\hfill                                     
$
{}-3(x_1+2x_3)\,\partial_{x_1}-2x_2\,\partial_{x_2}+3x_3\,\partial_{x_3}]=
{}-5x_2\partial_{x_1}+5x_2\partial_{x_3}\equiv {\frak p}_3(x)$
\ for all $x\in {\mathbb R}^3,
\hfill
$
\\[1.5ex]
\mbox{}\hfill
$
[{\frak p}_1(x),{\frak p}_3(x)]={}-{\frak p}_3(x),
\quad
[{\frak p}_2(x),{\frak p}_3(x)]={}-5\,{\frak p}_3(x)$
\ for all $x\in {\mathbb R}^3.
\hfill                                     
$
\\[1.5ex]
\indent
The system (2.38) has two common eigenvectors $\nu^{1}\!=\!(1,0,1)$ and $\nu^{2}\!=\!(0,1,0)$
cor\-re\-s\-po\-n\-ding to the eigenvalues $\lambda_1^1\!=\!{}-1,\, 
\lambda_1^2=\!{}-3,$ and $\lambda_2^1=1,\ \lambda_2^2=\!{}-2.$
The first integrals 
\\[1.5ex]
\mbox{}\hfill
$
F_{1}\colon (t,x)\to (x_1+x_3)\exp(t_1+3t_2),
\ \ 
F_{2}\colon (t,x)\to x_2\exp(2t_2-t_1)$
\ for all $(t,x)\in {\mathbb R}^5
\hfill
$
\\[1.5ex]
are a general integral of the linear system of equations in total differentials (2.38).
\\[2.5ex]
\indent
{\bf 2.2. Linear real nonhomogeneous systems in total differentials}
\\[1ex]
\indent
Consider now a nonhomogeneous system of equations in total differentials
\\[1.5ex]
\mbox{}\hfill                                          
$
\displaystyle
dx =\sum\limits_{j=1}^m(B_j\,x+f_j(t))\,dt_j
$
\hfill (2.39)
\\[1.5ex]
corresponding to the linear homogeneous system (2.1), where
the matrices $B_1,\ldots, B_m$ are transpose of $A_1,\ldots, A_m,$ respectively, 
and the vector functions
\\[1.5ex]
\mbox{}\hfill
$
f_j\colon t\to {\rm colon}(f_{j1}(t),\ldots, f_{jn}(t))$
\ for all $t\in {\mathcal T},
\quad
j=1,\ldots, m
\hfill
$
\\[1.5ex]
are continuously differentiable on a domain ${\mathcal T}\subset\R^{m}.$

The Frobenius conditions for the total solvability [4, p. 44] of system (2.39) are (2.2) and 
\\[1.5ex]
\mbox{}\hfill                                                          
$
\partial_{t_\zeta}f_j(t)-B_{\zeta}f_j(t)=
\partial_{t_j}f_\zeta(t)-B_{j}f_\zeta(t)$
\ for all $t\in {\mathcal T},
\quad 
j,\zeta=1,\ldots, m.
\hfill (2.40)
$
\\[2ex]
\indent
{\bf Theorem 2.7}.
{\it
Let $\nu$ be a real common eigenvector of  the matrices  $A_{j}$ cor\-res\-pon\-ding to the eigenvalues 
$\lambda^{j},\ j=1,\ldots, m,$ respectively.
Then the completely solvable system of equations in total differentials {\rm (2.39)} has the first integral
\\[1.75ex]
\mbox{}\hfill                                            
$
\displaystyle
F\colon (t,x)\to 
\nu x\!\cdot\!\varphi(t)-
\int\sum\limits_{j=1}^m \nu f_j(t)\!\cdot\! \varphi(t)\,dt_j$
\ \ for all $(t,x)\in \widetilde{\mathcal T}\times\R^n,
$
\hfill {\rm (2.41)}
\\[1.75ex]
where $\widetilde{\mathcal T}$ is a simply connected domain  from ${\mathcal T},$ the function 
\\[1.25ex]
\mbox{}\hfill                                            
$
\displaystyle
\varphi\colon t\to 
\exp\biggl({}-\sum\limits_{j=1}^m \lambda^j\,t_j\biggr)$
\ \ for all $t\in {\mathcal T}.
\hfill
$
\\[1ex]
}
\indent
{\sl Proof}. Let us introduce the 1-form
\\[1.75ex]
\mbox{}\hfill                                           
$
\displaystyle
\omega\colon t\to 
\sum\limits_{j=1}^m \nu f_j(t)\!\cdot\! \varphi(t)\,dt_j$
\ \ for all $t\in {\mathcal T},
\quad {\mathcal T}\subset {\mathbb R}^m.
\hfill
$
\\[1.5ex]
\indent
The external differential has the form
\\[1.5ex]
\mbox{}\hfill                                           
$
\displaystyle
d\,\omega(t)=\varphi(t)\sum\limits_{j=1}^m \sum\limits_{\zeta=1}^m
\bigl(\partial_{t_\zeta} \nu f_j(t) -
\lambda^{\zeta}\cdot\nu f_j(t)\bigr)\,dt_{\zeta}\wedge dt_j$
\ \ \ for all $t\in {\mathcal T}.
\hfill
$
\\[1.75ex]
\indent
Using (2.40) and 
\\[1.25ex]
\mbox{}\hfill
$
dt_j\wedge dt_j=0,\ \ \ dt_\zeta\wedge dt_j={}-dt_j\wedge dt_\zeta, 
\quad 
\zeta, j=1,\ldots, m,
\hfill
$
\\[1.5ex]
we get the differential 1-form $\omega$ is 
\vspace{0.35ex}
closed on ${\mathcal T}.$ By the Poincar\'{e} theorem (see [5, p. 14]), it follows that $\omega$ is 
exact in a simply connected domain $\widetilde{\mathcal T}\subset {\mathcal T}.$
\vspace{0.35ex}

Taking into account Lemma 2.1 and (2.15), we obtain
\\[1.75ex]
\mbox{}\hfill
$
\displaystyle
{\frak P}_j F(t,x)=
{}-\lambda^j\nu x\cdot\varphi(t)
-\nu f_j(t)\varphi(t)+{\frak p}_j\nu x\cdot\varphi(t)
+\nu f_j(t)\varphi(t)=0
\hfill
$
\\[2ex]
\mbox{}\hfill
for all 
$(t,x)\in \widetilde{\mathcal T}\times\R^n,
\quad
j=1,\ldots, m,
\hfill 
$
\\[1.5ex]
where the linear differential operators
\\[1.75ex]
\mbox{}\hfill
$
\displaystyle
{\frak P}_j(t,x)=
\partial_{t_j}+ {\frak p}_j(x) + f_j(t)\,\partial_{x}$ 
\ \ for all $(t,x)\in {\mathcal T}\times {\mathbb R}^{n},
\ \ \ 
j=1,\ldots, m.
\hfill
$
\\[1.75ex]
\indent\!
Therefore the function \!(2.41)\! is a 
\vspace{1.25ex}
first integral of the completely solvable system \!(2.39).\k

{\bf Corollary 2.4}.
{\it 
Let $\nu={\stackrel{*}{\nu}}+\widetilde{\nu}\,i\ 
({\stackrel{*}{\nu}}={\rm Re}\,\nu,\ \widetilde{\nu}={\rm Im}\,\nu)$
be a common complex eigenvector of the matrices $A_j$ cor\-res\-pon\-ding to the eigenvalues
\vspace{0.35ex}
$\lambda^j={\stackrel{*}{\lambda}}{}^{j}+\widetilde{\lambda}{}^{j}\,i\ 
({\stackrel{*}{\lambda}}{}^{j}={\rm Re}\,\lambda^j,\ \widetilde{\lambda}{}^j={\rm Im}\,\lambda^j),$
$j\!=\!1,\ldots,\! m,\!$ respectively.
Then the completely solvable sys\-tem {\rm (2.39)} has the first integrals
\\[2ex]
\mbox{}\hfill                                          
$
\displaystyle
F_\tau\colon (t,x)\to \alpha_\tau(t, x)-
\int \sum\limits_{j=1}^m\alpha_\tau(t, f_j(t))\,dt_j$
\ \ for all $(t,x)\in \widetilde{\mathcal T}\times\R^n,
\ \ \ \tau=1,2,
\hfill
$
\\[2ex]
where the scalar functions of the vector arguments
\\[2ex]
\mbox{}\hfill                                     
$
\displaystyle
\alpha_1(t,x)=
\biggl({\stackrel{*}{\nu}}x\cdot
\cos\sum\limits_{j=1}^{m}\widetilde{\lambda}{}^j\;\!t_j  +
\widetilde{\nu}x\cdot\sin\sum\limits_{j=1}^{m}\widetilde{\lambda}{}^j\;\!t_j\biggr)\!
\cdot\exp\biggl(-\sum\limits_{j=1}^{m}{\stackrel{*}{\lambda}}{}^j\;\!t_j\biggr)$
\ for all $(t,x)\in {\mathcal T}\times\R^n,
\hfill
$
\\[2.5ex]
\mbox{}\hfill                                     
$
\displaystyle
\alpha_2(t,x)=
\biggl(\widetilde{\nu}x\cdot
\cos\sum\limits_{j=1}^{m}\widetilde{\lambda}{}^j\;\!t_j  -
{\stackrel{*}{\nu}}x\cdot\sin\sum\limits_{j=1}^{m}\widetilde{\lambda}{}^j\;\!t_j\biggr)\!
\cdot\exp\biggl({}-\sum\limits_{j=1}^{m}{\stackrel{*}{\lambda}}{}^j\;\!t_j\biggr)$
\ for all $(t,x)\in {\mathcal T}\times\R^n.
\hfill                                     
$
\\[2.25ex]
}
\indent
For example, the completely solvable system of total differential equations 
\\[1.75ex]
\mbox{}\hfill                                      
$
dx_1 = x_1\,dt_1 + (x_2+1)\,dt_2,
\quad \
dx_2 = (x_2+1)\,dt_1 +({}- x_1+e^{t_1+at_2})\,dt_2,
\hfill                                      
$
\\[1.75ex]
\mbox{}\hfill                                      
$
dx_3 = (x_3+t_2-t_1)\,dt_1 +({}- x_3+t_1-t_2)\,dt_2
\hfill
$
\\[1.75ex]
cor\-res\-pon\-ding to the linear homogeneous system (2.8). This system has the first integrals
\\[2ex]
\mbox{}\hfill                                      
$
F_1\colon (t,x)\to 
(x_1\cos t_2-x_2\sin t_2-\sin t_2)e^{{}-t_1}+\dfrac{a\sin t_2-\cos t_2}{a^2+1}\ e^{at_2}$
\ \ (by Corollary 2.4),
\hfill\mbox{}                                      
\\[2.5ex]
\mbox{}\hfill                                      
$
F_2\colon (t,x)\to 
(x_1\sin t_2+x_2\cos t_2+\cos t_2)e^{{}-t_1}-\dfrac{a\cos t_2+\sin t_2}{a^2+1}\ e^{at_2}$
\ \ (by Corollary 2.4),
\hfill\mbox{}                                      
\\[2.5ex]
\mbox{}\hfill                                      
$
F_3\colon (t,x)\to 
(x_3+t_2-t_1-1)e^{t_2-t_1}$
\ \ for all $(t,x)\in {\mathbb R}^5$
\ \ (by Theorem 2.7),
\hfill\mbox{}
\\[2ex]
where $a$ is some real number. 
\vspace{1ex}

{\bf Theorem 2.8}.\!
{\it 
Let the assumptions of Lemma {\rm 2.2} hold. 
Then the completely solvable sys\-tem  of equations in total differentials {\rm (2.39)} has the first integrals
\\[1.75ex]
\mbox{}\hfill                                  
$
F_{\theta}\colon (t,x)\to \nu^{\,\theta l}x\cdot \varphi(t) -
{\displaystyle \sum\limits_{\tau=1}^{\theta} }
K_{\tau-1}^{\theta}(t)\cdot F_{\tau-1}(t,x)-C_{\theta}(t)
\hfill                                  
$
\\
\mbox{}\hfill  {\rm (2.42)}                                
\\
\mbox{}\hfill                                  
for all $(t,x)\in \widetilde{\mathcal T}\times {\mathscr X},
\ \ \ \
\theta=0,\ldots, s_{l}-1,
\quad
 {\mathscr X}\subset \{x\colon \nu^{0}x\ne 0\},
\hfill 
$
\\[2ex]
where the functions 
\\[1.5ex]
\mbox{}\hfill
$
\displaystyle
\varphi\colon t\to \exp\biggl({}-\sum\limits_{j=1}^m \mu_{0l}^{j\zeta}\,t_j\biggr)$
\ \ for all $t\in \widetilde{\mathcal T},
\quad
\widetilde{\mathcal T}\subset {\mathcal T},
\hfill
$
\\[2ex]
\mbox{}\hfill
$
K_{\tau-1}^{\theta}\colon t\to 
{\displaystyle \int} 
{\displaystyle \sum\limits_{j=1}^m}
\biggl(\binom{\theta}{\tau-1}\,\mu_{\theta-\tau+1,l}^{\,j\zeta}+
{\displaystyle \sum\limits_{\delta=1}^{\theta-\tau}}
\binom{\theta}{\delta}\,\mu_{\delta l}^{j\zeta}\cdot K_{\tau-1}^{\theta-\delta}(t)\biggr)\,dt_j,
\ \ \ 
\tau=1,\ldots, \theta,
\ \ \theta=1,\ldots, s_{l}-1,
\hfill
$
\\[2ex]
\mbox{}\hfill
$
C_{\theta}\colon t\to 
{\displaystyle \int\sum\limits_{j=1}^m} 
\biggl(\nu^{\,\theta l}f_j(t)\cdot \varphi(t)+
{\displaystyle \sum\limits_{\tau=1}^{\theta}}
\binom{\theta}{\tau}\,\mu_{\tau l}^{j\zeta}\cdot C_{\theta-\tau}(t)\biggr)\,dt_j$
\ for all $t\in \widetilde{\mathcal T},
\quad 
\theta=0,\ldots, s_{l}-1,
\hfill
$
\\[2ex]
the numbers $\mu_{0l}^{j\zeta}=\lambda^j_{l}, \ \mu_{\theta l}^{j\zeta}={\frak p}_{j}\,v_{\theta l}^{\zeta}(x),
\ \theta=1,\ldots, s_{l}-1, \ j=1,\ldots, m.$ 
}
\vspace{1.25ex}

{\sl Proof}. The proof is by induction on $\theta.$

The case $\theta=0$ was considered in Theorem 2.7.

Suppose $\theta=1.$ Using Lemma 2.3 with condition (2.16), we obtain  
\\[2ex]
\mbox{}\hfill                                                 
$
\displaystyle
{\frak P}_j F_1(t,x)= {\frak P}_j \biggl(\!\nu^{1l}x\cdot \varphi(t)-
\int\sum\limits_{\xi=1}^m \mu_{1l}^{\xi\zeta}\,dt_{\xi}\cdot F_{0}(t,x) -
\int\sum\limits_{\xi=1}^m
\bigl(\nu^{1l}f_{\xi}(t)\cdot \varphi(t)+\mu_{1l}^{\xi\zeta}\,C_{0}(t)\bigr)\,dt_{\xi}\!\biggr)=
\hfill                                                 
$
\\[2.25ex]
\mbox{}\hfill                                                 
$
\displaystyle
=\mu_{1l}^{j\zeta}\bigl(\nu^{0l}x\cdot \varphi(t)-C_{0}(t)-F_{0}(t,x)\bigr)=0$
\ \ for all $(t,x)\in \widetilde{\mathcal T}\times {\mathscr X},
\ \ \  j=1,\ldots, m.
\hfill 
$
\\[2ex]
\indent
Therefore $F_1\colon \widetilde{\mathcal T}\times {\mathscr X}\to {\mathbb R}$ 
is a first integral of the completely solvable system (2.39).
\vspace{0.35ex}

Suppose that the assertion of the theorem is valid for 
\vspace{0.35ex}
$\theta=\varepsilon-1,$ i.e., the scalar functions
$F_{\theta}\colon \widetilde{\mathcal T}\times {\mathscr X}\to {\mathbb R},\ \theta=1,\ldots,\varepsilon-1$
are first integrals of the system (2.39). Then 
\\[2ex]
\mbox{}\hfill                                 
$
{\frak P}_j F_{\varepsilon}(t,x)=
{\frak P}_j \biggl( \nu^{\,\varepsilon l}x\cdot \varphi(t) -
{\displaystyle \sum\limits_{\tau=1}^{\varepsilon} }
K_{\tau-1}^{\varepsilon}(t)\cdot F_{\tau-1}(t,x)-C_{\varepsilon}(t)\biggr)=
\hfill                                  
$
\\[2ex]
\mbox{}\hfill                                                 
$
={}-\mu_{0l}^{j\zeta}\,\nu^{\,\varepsilon l}x\cdot \varphi(t) +
{\displaystyle \sum\limits_{\tau=0}^{\varepsilon} }
\binom{\varepsilon}{\tau}\,\mu_{\tau l}^{j\zeta}\, \nu^{\,\varepsilon-\tau,l}x\cdot \varphi(t)
+ \nu^{\,\varepsilon l}f_{j}(t)\cdot \varphi(t)\ -
\hfill 
$
\\[2ex]
\mbox{}\hfill                                                 
$
-\ {\displaystyle \sum\limits_{\tau=1}^{\varepsilon}}\biggl(\!\!
\binom{\varepsilon}{\tau-1}\,\mu_{\varepsilon-\tau+1,l}^{\, j\zeta}+
{\displaystyle \sum\limits_{\delta=1}^{\varepsilon-\tau}}\binom{\varepsilon}{\delta}\,
\mu_{\delta l}^{j\zeta}\, K_{\tau-1}^{\varepsilon-\delta}(t)\!\!\biggr) F_{\tau-1}(t,x) -
\biggl(\!\nu^{\,\varepsilon l}f_j(t)\cdot \varphi(t) +
{\displaystyle \sum\limits_{\tau=1}^{\varepsilon} }
\binom{\varepsilon}{\tau}\,\mu_{\tau l}^{j\zeta}\, C_{\varepsilon-\tau}(t)\!\!\biggr)\!=
\hfill 
$
\\[2ex]
\mbox{}\hfill                                                 
$
={\displaystyle \sum\limits_{\tau=1}^{\varepsilon} }
\binom{\varepsilon}{\tau}\,\mu_{\tau l}^{j\zeta}
\biggl(\!
\Bigl( \nu^{\,\varepsilon-\tau, l}x\cdot \varphi(t) -
{\displaystyle \sum\limits_{\delta=1}^{\varepsilon-\tau} }
K_{\delta-1}^{\varepsilon-\tau}(t)\cdot F_{\delta-1}(t,x) -
C_{\varepsilon-\tau}(t)\Bigr)
-F_{\varepsilon-\tau}(t,x)\biggr)=0
\hfill 
$
\\[2ex]
\mbox{}\hfill                                                 
for all $(t,x)\in \widetilde{\mathcal T}\times {\mathscr X},
\quad
j=1,\ldots, m.
\hfill 
$
\\[2ex]
\indent
Therefore $F_{\varepsilon}\colon \widetilde{\mathcal T}\times {\mathscr X}\to {\mathbb R}$ 
is a first integral of the linear system (2.39).
\vspace{0.35ex}

Thus the functions (2.42) are first integrals of the completely solvable system (2.39). \ \k
\vspace{1ex}

{\bf Remark 2.3}. In complex case from Theorem 2.8, we get the following real-valued first in\-teg\-rals 
of the completely solvable system of equations in total differentials (2.39):
\\[2ex]
\mbox{}\hfill
$
F_{\theta}^1\colon (t,x)\to {\rm Re}\,F_{\theta}(t,x),
\ \, 
F_{\theta}^2\colon (t,x)\to {\rm Im}\,F_{\theta}(t,x)$
for all $(t,x)\in \widetilde{\mathcal T}\times {\mathscr X},
\ \theta=0,\ldots, s_{l}-1.
\hfill
$
\\[2ex]
\indent
Moreover, we have
\\[1.5ex]
\mbox{}\hfill                                         
$
F_{\theta}^1\colon (t,x)\to \alpha_{{}_{\scriptstyle \theta}}(t,x)-
{\displaystyle \sum\limits_{\tau=1}^{\theta}}
\Bigl({\rm Re}\,K_{\tau-1}^{\theta}(t)\cdot F_{\tau-1}^1(t,x) -
{\rm Im}\,K_{\tau-1}^{\theta}(t)\cdot F_{\tau-1}^2(t,x)\Bigr)
-{\rm Re}\,C_{\theta}(t),
\hfill                                  
$
\\[2ex]
\mbox{}\hfill                                         
$
F_{\theta}^2\colon (t,x)\to \beta_{{}_{\scriptstyle \theta}}(t,x)-
{\displaystyle \sum\limits_{\tau=1}^{\theta}}
\Bigl({\rm Re}\,K_{\tau-1}^{\theta}(t)\cdot F_{\tau-1}^2(t,x) +
{\rm Im}\,K_{\tau-1}^{\theta}(t)\cdot F_{\tau-1}^1(t,x)\Bigr)
-{\rm Im}\,C_{\theta}(t)
\hfill      
$
\\[1.75ex]
\mbox{}\hfill                                         
for all $(t,x)\in \widetilde{\mathcal T}\times {\mathscr X},
\quad 
\theta=0,\ldots, s_{l}-1,
\hfill                                  
$
\\[1.5ex]
where the functions 
\\[2.5ex]
\mbox{}\hfill
$
{\rm Re}\, K_{\tau-1}^{\theta}\colon t\to 
{\displaystyle \int} 
{\displaystyle \sum\limits_{j=1}^m}
\biggl(\binom{\theta}{\tau-1}\,{\stackrel{*}{\mu}}{}_{\theta-\tau+1,l}^{\,j\zeta} +
{\displaystyle \sum\limits_{\delta=1}^{\theta-\tau}}\binom{\theta}{\delta}\,
\bigl(\,{\stackrel{*}{\mu}}{}_{\delta l}^{\,j\zeta}\, {\rm Re}\,K_{\tau-1}^{\theta-\delta}(t) -
\widetilde{\mu}{}_{\delta l}^{\,j\zeta}\, {\rm Im}\,K_{\tau-1}^{\theta-\delta}(t)\bigr)\biggr)\,dt_j\,,
\hfill
$
\\[3ex]
\mbox{}\hfill
$
{\rm Im}\, K_{\tau-1}^{\theta}\colon t\to 
{\displaystyle \int} 
{\displaystyle \sum\limits_{j=1}^m}
\biggl(\binom{\theta}{\tau-1}\,\widetilde{\mu}{}_{\theta-\tau+1, l}^{\,j\zeta} +
{\displaystyle \sum\limits_{\delta=1}^{\theta-\tau}}\binom{\theta}{\delta}\,
\bigl({\stackrel{*}{\mu}}{}_{\delta l}^{\,j\zeta}\, {\rm Im}\,K_{\tau-1}^{\theta-\delta}(t) +
\widetilde{\mu}{}_{\delta l}^{\,j\zeta}\, {\rm Re}\,K_{\tau-1}^{\theta-\delta}(t)\bigl)\biggr)\,dt_j
\hfill
$
\\[2.75ex]
\mbox{}\hfill
for all 
$t\in \widetilde{\mathcal T},
\quad
\tau=1,\ldots,\theta,
\ \ \ \theta=1,\ldots, s_{l}-1,
\hfill
$
\\[3ex]
\mbox{}\hfill
$
{\rm Re}\,C_{\theta}\colon t\to 
{\displaystyle \int\sum\limits_{j=1}^m} 
\biggl(\alpha_{{}_{\scriptstyle \theta}}(t, f_j(t)) +
{\displaystyle \sum\limits_{\tau=1}^{\theta}}\binom{\theta}{\tau}\,
\bigl(\,{\stackrel{*}{\mu}}{}_{\tau l}^{\,j\zeta}\, {\rm Re}\,C_{\theta-\tau}(t) -
\widetilde{\mu}_{\tau l}^{\,j\zeta}\, {\rm Im}\,C_{\theta-\tau}(t)\bigr)\biggr)\,dt_j\,,
\hfill
$
\\[3ex]
\mbox{}\hfill
$
{\rm Im}\,C_{\theta}\colon t\to 
{\displaystyle \int\sum\limits_{j=1}^m} 
\biggl(\beta_{{}_{\scriptstyle \theta}}(t, f_j(t)) +
{\displaystyle \sum\limits_{\tau=1}^{\theta}}\binom{\theta}{\tau}\,
\bigl(\,{\stackrel{*}{\mu}}{}_{\tau l}^{\,j\zeta}\, {\rm Im}\,C_{\theta-\tau}(t)  +
\widetilde{\mu}{}_{\tau l}^{\,j\zeta}\, {\rm Re}\,C_{\theta-\tau}(t)\bigr)\biggr)\,dt_j
\hfill
$
\\[2.75ex]
\mbox{}\hfill
for all $t\in \widetilde{\mathcal T},
\quad
\theta=0,\ldots, s_{l}-1,
\ \ \ \widetilde{\mathcal T}\subset {\mathcal T},
\hfill
$
\\[2.75ex]
the real numbers 
\vspace{0.35ex}
${\stackrel{*}{\mu}}{}_{0l}^{\,j\zeta}={\rm Re}\,\lambda^j_l,\ \
\widetilde{\mu}{}_{0l}^{\,j\zeta}={\rm Im}\,\lambda^j_l,\ \
{\stackrel{*}{\mu}}{}_{\theta l}^{\,j\zeta}={\rm Re}\,\mu_{\theta l}^{\,j\zeta},\ \
\widetilde{\mu}{}_{\theta l}^{\,j\zeta}={\rm Im}\,\mu_{\theta l}^{\,j\zeta},\ j=1,\ldots, m,$ the vectors
${\stackrel{*}{\nu}}{}^{\,\theta l}={\rm Re}\,\nu^{\,\theta l},\ 
\widetilde{\nu}{}^{\,\theta l}={\rm Im}\,\nu^{\,\theta l},\ \theta=0,\ldots, s_{l}-1,$ and  
\\[2.5ex]
\mbox{}\hfill                                     
$
\displaystyle
\alpha_{{}_{\scriptstyle \theta}}(t,x) =
\biggl( {\stackrel{*}{\nu}}{}^{\,\theta l}x\cdot
\cos\sum\limits_{j=1}^{m}\widetilde{\mu}{}_{0 l}^{\,j\zeta}\;\!t_j  +
\widetilde{\nu}{}^{\,\theta l}x\cdot\sin\sum\limits_{j=1}^{m}\widetilde{\mu}{}_{0l}^{\,j\zeta}\;\!t_j\biggr)\cdot
\exp\biggl({}-\sum\limits_{j=1}^{m}{\stackrel{*}{\mu}}{}_{0l}^{\,j\zeta}\;\!t_j\biggr),
\hfill                                     
$
\\[2.5ex]
\mbox{}\hfill                                     
$
\displaystyle
\beta_{{}_{\scriptstyle \theta}}(t,x)=
\biggl(\widetilde{\nu}{}^{\,\theta l}x\cdot
\cos\sum\limits_{j=1}^{m}\widetilde{\mu}{}_{0 l}^{\,j\zeta}(t)\;\!t_j  -
{\stackrel{*}{\nu}}{}^{\,\theta l}x\cdot\sin\sum\limits_{j=1}^{m}\widetilde{\mu}{}_{0l}^{\,j\zeta}\;\!t_j\biggr)\cdot
\exp\biggl({}-\sum\limits_{j=1}^{m}{\stackrel{*}{\mu}}{}_{0l}^{\,j\zeta}\;\!t_j\biggr)
\hfill                                     
$
\\[2.25ex]
\mbox{}\hfill                                     
for all 
$(t,x)\in \widetilde{\mathcal T}\times {\mathbb R}^n,
\quad \theta=0,\ldots, s_{l}-1.
\hfill
$
\\[3ex]
\indent
{\bf Remark 2.4}.
Suppose the system (2.39) satisfies (2.2) and (2.40). Then 1-forms 
\\[2ex]
\mbox{}\hfill
$
{\stackrel{*}{\omega}}{}_{\tau-1}^{\theta}\colon t\to 
{\displaystyle \sum\limits_{j=1}^m}
\biggl(\binom{\theta}{\tau-1}\,\mu_{\theta-\tau+1,l}^{\,j\zeta}+
{\displaystyle \sum\limits_{\delta=1}^{\theta-\tau}}
\binom{\theta}{\delta}\,\mu_{\delta l}^{j\zeta}\cdot K_{\tau-1}^{\theta-\delta}(t)\biggr)\,dt_j,
\ \ \ 
\tau=1,\ldots, \theta,
\ \ \theta=1,\ldots, s_{l}-1,
\hfill
$
\\[2ex]
\mbox{}\hfill
$
\widetilde{\omega}_{\theta}\colon t\to 
{\displaystyle \sum\limits_{j=1}^m} 
\biggl(\nu^{\,\theta l}f_j(t)\cdot \varphi(t)+
{\displaystyle \sum\limits_{\tau=1}^{\theta}}
\binom{\theta}{\tau}\,\mu_{\tau l}^{j\zeta}\cdot C_{\theta-\tau}(t)\biggr)\,dt_j$
\ for all $t\in \widetilde{\mathcal T},
\quad 
\theta=0,\ldots, s_{l}-1
\hfill
$
\\[2ex]
are exact in a simply connected domain $\widetilde{\mathcal T}\subset {\mathcal T}.$

\newpage

\mbox{}
\\
\centerline{
\large\bf  
3.  First integrals of  linear systems of ordinary differential equations
}
\\[1.75ex]
\indent
{\bf  
3.1. Linear homogeneous systems of ordinary differential equations}
\\[1.25ex]
\indent
Let us consider a linear autonomous homogeneous system of  ordinary differential equations
\\[2ex]
\mbox{}\hfill                                          
$
\dfrac{dx}{dt} =  Ax,
$
\hfill (3.1)
\\[2ex]
where $x = \mbox{colon}(x_{1},\ldots,x_{n})\in{\mathbb R}^n,$ the $A = \bigl\|a_{ij}\bigr\|$ is a real $n\times n$ matrix.
\vspace{0.5ex}

Let B be the transpose of the  matrix $A.$
\vspace{0.5ex}

Using methods of Section 2 , we obtain the following statements. 
\vspace{0.75ex}

{\bf Theorem 3.1}.
{\it
Let $\nu\in{\mathbb C}^n$ be an eigenvector of the matrix $B.$ Then the linear function
$p\colon x\to \nu x$ for all $x\in {\mathbb R}^n$ is a partial integral of the system} (3.1).
\\[2ex]
\indent
{\bf 3.1.1. Autonomous first integrals}\footnote[4]{
This Subsubsection has been published in 
{\it Vestnik of the Yanka Kupala Grodno State Univ.}, 
2002, Ser. 2, 
No.2(11), 23 - 29; 2003, Ser. 2, No. 2(22), 50 - 60.
Using method of Jordan canonical form, the same result was obtain by M.~Falconi and J. Llibre in
{\it Qualitative theory of dynamical systems}, 2004, No. 4, 233 - 254.
} 
\vspace{1ex}

{\bf Theorem 3.2}.
{\it 
Let $\nu^{1},\ \nu^{2}$ be real eigenvectors of  the matrix $B$ 
corresponding to the eigenvalues $\lambda_{1}, \ \lambda_{2} \ (\lambda_{1}\ne\lambda_{2}),$ respectively. 
Then the system {\rm (3.1)} has the first integral
\\[1.5ex]
\mbox{}\hfill                                 
$
F\colon x\to \bigl|\nu^{1} x\bigr|^{h_1}\,
\bigl|\nu^{2}x{\bigr|}^{h_2}$
\ \ for all $x\in {\mathscr  X},
\quad {\mathscr  X}\subset {\rm D}(F),
\hfill 
$
\\[1.5ex]
where $h_1,\ h_2$ is a real solution to the equation $\lambda_1h_1+\lambda_2h_2=0$ 
with $|h_1|+|h_2|\ne 0.$
}
\vspace{1.25ex}

{\bf Corollary 3.1}.
{\it
If $\nu$ is a real eigenvector of the matrix $B$ corresponding to the eigenvalue $\lambda=0,$
then the linear function
$
F\colon x\to \nu x
$
\ for all 
$x\in {\mathbb R}^n
$
is an autonomous first integral of the system of ordinary differential equations {\rm(3.1)}.
}
\vspace{1.25ex}

{\bf Corollary 3.2}.
{\it
Let $\lambda\ne 0$ be an eigenvalue of the matrix $B$ corresponding to two 
real linearly independent eigenvectors $\nu^{1}, \ \nu^{2}.$ Then 
\vspace{0.5ex}
the system {\rm (3.1)} has the autonomous first integral
$
F\colon x\to \dfrac{\nu^1 x}{\nu^2 x}$
\ \ for all 
$
x\in {\mathscr X},
$
where a domain ${\mathscr X}\subset \bigl\{x\colon \nu^2 x\ne 0\bigr\}.$
}
\vspace{1.75ex}

For example, the autonomous system of ordinary differential equations 
\\[2ex]
\mbox{}\hfill                              
$
\begin{array}{ll}
\dfrac{dx_1}{dt}=   x_{1} - 2x_{2} -  x_{4},
\quad
&
\dfrac{dx_2}{dt} = {}- x_{1} + 4x_{2} - x_{3} + 2x_{4},
\\[4ex]
\dfrac{dx_3}{dt} =  2x_{2} +  x_{3} + x_{4}, &
\dfrac{dx_4}{dt} =  2x_{1} - 4x_{2} + 2x_{3} - 2x_{4}
\end{array}
\hfill (3.2)
$
\\[2.25ex]
has the eigenvectors 
\vspace{0.5ex}
$\nu^{1}\!=\!(1,-1,1,-1),
\nu^{2}\!=\!(2,2,1,1),
\nu^{3}\!=\!(1,0,1,0),\!$
$\nu^{4}=(0,2,0,1)$
corresponding to the eigenvalues  
\vspace{0.5ex}
$\lambda_1=0,\ \lambda_2=\lambda_3=1, \ \lambda_4=2,$ respectively. 
The functions 
\\[2ex]
\mbox{}\hfill                              
$
F_{1}\colon x\to x_{1}-x_{2}+x_{3}-x_{4}$
\ \ for all $x\in {\mathbb R}^4
$
\ \ \  (by Corollary 3.1),
\hfill (3.3)
\\[2.75ex]
\mbox{}\hfill                              
$
F_{23}\colon x\to
\dfrac{2x_{1}+2x_{2}+x_{3}+x_{4}}{x_{1}+x_{3}}
$
\ \ for all 
$x\in {\mathscr X}_1
$
\ \ \ (by Corollary 3.2),
\hfill (3.4)
\\[2.5ex]
\mbox{}\hfill
$
F_{24}\colon x\to
\dfrac{(2x_{1}+2x_{2}+x_{3}+x_{4})^{2}}{2x_{2}+x_{4}}$
\ \ for all $x\in {\mathscr X}_2
$
\ \ \ (by Theorem 3.2),
\hfill (3.5)
\\[2ex]
where ${\mathscr X}_1\subset \{ x\colon x_{1}+x_{3}\ne 0\},\ 
{\mathscr X}_2\subset \{ x\colon 2x_{2}+x_{4}\ne 0\},$ are 
autonomous first integrals of the system (3.2).
The set of functionally independent first integrals $F_1,\ F_{23},\ F_{24}$ 
is a general autonomous integral of the system of ordinary differential equations (3.2).
\vspace{1ex}

{\bf  Theorem 3.3}.
{\it 
Let $\nu={\stackrel{*}{\nu}}+\widetilde{\nu}\,i\ 
({\stackrel{*}{\nu}}={\rm Re}\,\nu,\ \widetilde{\nu}={\rm Im}\,\nu)$
be an eigenvector of the matrix $B$ cor\-res\-pon\-ding to the complex eigenvalue
$\lambda={\stackrel{*}{\lambda}}+\widetilde{\lambda}\,i\ 
({\stackrel{*}{\lambda}}={\rm Re}\,\lambda,\ \widetilde{\lambda}={\rm Im}\,\lambda\not=0).$
Then the sys\-tem of ordinary differential equations {\rm (3.1)} has the autonomous first integral
\\[1.5ex]
\mbox{}\hfill                    
$
F\colon x\to
\bigl( ({\stackrel{*}{\nu}}x)^2
+(\widetilde{\nu}x)^2\bigr)\cdot
\exp\biggl(
{}-2 \
\dfrac{{\stackrel{*}{\lambda}}}{\widetilde{\lambda}} \
{\rm arctg}\,\dfrac{\widetilde{\nu}x}{{\stackrel{*}{\nu}}x}
\,\biggr)$
\ for all $x\in {\mathscr X},
\hfill 
$
\\[1.5ex]
where a domain ${\mathscr X}\subset \bigl\{x\colon {\stackrel {*}{\nu}}x\ne 0\bigr\}.$ 
}
\vspace{1ex}

{\bf Theorem 3.4}.
{\it 
Let $\nu^1\!=\!{\stackrel{*}{\nu}}{}^{\,1}+\widetilde{\nu}{}^{\,1}\,i\ 
({\stackrel{*}{\nu}}{}^{\,1}\!={\rm Re}\,\nu^1,\ \widetilde{\nu}{}^{\,1}\!={\rm Im}\,\nu^1)$
be an eigenvector of the matrix $B$ cor\-res\-pon\-ding to the complex eigenvalue
$\lambda_1={\stackrel{*}{\lambda}}_{1}+\widetilde{\lambda}_{1}\,i\ 
({\stackrel{*}{\lambda}}_{1}={\rm Re}\,\lambda_1,\ \widetilde{\lambda}_{1}={\rm Im}\,\lambda_1\not=0),$
$\nu^2$ be an real eigenvector of the matrix $B$ cor\-res\-pon\-ding to the eigenvalue $\lambda_2\ne 0.$
Then the sys\-tem of ordinary differential equations {\rm (3.1)} has the autonomous first integral
\\[2ex]
\mbox{}\hfill
$
F\colon x\to
\nu^2x\cdot \exp\Biggl(
{}-\dfrac{\lambda_2}{\widetilde{\lambda}_1} \
{\rm arctg}\,\dfrac{\widetilde{\nu}{}^{\,1}x}
{{\stackrel{*}{\nu}}{}^{\,1}x}\,\Biggr)$
\ for all 
$x\in {\mathscr X},
\hfill
$
\\[2ex]
where a domain ${\mathscr X}\subset \bigl\{x\colon {\stackrel{*}{\nu}}{}^1x\ne 0\bigr\}.$
}
\vspace{1ex}

The autonomous linear system of ordinary differential equations 
\\[2ex]
\mbox{}\hfill                              
$
\dfrac{dx_1}{dt} =  2x_{1} + x_{2},
\quad
\dfrac{dx_2}{dt} = x_{1} + 3x_{2} -  x_{3},
\quad
\dfrac{dx_3}{dt} =  {}-x_{1} + 2x_{2} + 3x_{3}
$
\hfill (3.6)
\\[2ex]
has the eigenvalues $\lambda_{1}=3+i,\ \lambda_{2}=2$ corresponding to two eigenvectors
$\nu^{1}=(1,i,{}-1),$ $\nu^{2}=(3,{}-1,{}-1),$ respectively, 
and the first autonomous integrals
\\[1.75ex]
\mbox{}\hfill                          
$
F_{1}\colon x\to
\bigl( (x_{1}-x_{3})^2+x_{2}^2\,\bigr)
\exp\Bigl( {}-6\,{\rm arctg}\,\dfrac{x_2}{x_1-x_3}\,\Bigr)
$
for all $x\in {\mathscr X}$ \ 
(by Theorem 3.3),
\hfill (3.7)
\\[2.75ex]
\mbox{}\hfill                           
$
F_{2}\colon x\to
(3x_1-x_2-x_3)\exp\Bigl( {}-2\,{\rm arctg}\,\dfrac{x_2}{x_1-x_3}\Bigr)
$
for all $x\in {\mathscr X}$ \ 
\ (by Theorem 3.4),
\hfill (3.8)
\\[1.75ex]
where a domain ${\mathscr X}\subset \{x\colon x_1-x_3\ne 0\}.$ 
\vspace{1.25ex}

{\bf Theorem 3.5}.
{\it
Let $\nu^1={\stackrel{*}{\nu}}{}^{\,1}+\widetilde{\nu}{}^{\,1}\,i$
and  $\nu^2={\stackrel{*}{\nu}}{}^{\,2}+\widetilde{\nu}\,{}^{\,2}\,i$
be two eigenvectors of the matrix $B$ corresponding to the complex eigenvalues
$\lambda_1=
{\stackrel{*}{\lambda}}_1 +\widetilde{\lambda}_1\,i$  and
$\lambda_2={\stackrel{*}{\lambda}}_2 +
\widetilde{\lambda}_2\,i\ (\lambda_1\!\ne \overline{\lambda}_2),$ respectively.  
Then the sys\-tem {\rm (3.1)} has the autonomous first integral
\\[2ex]
\mbox{}\hfill
$
F\colon x\to\
{\stackrel{\sim}{\lambda}}_{1}\,
{\rm arctg}\,\dfrac{{\stackrel{\sim}{\nu}}{}^{\,2}x}
{{\stackrel{*}{\nu}}{}^{\,2}x}
\ -\ {\stackrel{\sim}{\lambda}}_{2}\,
{\rm arctg}\,\dfrac{{\stackrel{\sim}{\nu}}{}^{\,1}x}
{{\stackrel{*}{\nu}}{}^{\,1}x}$
\ \ \ for all 
$x\in {\mathscr X},
\hfill
$
\\[2ex]
where the vectors ${\stackrel{*}{\nu}}{}^{\,\tau}\!={\rm Re}\,\nu^\tau,\ \widetilde{\nu}{}^{\,\tau}\!={\rm Im}\,\nu^\tau,$
the numbers ${\stackrel{*}{\lambda}}_{\tau}\!=\!{\rm Re}\,\lambda_\tau,\, 
\widetilde{\lambda}_{\tau}\!=\!{\rm Im}\,\lambda_\tau\!\not=\!0,\, \tau\!=\!1,2,$
a domain ${\mathscr X}\subset \bigl\{x\colon  {\stackrel{*}{\nu}}{}^{\,2}x\ne 0\wedge
{\stackrel{*}{\nu}}{}^{\,1}x\ne 0\bigr\}.$
}
\vspace{1.25ex}

As an example, the linear autonomous system of ordinary differential equations 
\\[1.75ex]
\mbox{}\hfill                              
$
\begin{array}{ll}
\dfrac{dx_1}{dt} = {}-3x_1 + x_2 + 4x_3+ 2x_4,
\quad
&
\dfrac{dx_2}{dt} =  8x_1 - 3x_2 -2x_3+ 6x_4,
\\[3ex]
\dfrac{dx_3}{dt} = {}-9x_1 + 3x_2 + 4x_3- 4x_4,
&
\dfrac{dx_4}{dt} = 6x_1 - 3x_2 - 4x_3+ 2x_4
\end{array}
\hfill (3.9)
$
\\[1.75ex]
has the eigenvalues $\lambda_1\!=i,\, \lambda_2\!=2i$ corresponding to the eigenvectors
$\!\nu^{1}\!\!=\!(1-i,-1+2i,2i,2),$ $\nu^{2}=(i,{}-1,i,1+2i),$ respectively.  
The functionally independent first integrals 
\\[1.75ex]
\mbox{}\hfill                              
$
F_{1}\colon x\to
(x_{1}-x_{2}+2x_{4})^2+({}-x_{1}+2x_{2}+2x_{3})^2$
\ for all $x\in {\mathbb R}^4$
\ (by Theorem 3.3),
\hfill (3.10)
\\[2.5ex]
\mbox{}\hfill                              
$
F_{2}\colon x\to ({}-x_2+x_4)^{2}+(x_1+x_3+2x_4)^2$
\ for all $x\in {\mathbb R}^4$
\ (by Theorem 3.3),
\hfill (3.11)
\\[1.75ex]
and (by Theorem 3.5)
\\[1.5ex]
\mbox{}\hfill                              
$
F_{3}\colon x\to\
{\rm arctg}\,\dfrac{x_{1}+x_{3}+2x_{4}}{{}-x_{2}+x_{4}}
-2\,{\rm arctg}\,\dfrac{{}-x_{1}+2x_{2}+2x_{3}}{x_{1}-x_{2}+2x_{4}}$
\ \ for all 
$x\in {\mathscr X}
\hfill (3.12)
$ 
\\[1.5ex]
are a general autonomous integral on a domain 
${\mathscr X}\subset \{x\colon x_1-x_2+2x_4\ne 0\wedge x_2-x_4\ne 0\}$
of the linear system of ordinary differential equations (3.9). 
\\[1.5ex]
\indent
{\bf Definition 3.1}.
{\it
Let $\nu^{0}$ be an eigenvector of the matrix $B$ corresponding to the eigen\-va\-lue $\lambda$ with 
elementary divisor of multiplicity $m.$ 
A non-zero vector $\nu^{k}\in {\mathbb C}^n$ is called a
\textit{\textbf{generalized eigenvector of order}} {\boldmath $k$} for $\lambda$ if and only if
\\[1.5ex]
\mbox{}\hfill
$
(B-\lambda E)\,\nu^{k}=k \cdot \nu^{k-1},
\quad 
k=1,\ldots, m-1,
\hfill
$
\\[1.5ex]
where $E$ is the $n\times n$ identity matrix.
}
\vspace{1ex}

{\bf Theorem 3.6}.
{\it
Let $\lambda$ be an eigenvalue of the matrix $B$ with the elementary divisor of multiplicity $m\ (m\geq 2)$
corresponding to the real eigenvector $\nu^{0}$ and to the real order {\rm 1} generalized eigenvector $\nu^{1}.$ 
Then the sys\-tem {\rm (3.1)} has the autonomous first integral
\\[1.5ex]
\mbox{}\hfill
$
F\colon x\to
\nu^{0}x\,
\exp\biggl( {}-\lambda\, \dfrac{\nu^{1}x}{\nu^{0}x}\,\biggr)
$
\ \ for all $x\in {\mathscr X},
\hfill
$
\\[1.5ex]
where a domain ${\mathscr X}\subset \{x\colon  \nu^0x\ne 0\}.$ 
}
\vspace{1ex}

{\bf Corollary 3.3}.
{\it 
Let $\lambda={\stackrel{*}{\lambda}}+\widetilde{\lambda}\,i\ 
({\stackrel{*}{\lambda}}={\rm Re}\,\lambda,\ \widetilde{\lambda}={\rm Im}\,\lambda\not=0)$
be a complex eigenvalue of the matrix $B$ with the elementary divisor of multiplicity $m\ (m\geq 2)$
corresponding to the eigenvector 
$\nu^{0}={\stackrel{*}{\nu}}{}^{\,0}+\widetilde{\nu}{}^{\,0}\,i$
and to the order {\rm 1} generalized eigenvector 
$\nu^{1}={\stackrel{*}{\nu}}{}^{\,1}+\widetilde{\nu}{}^{\,1}\,i.$
Then the sys\-tem {\rm (3.1)} has the autonomous first integrals
\\[1.5ex]
\mbox{}\hfill
$
F_1\colon x\to\,
\Bigr(\bigl({\stackrel{*}{\nu}}\,{}^{0}x\bigr)^2+
\bigl({\stackrel{\sim}{\nu}}\,{}^{0}x\bigr)^2\,\Bigr)
\exp\biggl( {}-2\,\dfrac{{\stackrel{*}{\lambda}}\,\alpha(x)-
{\stackrel{\sim}{\lambda}}\,\beta(x)}
{\bigl({\stackrel{*}{\nu}}\,{}^{0}x\bigr)^2+
\bigl({\stackrel{\sim}{\nu}}\,{}^{0}x\bigr)^2}\,\biggr)
$ 
\ \ for all  
$x\in {\mathscr X}
\hfill
$
\\[1.5ex]
and
\\[1.5ex]
\mbox{}\hfill
$
F_2\colon x\to\ 
{\rm arctg}\,
\dfrac{{\stackrel{\sim}{\nu}}\,{}^{0}x}{{\stackrel{*}{\nu}}{}^{0}x}
\ - \
\dfrac{{\stackrel{\sim}{\lambda}}\,\alpha(x)+
{\stackrel{*}{\lambda}}\,\beta(x)}
{\bigl({\stackrel{*}{\nu}}\,{}^{0}x\bigr)^2+
\bigl(\,{\stackrel{\sim}{\nu}}\,{}^{0}x\bigr)^2}
$
\ \ for all $x\in {\mathscr X},
\hfill
$
\\[2ex]
where a domain ${\mathscr X}\subset \bigl\{x\colon {\stackrel{*}{\nu}}{}^{\,0}x\ne 0\bigr\},$
\vspace{0.5ex}
the vectors ${\stackrel{*}{\nu}}{}^{\,\tau}\!={\rm Re}\,\nu^\tau,\ \widetilde{\nu}{}^{\,\tau}\!={\rm Im}\,\nu^\tau,\ \tau=0,1,$
the po\-ly\-no\-mi\-als 
$
\alpha\colon x\to\,
{\stackrel{*}{\nu}}{}^{\,0}x\,{\stackrel{*}{\nu}}{}^{\,1}x +
{\stackrel{\sim}{\nu}}\,{}^{\,0}x\,{\stackrel{\sim}{\nu}}{}^{\,1}x,
\ \ \beta\colon x\to\,
{\stackrel{*}{\nu}}{}^{\,0}x\,{\stackrel{\sim}{\nu}}{}^{\,1}x -
{\stackrel{\sim}{\nu}}\,{}^{\,0}x\,{\stackrel{*}{\nu}}{}^{\,1} x$
\ for all $x\in {\mathbb R}^n.$
}
\vspace{1ex}

{\bf Corollary 3.4}.
{\it 
Let $\lambda_1=0$
be an eigenvalue  with the elementary divisor of multiplicity $m\ (m\geq 2)$ of the matrix $B$
corresponding to the real eigenvector $\nu^{0}$
and the real order {\rm 1} generalized eigenvector  $\nu^{1}.$
Let $\lambda_2$ be an eigenvalue of the matrix $B$ 
corresponding to the real eigenvector $\nu^{2}.$
Then the sys\-tem {\rm (3.1)} has the autonomous first integrals
\\[2ex]
\mbox{}\hfill
$
F\colon x\to
\nu^{2}x\,
\exp\biggl(-\,\lambda_{2}\ \dfrac{\nu^{1}x}{\nu^{0}x}\,\biggr)
$
\ \ \ for all $x\in {\mathscr X},
\hfill
$
\\[2ex]
where a domain ${\mathscr X}\subset \{x\colon \nu^0x\ne 0\}.$
}
\vspace{1.25ex}

The linear autonomous system of ordinary differential equations 
\\[2ex]
\mbox{}\hfill                              
$
\dfrac{dx_1}{dt} =4x_1 - 5x_2+ 2x_3,
\quad
\dfrac{dx_2}{dt} = 5x_1 - 7x_2+ 3x_3,
\quad
\dfrac{dx_3}{dt} = 6x_1 - 9x_2 + 4x_3
$
\hfill (3.13)
\\[2ex]
has the eigenvalue $\lambda_1=0$ with the elementary divisor $\lambda^2$ of multiplicity $2$ 
corresponding to the eigenvector $\nu^{0}=(1,{}-2,1)$ and the order {\rm 1} generalized eigenvector  $\nu^{1}=(0,{}-1,1),$
and the simple eigenvalue $\lambda_{3}=1$ with the elementary divisor $\lambda-1$ corresponding to the eigenvector 
$\nu^{3}=(3,{}-3,1).$ The functionally independent functions 
\\[2ex]
\mbox{}\hfill                     
$
F_1\colon x\to x_1-2x_2+x_3$
\ \ for all 
$x\in {\mathbb R}^3
$
\ \ \ \ (by Theorem 3.6)
\hfill (3.14)
\\[1.75ex]
and (by Corollary 3.4)
\\[1.5ex]
\mbox{}\hfill                     
$
F_{2}\colon x\to
(3x_1-3x_2+x_3)\exp\dfrac{x_2-x_3}{x_1-2x_2+x_3}$
\ \ \ for all $x\in {\mathscr X}, 
\hfill (3.15)
$
\\[2ex]
where ${\mathscr X}\subset\{x\colon x_1-2x_2+x_3\ne 0\},$ are autonomous first integrals of the system (3.13).
\vspace{1.25ex}

{\bf Corollary 3.5}.
{\it 
Let $\lambda_1=0$
be an eigenvalue  with elementary divisor of multiplicity $m$ $(m\geq 2)\!$ of the matrix $B$
corresponding to the real eigenvector $\nu^{0}$
and to the real order {\rm 1} generalized eigenvector  $\nu^{1}.$
Let $\lambda_2={\stackrel{*}{\lambda}}_2+\widetilde{\lambda}_2\,i\ 
({\stackrel{*}{\lambda}}_2={\rm Re}\,\lambda_2,\ \widetilde{\lambda}_2={\rm Im}\,\lambda_2\not=0)$
be a complex eigenvalue of the matrix $B$ 
corresponding to the eigenvector 
$\nu^{2}={\stackrel{*}{\nu}}{}^{\,2}+\widetilde{\nu}{}^{\,2}\,i.$
Then the sys\-tem of differential equations {\rm (3.1)} has two autonomous first integrals
\\[2ex]
\mbox{}\hfill
$
F_1\colon x\to
\Bigl( \bigl(\,{\stackrel{*}{\nu}}{}^{\,2}x\bigr)^{2} +
\bigl(\,{\stackrel{\sim}{\nu}}{}^{\,2}x\bigr)^2\,\Bigr)
\exp\biggl({}-2\,
{\stackrel{*}{\lambda}}_{2} \ \dfrac{\nu^{1}x}{\nu^{0}x}\,\biggr)$
\ \ for all $x\in {\mathscr X}
\hfill
$
\\[1.5ex]
and
\\[1.5ex]
\mbox{}\hfill
$
F_2\colon x\to \
{\rm arctg}\,
\dfrac{{\stackrel{\sim}{\nu}}{}^{\,2}x}{{\stackrel{*}{\nu}}{}^{\,2}x}
\ - \
{\stackrel{\sim}{\lambda}}_2 \ \dfrac{\nu^{1}x}{\nu^{0}x}$
\ \ \ for all $x\in {\mathscr X},
\hfill
$
\\[1.75ex]
where a domain ${\mathscr X}\subset \{x\colon \nu^0x\ne 0,\ {\stackrel{*}{\nu}}{}^{\,2}x\ne 0\}.$
}
\vspace{1.25ex}

{\bf Theorem 3.7}.
{\it 
Let $\lambda$ be an eigenvalue  with elementary divisor of multiplicity $m\geq 2$ of the matrix $B$
corresponding to the real eigenvector $\nu^{0}$
and to the real generalized eigenvectors  $\nu^{k},\ k=\overline{1,m-1}.$ 
Then the system {\rm (3.1)} has the first integrals
\\[2ex]
\mbox{}\hfill                                       
$
F_g\colon x\to \Psi_g(x)$
\ \ for all $x\in {\mathscr X},
\quad g=2,\ldots,m-1,
$
\hfill {\rm (3.16)}
\\[2ex]
where the functions $\Psi_{g}\colon {\mathscr X}\to {\mathbb R}$ is the solution to system
\\[2ex]
\mbox{}\hfill
$
\displaystyle
\nu^{k}x =
\sum\limits_{i=1}^{k}{\textstyle\binom{k-1}{i-1}}
\Psi_{i}(x)\nu^{k-i}x,
\ \ \ k=1,\ldots,m-1,
\quad
{\mathscr X}\subset \{x\colon \nu^0x\ne 0\}.
\hfill
$
}
\\[2.25ex]
\indent
{\bf Remark 3.1}. 
Let $\lambda={\stackrel{*}{\lambda}}+\widetilde{\lambda}\,i\ 
({\stackrel{*}{\lambda}}={\rm Re}\,\lambda,\ \widetilde{\lambda}={\rm Im}\,\lambda\not=0)$ 
be a complex eigenvalue  with ele\-men\-ta\-ry divisor of multiplicity $m\ (m\geq 2)$ of the matrix $B$
corresponding to the complex eigen\-vec\-tor $\nu^{0}={\stackrel{*}{\nu}}{}^{\,0}+\widetilde{\nu}{}^{\,0}\,i.$
Taking into account (3.16), we obtain the real-valued first integrals of the system of ordinary 
differential equations (3.1): 
\\[2ex]
\mbox{}\hfill
$
F_{g,1}\colon x\to \mbox{Re}\,\Psi_g(x),
\ \ \ 
F_{g,2}\colon x\to \mbox{Im}\,\Psi_g(x)$
\ for all $x\in {\mathscr X}, \ \ \  g=2,\ldots, m-1,
\hfill
$
\\[2ex]
where a domain ${\mathscr X}\subset 
\bigl\{x\colon (\,{\stackrel{*}{\nu}}{}^{\,0}x)^2 +
(\,{\widetilde{\nu}}{}^{\,0}x)^2\ne 0\bigr\}.$
\vspace{1.25ex}

For example, the system of ordinary differential equations 
\\[2ex]
\mbox{}\hfill                                   
$
\dfrac{dx_1}{dt} =   4x_1 - x_2,
\quad \ \ 
\dfrac{dx_2}{dt} =  3x_1 + x_2 - x_3,
\quad \ \
\dfrac{dx_3}{dt} =    x_1 + x_3
\hfill (3.17)
$
\\[2ex]
has the eigenvalue $\lambda_{1}\!=\!2$ with the elementary divisor $(\lambda-2)^3\!$ of multiplicity 3 
corresponding to the eigenvector $\!\nu^{0}\!=\!(1,-1,1)\!$
and the generalized eigenvectors  
$\!\nu^{1}\!=\!(1,0,-1),\, \nu^{2}\!=\!(0,0,2).$
First integrals of the system of differential equations (3.17) are the functions
\\[2.25ex]
\mbox{}\hfill                                
$
F_{1}\colon (x_1,x_2,x_3)\to
(x_1-x_2+x_3)\exp\Bigl({}-2\, \dfrac{x_1-x_3}{x_1-x_2+x_3}\,\Bigr)
$
\ \ \ (by Theorem 3.6)
\hfill (3.18)
\\[2ex]
and (by Theorem 3.7)
\\[2.5ex]
\mbox{}\hfill                                
$
F_{2}\colon  (x_1,x_2,x_3)\to
\dfrac{(x_1-x_3)^2-2x_3(x_1-x_2+x_3)}{(x_1-x_2+x_3)^2}$
\ \ for all $(x_1,x_2,x_3)\in {\mathscr X},
$
\hfill (3.19)
\\[1.5ex]
where a domain ${\mathscr X}\subset \{(x_1,x_2,x_3)\colon x_1-x_2+x_3\ne 0\}.$
\vspace{1.25ex}

The autonomous system of ordinary differential equations 
\\[1.75ex]
\mbox{}\hfill                                   
$
\begin{array}{ll}
\dfrac{dx_1}{dt} = x_1 - 2x_2+x_3-2x_6,
\ &
\dfrac{dx_2}{dt} = 3x_2-x_3-x_5+2x_6,
\\[2.5ex]
\dfrac{dx_3}{dt} = {}-x_1 +x_3+2x_4+2x_5,
\ &
\dfrac{dx_4}{dt} = {}-x_1+x_4+x_5+x_6,
\\[2.5ex]
\dfrac{dx_5}{dt} = x_1 +x_2+x_5,
\ &
\dfrac{dx_6}{dt} = x_1-x_2+x_3-x_4-x_6
\end{array}
\hfill (3.20)
$
\\[2ex]
has the complex eigenvalue $\lambda_{1}=1+i$ with the elementary divisor $(\lambda-1-i)^3$ 
corresponding to the eigen\-vec\-tor $\nu^{0}=(1,1,0,0,i,0)$
and to the generalized eigenvectors  $\nu^{1}=(0,1,0,i,i,1),$ 
\linebreak
$\nu^{2}=(0,1,i,0,i,0).$
A general autonomous integral of the system (3.20) is the  functions
\\[2.25ex]
\mbox{}\hfill                                
$
F_{1}\colon x\to P(x)\,\exp({}-2\varphi (x))$
\ \ for all $x\in {\mathscr X}
$
\ \ \ (by Theorem 3.3),
\hfill (3.21)
\\[3ex]
\mbox{}\hfill                                
$
F_{2}\colon x\to
P(x)\exp\biggl({}-2\,\dfrac{\alpha(x)-\beta(x)}{P(x)}\,\biggr)$
\ \ for all $x\in {\mathscr X}
$
\ \ \ (by Corollary 3.3),
\hfill (3.22)
\\[3ex]
\mbox{}\hfill                                
$
F_{3}\colon x\to
\varphi(x)- \dfrac{\alpha(x)+\beta(x)}{P(x)}$
\ \ for all $x\in {\mathscr X}
$
\ \ \ (by Corollary 3.3),
\hfill (3.23)
\\[3ex]
\mbox{}\hfill                                
$
F_{4}\colon x\to
\dfrac{\gamma(x)P(x)+\beta^2(x)-\alpha^2(x)}{P^2(x)}$
\ \ \  for all $x\in {\mathscr X}
$
\ (by Theorem 3.7),
\hfill (3.24)
\\[2.25ex]
and  
\\[2.25ex]
\mbox{}\hfill                                
$
F_{5}\colon x\to
\dfrac{\delta(x)P(x)-2\alpha(x)\beta(x)}{P^2(x)}$
\ \ \ for all $x\in {\mathscr X}$
\ \ (by Theorem 3.7),
\hfill (3.25)
\\[2ex]
where
\\[2ex]
\mbox{}\hfill
$
P\colon x\to (x_1+x_2)^2 + x_5^2, 
\quad
\alpha\colon x\to (x_1+x_2)(x_2+x_6) + x_5(x_4+x_5), 
\hfill
$
\\[2.25ex]
\mbox{}\hfill
$
\beta\colon x\to (x_1+x_2)(x_4+x_5) - x_5(x_2+x_6), \
\quad
\gamma\colon x\to x_2(x_1+x_2) + x_5(x_3+x_5),
\hfill
$
\\[2.25ex]
\mbox{}\hfill
$
\delta\colon x\to (x_1+x_2)(x_3+x_5) - x_2x_5$
\ \ \ for all $x\in {\mathbb R}^6,
\hfill
$
\\[2.25ex]
\mbox{}\hfill
$
\varphi\colon x\to\  {\rm arctg}\,\dfrac{x_5}{x_1+x_2}$
\ \ \ for all $x\in {\mathscr X},
\quad {\mathscr X}\subset \{x\colon x_1+x_2\ne 0\}.
\hfill
$
\\[3.75ex]
\indent
{\bf 3.1.2. Nonautonomous first integrals}
\\[1.5ex]
\indent
{\bf Theorem 3.8}.
{\it 
Let $\nu$ be a real eigenvector of  the matrix $B$ 
corresponding to the eigenvalue $\lambda.$  
Then the system {\rm (3.1)} has the first integral
\\[1.5ex]
\mbox{}\hfill
$
F\colon (t,x)\to\ \nu x\,\exp({}-\lambda\,t)$
\ \ for all $(t,x)\in {\mathbb R}^{n+1}.
\hfill
$
}
\\[2ex]
\indent
For example, the four dimensional system (3.2) has 
the eigenvalue $\lambda_2=1$ corresponding to the eigenvector $\nu^2=(2,2,1,1)$ and the first integral
\\[1.75ex]
\mbox{}\hfill
$
F\colon (t,x)\to (2x_1+2x_2+x_3+x_4)\,e^{{}-t}$
\ \ for all $(t,x)\in {\mathbb R}^5
$
\ \ (by Theorem 3.8).
\hfill\mbox{}
\\[1.75ex]
\indent
The first integrals (3.3), (3,4),  (3.5), and $F$ are a general integral on a domain ${\mathbb R}\times {\mathscr X}$
of the system (3.2),  where 
${\mathscr X}\subset \{x\colon x_1+x_3\ne 0\wedge 2x_2+x_4\ne 0\}\subset {\mathbb R}^4.$
\vspace{0.5ex}

Using the eigenvector $\nu^2=(3,{}-1,{}-1)$ corresponding to the eigenvalue $\lambda_2=2,$
we can build the first integral of the linear system (3.6):
\\[1.75ex]
\mbox{}\hfill
$
F\colon (t,x_1,x_2,x_3)\to (3x_1-x_2-x_3)\,e^{{}-2t}$
\ \ for all $(t,x_1,x_2,x_3)\in {\mathbb R}^4
$
\ \ (by Theorem 3.8).
\hfill\mbox{}
\\[1.75ex]
\indent
The functionally independent first integrals (3.7), (3,8), and $F$ are 
a general integral on a do\-main ${\mathbb R}\times {\mathscr X}$
of the system (3.6),  where a domain ${\mathscr X}\subset\{(x_1,x_2,x_3)\colon x_1-x_3\ne 0\}.$
\vspace{1.25ex}

{\bf  Corollary 3.6}.
{\it 
Let $\nu={\stackrel{*}{\nu}}+\widetilde{\nu}\,i\ 
({\stackrel{*}{\nu}}={\rm Re}\,\nu,\ \widetilde{\nu}={\rm Im}\,\nu)$
be an eigenvector of the matrix $B$ cor\-res\-pon\-ding to the complex eigenvalue
\vspace{0.35ex}
$\lambda={\stackrel{*}{\lambda}}+\widetilde{\lambda}\,i\ 
({\stackrel{*}{\lambda}}={\rm Re}\,\lambda,\ \widetilde{\lambda}={\rm Im}\,\lambda\not=0).\!$
Then the sys\-tem of ordinary differential equations {\rm (3.1)} has the first integrals
\\[1.75ex]
\mbox{}\hfill
$
F_1\colon (t,x)\to
\Bigl(\bigl(\,{\stackrel{*}{\nu}}x\bigr)^2
+ \bigl(\,{\stackrel{\sim}{\nu}}x\bigr)^2\,\Bigr)
\exp\bigl({}-2\,{\stackrel{*}{\lambda}}\,t\bigr)$
\ \ for all $(t,x)\in {\mathbb R}^{n+1}
\hfill
$
\\[1.5ex]
and
\\[1.5ex]
\mbox{}\hfill
$
F_2\colon (t,x)\to\
{\rm arctg}\,\dfrac{{\widetilde{\nu}}x}{{\stackrel{*}{\nu}}x}\ -\
{\widetilde{\lambda}}\,t$
\ \ for all $(t,x)\in {\mathbb R}\times {\mathscr X},
\hfill
$
\\[1.75ex]
where a domain ${\mathscr X}\subset \bigl\{x\colon {\stackrel{*}{\nu}}x\ne 0\bigr\}.$ 
}
\vspace{1.25ex}

As an example, the system of linear differential equations (3.9) has the eigenvalue 
$\lambda_1=i$ corresponding to the eigenvector $\nu^1=(1-i,{}-1+2i,2i,2)$ and the first integral
\\[1.75ex]
\mbox{}\hfill
$
F\colon (t,x)\to\
{\rm arctg}\,\dfrac{{}-x_1+2x_2+2x_3}{x_1-x_2+2x_4} \  -\ t$
\ \ for all $(t,x)\in {\mathbb R}\times {\mathscr X}_1
$
\ \ (by Corollary 3.6),
\hfill\mbox{}
\\[1.75ex]
where a domain ${\mathscr X}_1\subset \{x\colon x_1-x_2+2x_4\ne 0\}.$
\vspace{0.5ex}

The first integrals (3.10), (3.11), (3.12) and $F$ are 
\vspace{0.35ex}
a general integral on a domain ${\mathbb R}\times {\mathscr X}$
of the system (3.9),  where a domain 
${\mathscr X}\subset \{x\colon x_1-x_2+2x_4\ne 0\wedge x_4-x_2\ne 0\}\subset {\mathbb R}^4.$
\vspace{1.5ex}

{\bf Theorem 3.9}.
{\it
Let $\lambda$ be an eigenvalue of the matrix $B$ with the elementary divisor of multiplicity $m\ (m\geq 2)$
corresponding to the real eigenvector $\nu^{0}$ and the real order {\rm 1} generalized eigenvector $\nu^{1}.$ 
Then the sys\-tem {\rm (3.1)} has the first integral
\\[1.75ex]
\mbox{}\hfill                                
$
F\colon (t,x)\to \dfrac{\nu^1x}{\nu^0x} - t$
\ \ for all $(t,x)\in {\mathbb R}\times {\mathscr X},
\hfill (3.26)
$
\\[1.75ex]
where a domain ${\mathscr X}\subset \{x\colon \nu^0x\ne 0\}.$
}
\vspace{1.25ex}

The system of ordinary differential equations \!(3.13)\! has the eigenvalue $\!\lambda_1\!=\!0\!$  
corresponding to the eigenvector $\!\nu^0\!=\!(1,-2,1)\!$ and to the
order {\rm 1} generalized eigenvector $\nu^1\!=\!(0,-1,1).$ From Theorem 3.9 it follows that the function
\\[2ex]
\mbox{}\hfill
$
F\colon (t,x_1,x_2,x_3)\to \dfrac{x_3-x_2}{x_1-2x_2+x_3}\ - \ t$
\ \ for all $(t,x_1,x_2,x_3)\in {\mathbb R}\times {\mathscr X}
\hfill
$
\\[2ex]
is a first integral of the system (3.13), 
\vspace{0.5ex}
where a domain ${\mathscr X}\subset \{(x_1,x_2,x_3)\colon x_1-2x_2+x_3\ne 0\}.$

A general integral on a domain ${\mathbb R}\times {\mathscr X}$ of the system 
of ordinary differential equations (3.13) is 
the functionally independent first integrals (3.14), (3.15), and $F.$ 
\vspace{1.25ex}

Consider the system (3.17). 
Using the eigenvalue $\lambda_1=2$ corresponding to the eigenvector $\nu^0=(1,{}-1,1)$ and to the 
order {\rm 1} generalized eigenvector $\nu^1=(1,0,{}-1),$
we can build the first integral of the system of ordinary differential equations (3.17):
\\[2ex]
\mbox{}\hfill
$
F\colon (t,x_1,x_2,x_3)\to \dfrac{x_1-x_3}{x_1-x_2+x_3}\ - \ t$
\ \ for all $(t,x_1,x_2,x_3)\in {\mathbb R}\times {\mathscr X}
$
\ \ (by Theorem 3.9).
\hfill\mbox{}
\\[2.25ex]
\indent
The first integrals (3.18), (3.19), and $F$ are 
\vspace{0.35ex}
a general integral on a domain ${\mathbb R}\times {\mathscr X}$
of the system (3.17),  where a domain 
${\mathscr X}\subset \{(x_1,x_2,x_3)\colon x_1-x_2+x_3\ne 0\}\subset{\mathbb R}^3.$
\vspace{1.5ex}

{\bf Remark 3.2}. 
Suppose $\!\lambda\!=\!{\stackrel{*}{\lambda}}+\widetilde{\lambda}\,i\, 
({\stackrel{*}{\lambda}}\!=\!{\rm Re}\,\lambda,\, \widetilde{\lambda}\!=\!{\rm Im}\,\lambda\not\!=\!0)\!$ 
is a complex eigenvalue  of the matrix $B$ corresponding to the eigenvector 
$\nu^{0}={\stackrel{*}{\nu}}{}^{\,0}+\widetilde{\nu}{}^{\,0}\,i$ and 
to the generalised eigenvector $\nu^{1}={\stackrel{*}{\nu}}{}^{\,1}+\widetilde{\nu}{}^{\,1}\,i.$
Using (3.26), we get the real-valued first integrals of the system (3.1): 
\\[2ex]
\mbox{}\hfill
$\!
F_1\colon\! (t,x)\to
\dfrac{{\stackrel{*}{\nu}}\,{}^{0}x\,{\stackrel{*}{\nu}}{}^{1}x +
{\widetilde{\nu}}\,{}^{\,0}x\,{\widetilde{\nu}}{}^{\,1}x}
{\bigl(\,{\stackrel{*}{\nu}}\,{}^{0}x\bigr)^{2} +
\bigl(\,{\widetilde{\nu}}\,{}^{\,0}x\bigr)^{2}} \, - t,
\
F_2\colon\! (t,x)\to
\dfrac{{\stackrel{*}{\nu}}\,{}^{0}x\,{\widetilde{\nu}}{}^{\,1}x
-{\widetilde{\nu}}\,{}^{\,0}x\,{\stackrel{*}{\nu}}{}^{1}x}
{\bigl(\,{\stackrel{*}{\nu}}\,{}^{0}x\bigr)^{2} +
\bigl(\,{\widetilde{\nu}}\,{}^{\,0}x\bigr)^{2}}$
for all $(t,x)\!\in\! {\mathbb R}\!\times\!{\mathscr X},
\hfill
$
\\[1.75ex]
where a domain ${\mathscr X}\subset 
\bigl\{x\colon (\,{\stackrel{*}{\nu}}{}^{\,0}x)^2 +
(\,{\widetilde{\nu}}{}^{\,0}x)^2\ne 0\bigr\}.$
\vspace{1.25ex}

Consider the system (3.20). 
Using the eigenvalue $\lambda_1=1+i$ corresponding to the eigenvector $\nu^{0}=(1,1,0,0,i,0)$ and to the 
order {\rm 1} generalized eigenvector $\nu^{1}=(0,1,0,i,i,1),$
we can build the first integral of the system of ordinary differential equations (3.20):
\\[2ex]
\mbox{}\hfill
$
F\colon (t,x)\to
\dfrac{(x_1+x_2)(x_2+x_6)+x_5(x_4+x_5)}{(x_1+x_2)^2+x_5^2}\ - \ t$
\ \ for all $(t,x)\in {\mathbb R}\times {\mathscr X},
\hfill
$
\\[2ex]
where a domain ${\mathscr X}\subset \{x\colon x_1+x_2\ne 0\}\subset {\mathbb R}^6.$

The functionally independent first integrals (3.21), (3.22), (3.23), (3.24), (3.25), and $F$ 
are a general integral of the system of ordinary differential equations (3.20).
\\[3ex]
\indent
{\bf 3.2. Linear nonhomogeneous systems of ordinary differential equations}
\\[1ex]
\indent
We consider a system of $n$ first order constant-coefficient linear nonhomogeneous 
ordinary differential equations 
\\[1.25ex]
\mbox{}\hfill                                          
$
\dfrac{dx}{dt}=Ax+f(t),
$
\hfill (3.27)
\\[2.25ex]
where $x = \mbox{colon}(x_{1},\ldots,x_{n})\in{\mathbb R}^n,$ 
\vspace{0.5ex}
$A = \bigl\|a_{ij}\bigr\|$  is a real $n\times n$ matrix, 
and the vector function $f\colon t\to {\rm colon}(f_1(t),\ldots,f_n(t))$ for all $t\in J$ 
is continuous on an interval $J\subset {\mathbb R}.$
\vspace{0.5ex}

Let B be the transpose of the  matrix $A.$
\vspace{0.75ex}

{\bf Theorem 3.10}.
{\it
Suppose $\nu$ is a real eigenvector of the matrix $B$ corresponding to the eigenvalue $\lambda.$
Then a first integral of the system {\rm(3.27)} is the function
\\[2ex]
\mbox{}\hfill                                           
$
\displaystyle
F\colon (t,x)\to \nu x\cdot\exp({}-\lambda t)
-\int \nu f(t)\cdot\exp({}-\lambda t)\,dt$
\ \, for all $(t,x)\in \widetilde{J}\times {\mathbb R}^n, 
\ \ \widetilde{J}\subset J.
\hfill 
$
\\[2.5ex]
}
\indent
{\bf Corollary 3.7}.
{\it 
Let $\nu={\stackrel{*}{\nu}}+\widetilde{\nu}\,i\ 
({\stackrel{*}{\nu}}={\rm Re}\,\nu,\ \widetilde{\nu}={\rm Im}\,\nu)$
be an eigenvector of the matrix $B$ cor\-res\-pon\-ding to the complex eigenvalue
\vspace{0.35ex}
$\lambda={\stackrel{*}{\lambda}}+\widetilde{\lambda}\,i\ 
({\stackrel{*}{\lambda}}={\rm Re}\,\lambda,\ \widetilde{\lambda}={\rm Im}\,\lambda\not=0).$
Then first integrals of the sys\-tem of ordinary differential equations {\rm (3.27)} are the functions 
\\[2ex]
\mbox{}\hfill                                           
$
\displaystyle
F_{\tau}\colon (t,x)\to \alpha_{\tau}(t, x)-
\int \alpha_{\tau}(t, f(t))\,dt$
\ \ for all $(t,x)\in \widetilde{J}\times {\mathbb R}^n, 
\ \ \  \tau=1,2,\ \ \ \widetilde{J}\subset J,
\hfill 
$
\\[1.75ex]
where the functions 
\vspace{0.35ex}
$
\displaystyle
\alpha_1(t,x)=
\bigl(\,{\stackrel{*}{\nu}}x\cdot\cos\widetilde{\lambda}\,t +
\widetilde{\nu}x\cdot\sin\widetilde{\lambda}\,t \bigr)\cdot
\exp\bigl({}-{\stackrel{*}{\lambda}}\,t\bigr)$
for all $(t,x)\in J\times {\mathbb R}^n,
$ 
$
\displaystyle
\alpha_2(t,x)=
\bigl(\,\widetilde{\nu}x\cdot\cos\widetilde{\lambda}\,t -
{\stackrel{*}{\nu}}x\cdot\sin\widetilde{\lambda}\,t\bigr)\cdot
\exp\bigl({}-{\stackrel{*}{\lambda}}\,t\bigr)$
for all $(t,x)\in J\times {\mathbb R}^n.$
\\[2ex]
}
\indent
{\bf Remark 3.3}. Under the conditions of  Corollary 3.7, we have the function
\\[1.75ex]
\mbox{}\hfill                                     
$
\displaystyle
F\colon (t,x)\to 
\Bigr(\!({\stackrel{*}{\nu}}x)^2+(\widetilde{\nu}x)^2\Bigr)\!\exp\bigl(-2\,{\stackrel{*}{\lambda}}\,t\bigr) -
 2\biggl(\!\!\alpha_1(t,x)\!\int\!\! \alpha_1(t,f(t))dt+\alpha_2(t,x)\!\int\!\! \alpha_2(t,f(t))dt\!\!\biggr) +
\hfill                                     
$
\\[2.25ex]
\mbox{}\hfill                                     
$
\displaystyle
+\ \biggl(\int \alpha_1(t,f(t))\,dt\biggr)^2+
\biggl(\int \alpha_2(t,f(t))\,dt\biggr)^2$
\ \ \ for all $(t,x)\in \widetilde{J}\times {\mathbb R}^n
\hfill
$
\\[1.75ex]
is also a first integral of the system (3.27).

For example, the linear nonhomogeneous system of ordinary differential equations 
\\[1.75ex]
\mbox{}\hfill                             
$
\dfrac{dx_1}{dt} =  2x_{1} + x_{2}+2e^{2t},
\quad
\dfrac{dx_2}{dt} = x_{1} + 3x_{2} -  x_{3}+10,
\quad
\dfrac{dx_3}{dt} =  {}-x_{1} + 2x_{2} + 3x_{3}+e^{3t}
\hfill 
$
\\[1.75ex]
has the eigenvalue $\lambda_{1}=2$ corresponding to the eigenvector $\nu^{1}=(3,{}-1,{}-1),$ 
the eigenvalue $\lambda_{2}=3+i$ corresponding to the eigenvector $\nu^{2}=(1,i,{}-1),$ and the first integrals
\\[1.75ex]
\mbox{}\hfill                          
$
F_{1}\colon (t,x_1,x_2,x_3)\to
(3x_1-x_2-x_3-5)e^{{}-2t}+e^{t}-6t
$
\ \ \ (by Theorem 3.10),
\hfill\mbox{}
\\[2.25ex]
\mbox{}\hfill                         
$
F_{2}\colon (t,x_1,x_2,x_3)\to
\bigl( (x_{1}-x_{3}+1)\cos t+(x_{2}+3)\sin t\bigr) e^{{}-3t}+
(\cos t-\sin t)e^{{}-t}+\sin t,
\hfill
$
\\[2.5ex]
\mbox{}\hfill                         
$
F_{3}\colon (t,x_1,x_2,x_3)\to
\bigl( (x_{2}+3)\cos t+(x_{3}-x_{1}-1)\sin t\bigr) e^{{}-3t}-
(\cos t+\sin t)e^{{}-t}+\cos t
\hfill
$
\\[2ex]
\mbox{}\hfill                         
for all $(t,x_1,x_2,x_3)\in {\mathbb R}^4$ 
\ \ \ 
(by Corollary 3.7).
\hfill\mbox{}
\\[2.5ex]
\indent
{\bf Theorem 3.11}.
{\it 
Let $\lambda$ be an eigenvalue  with the elementary divisor of multiplicity $m$ $(m\geq 2)$ of the matrix $B$
corresponding to the real eigenvector $\nu^{0}$
and to the real generalized eigenvectors  $\nu^{k},\ k=1,\ldots, m-1.$ 
Then the system {\rm (3.27)} has the first integrals
\\[1.75ex]
\mbox{}\hfill                                  
$
\displaystyle
F_{k+1}\colon (t,x)\to 
(\nu^kx)\cdot\exp({}-\lambda t)-
\sum\limits_{i=0}^{k-1}
{\textstyle \binom{k}{i}}\, t^{k-i}\cdot F_{i+1}(t,x) \ -
\hfill                              
$
\\[2.75ex]
\mbox{}\hfill                                  
$
\displaystyle
-\ \int\biggl(\bigl(\nu^kf(t)\bigr)\cdot \exp({}-\lambda t)+
k\cdot \int\biggl(\bigl(\nu^{k-1}f(t)\bigr)\cdot \exp({}-\lambda t) \ +
\hfill                              
$
\\[2.75ex]
\mbox{}\hfill                                  
$
\displaystyle
+\, (k-1)\cdot\int\biggl(\bigl(\!\nu^{k-2}f(t)\!\bigr)\cdot \exp({}-\lambda t)+\ldots+
2\int\biggl(\bigl(\!\nu^{1}f(t)\!\bigr)\cdot \exp({}-\lambda t)\, +
\hfill                              
$
\\[2.5ex]
\mbox{}\hfill                                  
$
\displaystyle
+\ \int\bigl(\nu^{0}f(t)\bigr)\cdot \exp({}-\lambda t)\,dt\biggr)\,dt\ldots\biggr)dt\biggr)dt\biggr)dt$
\ for all $(t,x)\in \widetilde{J}\times {\mathbb R}^n,
\ \ \ k=1,\ldots,m-1,
\hfill                              
$
\\[2.25ex]
where the first integral {\rm(}by Theorem {\rm 3.10)}
\\[2ex]
\mbox{}\hfill                                                               
$
\displaystyle
F_1\colon\! (t,x)\!\to\! (\nu^0 x)\cdot\exp({}-\lambda t)
-\int\!\! \bigl(\nu^0 f(t)\bigr)\cdot\exp({}-\lambda t)\,dt$
\ for all $(t,x)\in \widetilde{J}\times {\mathbb R}^n,
\ \ \widetilde{J}\subset J.
\hfill 
$
\\[2.5ex]
}
\indent
As an example, the system of ordinary differential equations 
\\[2ex]
\mbox{}\hfill                                   
$
\dfrac{dx_1}{dt} =   4x_1 - x_2+e^{3t},
\quad
\dfrac{dx_2}{dt} =  3x_1 + x_2 - x_3+8t,
\quad
\dfrac{dx_3}{dt} =    x_1 + x_3+4
\hfill 
$
\\[2ex]
has the eigenvalue $\!\lambda_{1}\!=\!2\!$ with the elementary divisor $\!(\lambda-2)^3\!$ of multiplicity 3 
corresponding to the eigenvector $\!\nu^{0}\!=\!(1,-1,1)\!$
and the generalized eigenvectors  
$\!\nu^{1}\!=\!(1,0,-1), \nu^{2}\!=\!(-2,2,0).$
First integrals of this system of differential equations are the functions
\\[2ex]
\mbox{}\hfill                           
$
F_{1}\colon (t, x)\to
(x_1-x_2+x_3-4t)e^{{}-2t}-e^{t}
$
\ for all $(t, x)\in {\mathbb R}^4$
\ \ \ (by Theorem 3.10),
\hfill\mbox{}
\\[2.5ex]
\mbox{}\hfill                           
$
F_{2}\colon (t, x)\to
(x_1-x_3+2t-1) e^{{}-2t}-t\, F_1(t,x)-2e^{t}
$
 for all $(t, x)\in {\mathbb R}^4$
\ (by Theorem 3.11),
\hfill\mbox{}
\\[2.25ex]
and (by Theorem 3.11)
\\[2ex]
\mbox{}\hfill                           
$
F_{3}\colon (t, x)\to
2(x_2-x_1+3t+2) e^{{}-2t}-t^2\,F_1(t,x)-2t\,F_2(t,x)-2e^{t}
$
for all $(t, x)\in {\mathbb R}^4.$
\hfill\mbox{}
\\[2.25ex]
\indent
{\bf Theorem 3.12}.
{\it 
Let $\lambda={\stackrel{*}{\lambda}}+\widetilde{\lambda}\,i\ (\widetilde{\lambda}\ne 0)$ 
be an eigenvalue  with the elementary divisor of multiplicity $m\ (m\geq 2)$ of the matrix $B$
corresponding to the eigenvector $\nu^0={\stackrel{*}{\nu}}\;\!{}^{0}+\widetilde{\nu}\;\!{}^{0}\,i$
and to the generalized eigenvectors  
$\nu^k={\stackrel{*}{\nu}}\;\!{}^{k}+\widetilde{\nu}\;\!{}^{k}\,i,\, k\!=\!1,\ldots,m\!-\!1.\!$ 
Then the system of  ordinary differential equations {\rm (3.27)} has the first integrals
\\[2ex]
\mbox{}\hfill                                  
$
\displaystyle
F_{{}_{\scriptstyle \tau,k+1}}\colon (t,x)\to 
\alpha_{{}_{\scriptstyle \tau k}}(t,x)-
\sum\limits_{\xi=0}^{k-1}
{\textstyle \binom{k}{\xi}}\, t^{k-\xi}\cdot F_{{}_{\scriptstyle \tau, \xi+1}}(t,x) \ -
\hfill                              
$
\\[3ex]
\mbox{}\hfill                                  
$
\displaystyle
-\ \int\biggl(\alpha_{{}_{\scriptstyle \tau k}}(t, f(t)) +
k\cdot \int\biggl(\alpha_{{}_{\scriptstyle \tau,k-1}}(t, f(t))  +
(k-1)\cdot\int\biggl(\alpha_{{}_{\scriptstyle \tau,k-2}}(t, f(t))\ + 
\hfill                              
$
\\[3.5ex]
\mbox{}\hfill                                  
$
\displaystyle
+\ \ldots+
2\int\biggl(\alpha_{{}_{\scriptstyle \tau 1}}(t, f(t))
+\int \alpha_{{}_{\scriptstyle \tau 0}}(t, f(t))\,dt\biggr)\,dt\ldots\biggr)dt\biggr)dt\biggr)dt
\hfill                              
$
\\[3ex]
\mbox{}\hfill                                  
for all $(t,x)\in \widetilde{J}\times {\mathbb R}^n,
\ \ \ k=1,\ldots, m-1,
\ \ \ \ \tau=1,2,
\quad
\widetilde{J}\subset J,
\hfill                              
$
\\[2.25ex]
where the functions 
\\[2ex]
\mbox{}\hfill                                     
$
\displaystyle
\alpha_{{}_{\scriptstyle 1k}}(t,x)=
\bigl(\,{\stackrel{*}{\nu}}\;\!{}^kx\cdot\cos\widetilde{\lambda}\,t +
\widetilde{\nu}\;\!{}^kx\cdot\sin\widetilde{\lambda}\,t \bigr)\cdot
\exp\bigl({}-{\stackrel{*}{\lambda}}\,t\bigr)$
\ for all $(t,x)\in {\mathbb R}^{n+1},
\ k=0,\ldots, m-1,
\hfill
$
\\[2.75ex]
\mbox{}\hfill                                     
$
\displaystyle
\alpha_{{}_{\scriptstyle 2k}}(t,x)=
\bigl(\,\widetilde{\nu}x\;\!{}^k\cdot\cos\widetilde{\lambda}\,t -
{\stackrel{*}{\nu}}\;\!{}^kx\cdot\sin\widetilde{\lambda}\,t\bigr)\cdot
\exp\bigl({}-{\stackrel{*}{\lambda}}\,t\bigr)$
\ for all $(t,x)\in {\mathbb R}^{n+1},
\ k=0,\ldots, m-1,
\hfill
$
\\[2.5ex]
and the first integrals {\rm(}by Corollary {\rm 3.7)}
\\[2.5ex]
\mbox{}\hfill                                          
$
\displaystyle
F_{\tau 1}\colon (t,x)\to 
\alpha_{{}_{\scriptstyle \tau 0}}(t, x)-\int \alpha_{{}_{\scriptstyle \tau 0}}(t, f(t))\,dt$
\ \ for all $(t,x)\in \widetilde{J}\times {\mathbb R}^n,
\quad \tau=1,2.
\hfill 
$
\\[3ex]
}


}

\begin{thebibliography}{99}

\bibitem{1}
B.V. Shabat, {\it An introduction to complex analysis}: 
{\it functions of several variables} (Russian), P. 2, Nauka, Moscow, 1985.    

\bibitem{2}
V.N. Gorbuzov and V.Yu. Tushchenco, 
On R-holomorphic solutions of a system of total differential equations (Russian), 
{\it Vestsi Nats. Akad. Navuk Belarusi}, Ser. Fiz.-Mat. Navuk, 1999, No. 3, 124-126.   

\bibitem{3}
V.N. Gorbuzov and V.Yu. Tushchenco, 
R-holomorphic solutions of a total differential equation, 
{\it Differential Equations} 35 (1999), No. 4, 447-452.   

\bibitem{4}
I.V. Gaishun, 
{\it Completely solvable Multidimensional differential equations} (Russian), 
Nauka and Texnika, Minsk, 1983.

\bibitem{5}
V.N. Gorbuzov, 
{\it Integrals of  differential systems} (Russian), Grodno State University, Grodno, 2006. 

\bibitem{6}
V.N. Gorbuzov, 
Autonomy of a system of equations in total differentials, 
{\it Differential Equations} 34 (1998), No. 2, 153-160.

\bibitem{7}
G. Darboux,
M\'{e}moire sur les \'{e}quations differentielles
algebriques du premier ordre et du premier degr\'{e},
{\it Bull. Sc. Math.},  Vol. 2 (1878), 60-96, 123-144, 151-200.

\bibitem{8}
V.N. Gorbuzov, 
{\it Integrals of systems of equations in total differentials} (Russian), 
Grodno State University, Grodno, 2005. 

\bibitem{9}
V.N. Gorbuzov, 
Particular integrals of real autonomous polynomial systems of exact differential equations (Russian),  
{\it J. Differential equations and control processes}, 2000, No. 2, 1-36 (http://www.neva.ru/journal).

\bibitem{10}
V.N. Gorbuzov and A.F. Pronevich, 
Spectral method of jacobian systems integral basic in partial equations construction (Russian),  
{\it J. Differential equations and control processes}, 2001, No. 3, 17-45 (http://www.neva.ru/journal).

\bibitem{11}
A.F. Pronevich, 
Autonomous integrals of linear total differential systems (Russian), 
Deponent VINITI of {\it Differentsial Uravneniya} 
02.10.2002, No. 1667-B2002, 24 p. 

\bibitem{12}
V.N. Gorbuzov and A.F. Pronevich, 
Building of  the integrals of  linear differential systems (Russian), 
{\it Vestnik of the Grodno State Univ.}, 2003, Ser. 2, No. 2(22), 50-60.   

\bibitem{13}
A.F. Pronevich, 
Integrals of linear multidimensional system of non-derogatory matrix structure (Russian), 
{\it Mathematical research}, Vol. 10 (2003), 143-152.

\bibitem{14}
V.N. Gorbuzov and A.F. Pronevich, 
Integrals of R-linear systems of exact differential (Russian), {\it Dokl. Akad. Nauk Belarusi} 48 (2004), No. 1, 49-52.   


\bibitem{15}
A.F. Pronevich, 
Integrals of partial differential systems with {\rm R}-linear coefficients (Russian), 
{\it Vestnik of the Yanka Kupala Grodno State Univ}, 2005, Ser. 2,  No. 1(31),  45-52.

\bibitem{16}
A.F. Pranevich, {\it {\rm R}-differentiate integrals of exact differential systems} (Russian), 
Dissert. Kand. fiz-mat. nauk (Ph.D Thesis), Grodno State Univ.,  Grodno, 2005. 

\bibitem{17}
A.F. Pranevich,
Building of first integrals for linear non-homogeneous system of ordinary differential equations 
with constant coefficients (Russian), 
{\it Vestnik of the Yanka Kupala Grodno State Univ.}, 2008, Ser. 2,  No. 2(68), 5-10. 


\bibitem{18}
C. Christopher and J. Llibre, Algebraic aspects of integrability for polynomial systems, 
{\it Qua\-li\-ta\-tive Theory of Dynamical Systems}, 1999, No. 1, 71-95.

\bibitem{19}
M.V. Dolov and S.A. Chistyakova,
On the structure of the general solution and integrating factor in a neighborhood 
of a simple singular point,
{\it Differential Equations}, 37 (2001),  No.5, 747-750.

\bibitem{20}
N.P. Erugin, 
{\it A book for reading the general course of differential equations} (Russian), 
Nauka and Technika, Minsk, 1972.

\bibitem{21}
M. Falconi and J. Llibre, $n-1$ independent first integrals for linear differential systems
in ${\mathbb R}^n$ and ${\mathbb C}^n,$
{\it Qualitative theory of dynamical systems}, 2004, No. 4, 233-254.

\bibitem{22}
S.D. Furta, On non-integrability of general systems of differential equations, 
{\it Z. Angew. Math. Phys.}, 1996, No. 47, 112-131.

\bibitem{23}
V.N. Gorbuzov, 
Autonomous integrals and Jacobi last factors for systems of ordinary 
differential equations, {\it Differential Equations} 30 (1994), No. 6, 868-875.

\bibitem{24}\!\!
V.N.Gorbuzov and S.N.Daranchuk, A basis of the autonomous first integrals of the Jaco\-bi-Fourier's system
(Russian), {\it Vestnik of the Belarusian State Univ.}, 2005, No.3, 70-74.

\bibitem{25}\!
V.N. Gorbuzov and P.B. Pauliuchyk,  
On the solutions, integrals and limit cycles of  
{\it n}-Dar\-boux's system (Russian),  
{\it J. Differential equations and control processes}, 2002, No. 2, 26-46 (http://www.neva.ru/journal).

\bibitem{26}
V. N. Gorbuzov and V. Yu. Tyshchenko,
Particular integrals of systems of ordinary differential equations, 
{\it Mat. Sb.} 183 (1992), No. 3, 76-94.

\bibitem{27}
A. Goriely,
Integrability, partial integrability and nonintegrability for systems of ordinary 
differential equations,
{\it J. Math. Phys.} 37 (1996), 1871-1893.

\bibitem{28}
A. Goriely,
A brief history of Kovalevskaya exponents and modern developments,
{\it Regular and Chaotic Dynamics}, Vol. 5 (2000), No. 1, 3-15.

\bibitem{29}
B. Grammaticos, J. Moulin Ollagnier, A. Ramani, J.-M. Strelcyn, and S. Wojciechowski,
Integrals of quadratic ordinary differential equations in ${\mathbb R}^3\colon\!\!$ 
The Lotka-Volterra system, {\it Phy\-sica A} 163 (1990), 683 -722.

\bibitem{30}
A.I. Jablonskii,
Algebraic integrals of system of differential equations,
{\it Differential Equations} 6 (1970), No. 11, 1326-1333.
 
\bibitem{31}
C.G.J. Jacobi, 
Theoria nova multiplicatoris systemati
aequationum differentialium vulgarium applicandi,
{\it J. f\"{u}r reine und angew. Math.} Vol. 27 (1844), 199-268; 
Vol. 29 (1845), 213-279, 333-376.

\bibitem{32}
J.P. Jouanolou, 
{\it Equations de Pfaff alg\'{e}briques}, Springer-Verlag, New York, 1979.

\bibitem{33}
B.M. Koialovich,
{\it Researches on the differential equations} $ydy-ydx=Rdx$ (Russian),
Sankt Peterburg, 1894.

\bibitem{34}
A.N. Korkine, 
Thoughts about multipliers of differential equations of first degree (Russian),
{\it Math. Sbornik} 24 (1904), 194-350 and 351-416.

\bibitem{35}
S. Kowalevski,
Sur le probleme de la rotation d'un corps solide autour d'un point fixe,
{\it Acta Mathematica} 12 (1889), 177-232.

\bibitem{36}
V.V. Kozlov and D.V. Treschchev,
Kovalevskya  numbers of generalized Toda chains,
{\it Math. Notes} 46 (1989), 840-848.

\bibitem{37}
M.N. Lagutinskii,
{\it Partial Algebraic Integrals} (Russian),
Kharkov, 1908.

\bibitem{38}\!\!
K.Ja. Latysheva, 
On the works of  V.P. Jermakov on the theory of differential equations (Rus\-sian),
{\it Istoriko-Matematicheskije Issledovanija} 9, Gostekhizdat, Moscow (1956), 
\linebreak
691-722.

\bibitem{39}
P.D. Lax, Integrals of nonlinear equations of evolution and solitary waves, 
{\it Commun. Pure Appl. Math.}, 1968, No. 221, 467-490.

\bibitem{40}
J. Liouville,
M\'{e}moire sur l'int\'{e}gration d'une classe
d'\'{e}quations diff\'{e}rentielles du second ordre en quantit\'{e}s finies explicites,
{\it J. math. pures et appl}, Vol. 4 (1839),  423-456.

\bibitem{41}
A.J. Maciejewski,
About algebraic integrability and non-integrability of ordinary differential equations,
{\it Chaos, Solitons, Fractals} 9 (1998), 51-65.


\bibitem{42}
F. Minding,
Beitr\"{a}ge zur Integration der Differentialgleichungen erster Ordnung,
{\it Mem. de l'Acad. des Sci. de St.-Petersbourg VII-me s\'{e}rie},
Vol. 5 (1862), No. 1,  1-95.

\bibitem{43}
V.I. Mironenko, 
{\it Linear dependence of functions along solution of differential equations} (Russian),
Belarusian State University, Minsk, 1981.

\bibitem{44}
D. Morduhai-Boltovskoi,
{\it On integration of linear differential equations in finite terms} (Russian),
Warsaw, 1910.

\bibitem{45}
J. Moulin Ollagnier and J.-M. Strelcyn, 
On first integrals of linear systems, Frobenius integrability theorem and linear 
representations of Lie algebras, 
in Bifurcations of Planar Vector Fields, Proceedings, Luminy 1989, J.-P. Francoise and R. Roussarie (Eds.),
{\it Lecture Notes in Mathematics} 1455, Springer-Verlag, Berlin, Heidelberg, New-York, 1991.


\bibitem{46}
A.Nowicki,
On the nonexistrence of rational first integrals for systems of linear differential equations,
{\it Linear Algebra Appl.} 235 (1996), 107-120.

\bibitem{47}
M.J. Prelle and M.F. Singer,
Elementary first integrals of differential equations,
{\it Trans. Amer. Math. Soc.} 279 (1983), 215-229.

\bibitem{48}
M.F. Singer,
Liouvillian first integrals of differential equations,
{\it Trans. Amer. Math. Soc.} 333 (1992), 673-688.

\bibitem{49}
D. Schlomiuk, Algebraic particular integrals, integrability and the problem of the center, 
{\it Trans. Amer. Math. Soc.}, 338 (1993), No. 2, 799-841.

\bibitem{50}
D. Schlomiuk and N. Vulpe, 
Geometry of quadratic differential systems in the neighbourhood of the line
at infinity, {\it Preprint of the Cornell University}, arXiv:math/0405026v1 [math.DS] 
 3 May 2004.

\bibitem{51}
A. Tsygvintsev,
On the existence of polynomial first integrals of quadratic homogeneous systems of 
ordinary differential equations,
{\it J. Phys. A: Math. Gen.} 34 (2001), 2185-2193.

\bibitem{52}
J. A. Weil,  {\it Constantes et polynomes de Darboux en algebre differentielle: applications aux systemes
differentiesl lineaires} (Ph.D Thesis), Ecole polytechnique, 1995.

\bibitem{53}
H. Yoshida,
Necessary conditions for the existence of algebraic first integrals I, II,
{\it Celest. Mech.} 31 (1983), 363-379 and 381-399.

\bibitem{54}
S.L. Ziglin,
Branching of solutions and nonexistence of first integrals in Hamiltonian
mechanics,
{\it Functional Anal. Appl.} 16 (1983), 181-189.




\bibitem{55}
D.V. Buslyuk,
Integrals and last multipliers of partial differential systems,
{\it Differential Equations} 35 (1999), No. 3, 421-422.

\bibitem{56}
N.M. Gjunter,  
{\it Integrirovanie uravnenij v chastnyh proizvodnyh {\rm 1}-go poryadka} (Russian),
GTTI, Moscow-Leningrad, 1934.

\bibitem{57}
V.N. Gorbuzov and S.N. Daranchuk,
The integrals and last multipliers of one class of partial differential systems (Russian),  
{\it J. Differential equations and control processes}, 2007, 
\linebreak
No. 4, 1-16 (http://www.neva.ru/journal).

\bibitem{58}
E. Kamke, 
{it Differentialgleichungen: L\'{o}sungsmethoden und 
L\'{o}sungen, II, Partielle 
Differentialgleichungen Erster Ordnung f\'{u}r eine gesuchte Funktion}, 
Akad. Verlagsgesellschaft Geest \& Portig, Leipzig, 1965.

\bibitem{59}
A.D. Polanin, V.F. Zaitsev, and A. Moussiaux,
{\it Handbook of first order partial differential equations},
Taylor \& Francis, London, 2002.

\bibitem{60}
V.V. Amel'kin, 
{\it Autonomous and linear many-dimensional differential equations} (Russian),
Universitetskae, Minsk, 1985.

\bibitem{61}
A.R. Forsyth,
{\it Theory of differential equations}, Vol. 1-6,
Cambridge University Press (1890-1906),
reprinted by Dover Public., New York, 1959.

\bibitem{62}
V.N. Gorbuzov and A.F. Pranevich, 
Autonomy and cylindricality of R-differentiable integrals for systems in total differentials (Russian), 
{\it J. Differential equations and control processes}, 2008, No. 1, 35-49 
(http://www.neva.ru/journal). 


\bibitem{63}
E. Goursat, {\it Cours d'analyse mathematique}, Vol. 2, Paris (4th. ed. 1924). English translation:
{\it A Course of Mathematical Analysis}, Vol II, Part Two: Differential Equations, Dover Public.,
New York, 1959.

\bibitem{64}
F.R. Gantmacher,
{\it The theory of matrices} (Russian), Nauka, Moscow, 1988.

\bibitem{65}
R.A. Horn and C.R. Johnson,
{\it Matrix analysis},
Cambridge Univ. Press, Cambridge, 1990.

\end{thebibliography}
\end{document}